\documentclass[11pt]{article}

\usepackage{geometry}
\RequirePackage[OT1]{fontenc}
\RequirePackage{amsthm,amsmath}
\RequirePackage[authoryear]{natbib}
\RequirePackage[colorlinks,citecolor=blue,urlcolor=blue]{hyperref}

\usepackage{amssymb}
\usepackage{enumitem}  
\usepackage[hang,flushmargin]{footmisc}
\usepackage{lipsum}
\usepackage{changepage}
\usepackage{colortbl}
\usepackage[dvipsnames]{xcolor}
\usepackage{bbm} 
\usepackage{graphicx,adjustbox}
\usepackage{multirow,makecell}
\usepackage[normalem]{ulem}
\usepackage[linesnumbered,ruled,vlined]{algorithm2e}

\newcommand{\biblist}{
\bibliographystyle{apalike}
\bibliography{weightsRAMP}
}
\makeatletter
\newcommand{\algorithmfootnote}[2][\footnotesize]{%
  \let\old@algocf@finish\@algocf@finish
  \def\@algocf@finish{\old@algocf@finish
    \leavevmode\rlap{\begin{minipage}{\linewidth}
    #1#2
    \end{minipage}}%
  }%
}
\makeatother

\newcommand{\Cov}{\mathrm{Cov}}
\newcommand{\Var}{\mathrm{Var}}

\SetKwRepeat{Do}{do}{while}%
\SetKwInOut{Output}{output}
\SetKwInput{Initialization}{Initialization}
\SetKwProg{Fn}{Function}{:}{}
\def\bm#1{\mbox{\boldmath $#1$}}
\def\thick#1{\hbox{\rlap{$#1$}\kern0.25pt\rlap{$#1$}\kern0.25pt$#1$}}
\def\sg#1{{\scriptstyle{#1}}} 

\def\AMSE{{\rm AMSE}}
\def\MA{{\rm MA}}
\newtheorem{theorem}{Theorem}

\newtheorem{lemma}{Lemma}
\newtheorem{corollary}{Corollary}
\usepackage{authblk}
\usepackage{color}

\definecolor{coquelicot}{rgb}{0.20, 0.12, 0.72}

\usepackage{soul}

\title{Detangling robustness in high dimensions: composite versus model-averaged estimation}

\author[1]{ {Jing} Zhou \thanks{Jing.Zhou@kuleuven.be}}
\author [1]{Gerda  {Claeskens}  \thanks{Gerda.Claeskens@kuleuven.be}}
\affil[1]{ORStat and Leuven Statistics Research Center,  Katholieke Universiteit Leuven
         }
\author[2]{ {Jelena} {Bradic} \thanks{Contact author, jbradic@ucsd.edu}}
\affil[2]{Department of Mathematics and Halicio\u{g}l{u} Data Science Institute, University of California San Diego }
\date{}

\begin{document}

 \maketitle

\begin{abstract}
Robust methods, though ubiquitous in practice, are yet to be fully understood in the context of regularized estimation and high dimensions. Even simple questions become challenging very quickly. For example, classical statistical theory identifies equivalence between model-averaged and composite quantile estimation. However, little to nothing is known about such equivalence between methods that encourage sparsity. This paper provides a toolbox to further study
robustness in these settings and focuses on prediction.
In particular, we study optimally weighted model-averaged as well as composite $l_1$-regularized estimation. Optimal weights are determined by minimizing the asymptotic mean squared error. This approach incorporates the effects of regularization, without the assumption of perfect selection, as is often used in practice. Such weights are then optimal for prediction quality. Through an extensive simulation study, we show that no single method systematically outperforms others. We find, however, that model-averaged and composite quantile estimators often outperform least-squares methods, even in the case of Gaussian model noise.
Real data application witnesses the method's practical use through the reconstruction of compressed audio signals.
\end{abstract}

\section{Introduction}

 We investigate the benefits of model-averaged as well as composite estimators in high-dimensional problems where the underlying goal is superior prediction quality. Robustness in data analysis with potentially more parameters than samples is a critical practical question and is of particular interest in constructing recoveries of compressed images and signals which should have high precision.

Model averaging, often used as a first tool to improve estimation quality, forms a weighted average of estimators and is here utilized for regularized sparsity-encouraging estimation in a high-dimensional regression setting.
Model averaging is also well-known in the Bayesian setting \citep{HoetingMadiganRafteryVolinsky99}, though we focus on its frequentist version in which a user determines the weights assigned to the separate estimators \citep[e.g., see][]{ClaeskensHjort08, HjortClaeskens03, YuanYang2005}.
Model averaging enjoys a wide application, see, for example, the recent overview paper for model averaging in ecology by \citet{Dormannetal2018} and for application to hydrology by \citet{HogeGuthkeNowak2019}. In econometrics, the terminology ``forecast combinations" appears \citep[e.g., in][]{ChengWangYang2015,BatesGranger1969}; whereas ``multimodel inference'' is another commonly used term for this procedure \citep{BurnhamAnderson2002}.

While the technique is quite thoroughly investigated for low-dimensional models, far fewer results have been obtained in high dimensions.
\citet{ando2014model} consider high-dimensional linear regression. By computing the marginal correlation between each covariate and the response and forming groups according to the obtained values, regularized estimation is avoided. The authors fit a fixed number of low-dimensional models by the least squares method and subsequently average them.
\citet{ZhaoZhouLi2016} extend this method to dependent data, while \citet{ando2017weight} extend this approach to generalized linear models, again by only fitting low-dimensional models, this time via maximum likelihood estimation. In these papers, the weights are obtained via cross-validation; see also \citet{Hansen2007} and \citet{Hansen2012} for similar weight finding approaches in low-dimensional models.

Our setting is different and is theoretically valid (see Theorem 1 below). We explicitly work with $l_1$-regularized estimators that are averaged, and we do not rely on the correct low-dimensional representation of the model. When designing the optimal weights, we explicitly take variable selection effects (of regularization itself) into account. Is the dependence among regularized estimators an impediment or a hidden benefit in obtaining robust predictions, i.e., predictions that do not change much when the data is changed a little?

A second approach to robustness is through composite estimation. While model averaging combines estimators after optimization of their respective loss functions, composite estimation weights the loss functions directly (before optimization).
For quantile regression in low dimensions, \citet[][Theorem 5.2]{Koenker2005quantile} stated the asymptotic equivalence of model-averaged and composite quantile regression estimators, provided each method uses its own, optimal set of weights that minimize the asymptotic variance. Hence, with optimal weights, there is no asymptotic preference between the two methods in low dimensions. For high-dimensional quantile regression, when one restricts the
attention to inference regarding the true nonzero part of the regression coefficient and ignores the variable selection effect, \citet{BloznelisClaeskensZhou2019} obtained the same equivalence for high-dimensional quantile regression using different types of regularizations (SCAD, lasso, adaptive lasso).

In practice, however, one works with an estimated coefficient vector for which one is not sure that the regularization has led to the correct selection. Therefore, incorporating imperfections of variable selection is especially important for achieving robustness. This is where our approach differs from \citet{BloznelisClaeskensZhou2019} or \citet{bradic2011penalized}, where an irrepresentable condition (needed for consistent model or asymptotically perfect selection)
has been used to specify weights and analyze robustness.

The approximate message passing (AMP) algorithm is crucial in our approach to take the variable selection into account when studying the estimators' asymptotic mean squared errors.
 The use of such algorithms has been investigated by \citet{donoho2009message} and \citet{bayati2011dynamics} for compressed sensing.
\citet{donoho2016high} explain the use of AMP algorithms for obtaining the variance of high-dimensional M-estimators for which $n/p\to\delta\in(1,\infty)$. Here, $n$ denotes the sample size and $p$ the number of regression coefficients. However, the robustness of sparsity encouraging AMP estimators is still largely unknown.

In this paper, we first extend the robust AMP (RAMP) of \citet{bradic2016robustness} to regularized composite estimation. Second, we construct estimators and develop new theory for the asymptotic mean squared error (AMSE) both for model-averaged and for composite estimators. Note that model-averaged AMSE required an extension of AMP theory for a challenging case of dependent estimates. Besides, we establish new Stein-type risk estimates of the AMSE in both cases.

The new estimates of the AMSE of the model-averaged and composite estimators enable a theoretically justified and data-driven optimal weight choice by minimizing the estimated AMSE (without relying on perfect variable selection).
The estimated AMSE provides more information regarding the estimators than merely considering which variables have been selected.

Organization of the paper.
First, in Section~\ref{sec:estimation}, we detail the model-averaged and composite estimators in a high-dimensional setup. Next, we explain the model-averaged robust message passing algorithm in Section~\ref{sec:RAMP}. The limiting behavior of the estimators in the algorithm is studied by state evolution parameters in Section~\ref{section:state_evolution}. We obtain the estimators' asymptotic mean squared error as well as an estimator of that quantity in Section~\ref{sec:AMSE}. We showcase the procedure for high-dimensional regularized quantile regression in Section~\ref{sec:examples} and present numerical results in Section~\ref{sec:numerical}. Section~\ref{sec:discussion} concludes. All proofs, together with the assumptions and some technical lemmas, are collected in the Appendix.

\section{Model-averaged and composite estimation} \label{sec:estimation}

We consider a high-dimensional linear model $Y = X {\beta} + \bm{\varepsilon}$ with $Y \in \mathbbm{R}^n$, the design matrix $X \in \mathbbm{R}^{n \times p}$ and the parameter vector ${\beta} \in \mathbbm{R}^{p}$. The $i$th row of $X$ is denoted $X_{i\cdot}$, $i=1,\ldots,n$, the $j$th column of $X$ is denoted by $X_{\cdot j}$, $j=1,\ldots,p$. We assume the components of $\bm{\varepsilon}$ to be independent and identically distributed with mean zero, cumulative distribution function $F_\varepsilon$ and probability density function $f_\varepsilon$.
We allow for a sparse high-dimensional setup. Denote by $s$ the $l_0$ norm of the parameter vector, $s=\|\beta \|_0$, which counts the number of nonzero components of the vector $\beta$. We assume that the ratios $n/p \to \delta \in (0,1)$ and $n/s \to a \in (1, \infty)$ when $p, n, s$ tend to $\infty$.

We consider two types of weighted estimation methods. First, model-averaged estimation where estimators from different models or estimation methods are weighted and summed to arrive at a final estimator, see \eqref{eq:MA}. Second, composite estimation where a weighted average of loss functions is minimized; see \eqref{eq:Composite}.

For model-averaged estimation of the parameter $\beta$, define for $k=1,\ldots,K$ the regularized estimators
\begin{eqnarray}
\label{eq:MA_single}
\widehat{{\beta}}_k(\lambda_k) &=& \arg\min_{{\beta}\in\mathbbm{R}^{p}}\left\{ \sum_{i=1}^n \rho_k(Y_i-X_{i\cdot}{\beta})+\lambda_k\|{\beta}\|_1 \right\},
\end{eqnarray}
where $\rho_1,\ldots,\rho_K$ are nonnegative convex loss functions and ${\bm\lambda}=(\lambda_1,\ldots,\lambda_K)^\top$ is a vector of possibly different nonnegative regularization parameters.
For a set of weights $w=(w_1,\ldots,w_K)^\top$, the model-averaged estimator is defined as
\begin{eqnarray}
\label{eq:MA}
\widehat{{\beta}}_{{\rm MA}}({\bm\lambda}) &=& \sum_{k=1}^K w_k \widehat{{\beta}}_k(\lambda_k).
\end{eqnarray}
Often one assumes that the weights $w_1,\ldots,w_K$ are all nonnegative and sum to 1, although this is not necessary for the computation of the estimator.

For composite estimation we consider again $K$ loss functions, though only with a single nonnegative regularization parameter $\lambda$, such that the regularized composite estimator is defined as
\begin{eqnarray}
\label{eq:Composite}
\widehat{{\beta}}_{{\rm C}}(\lambda) &=& \arg\min_{\sg{\beta}\in\mathbbm{R}^{p}}\left\{ \sum_{k=1}^K \sum_{i=1}^n w_k \rho_k(Y_i-X_{i\cdot}{\beta})+\lambda\|{\beta}\|_1 \right\}.
\end{eqnarray}
Computationally, composite estimation is harder than model-averaged estimation and requires that all weights are positive to ensure a nonnegative and convex weighted loss function, even when all $\rho_k$ are nonnegative and convex. Hence, for composite estimation it is required that the weight vector $w \in [0,1]^K$ such that $\sum_{k=1}^K w_k=1$.

As a worked-out scenario throughout the paper, we consider quantile loss functions $\rho_k(\cdot), k = 1, \ldots, K$ that are defined below.
For more information about quantile regression with i.i.d. errors, see \citet[Sec.~3.2.2]{Koenker2005quantile}. In this paper, we assume that the design matrix $X$ does not contain a column of ones; see assumption \ref{Afirst} in the Appendix. With $\tau \in (0,1)$, the $\tau$-quantile of the response $Y$ is obtained as $X\beta+F_\varepsilon^{-1}(\tau) = X\beta + u_\tau$.

\begin{figure}[!t]
\includegraphics[width=\textwidth]{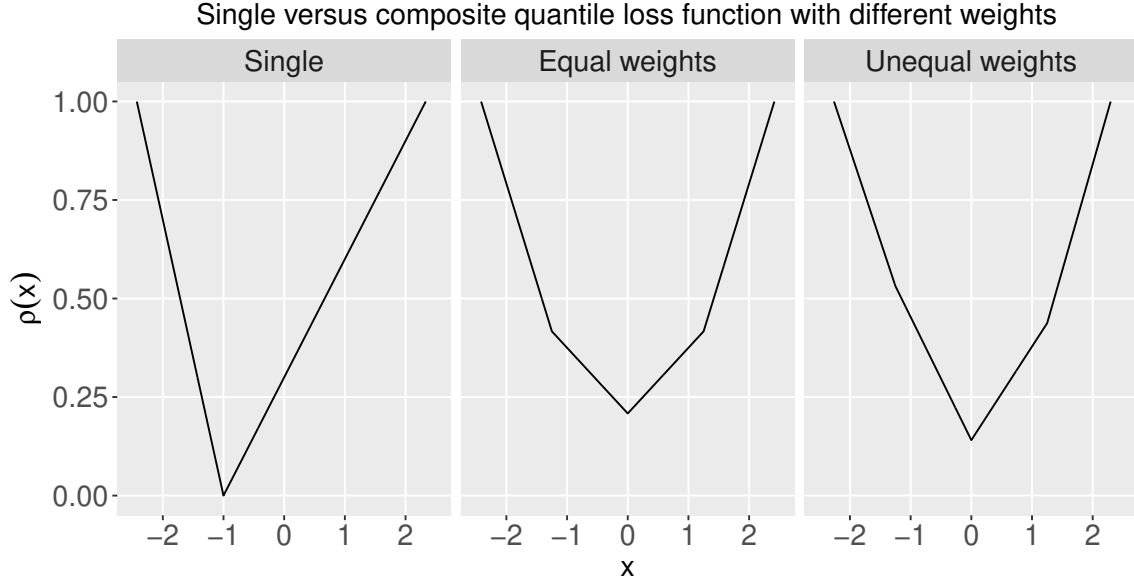}
\caption{Examples of quantile loss functions. Left: $\tau=0.3$ quantile loss function. Middle: Composite quantile loss function at quantile levels 0.25, 0.5, 0.75 with equal weights $w=(1/3,1/3,1/3)^\top$. Right: Composite quantile loss function at quantile levels 0.25, 0.5, 0.75 with weights $w=(0.15, 0.55, 0.3)^\top$.
\label{Fig:quantileloss}}
\end{figure}
Figure~\ref{Fig:quantileloss} presents first a single quantile loss function with $\tau=0.3$,
$$
\rho(x) = (x - u_{\tau}) (\tau - I\{x \leq u_{\tau}\}).
$$
For model averaging we specify $K$ different quantile levels and use $K$ different such quantile loss functions for estimation of $\beta$:
$$
\rho_k (x) = (x - u_{\tau_k}) (\tau_k - I\{x \leq u_{\tau_k}\}), \quad k \in \{1, \ldots, K\}.
$$
For composite quantile estimation we assume that the $K$ quantile levels $\tau_1 < \cdots < \tau_K$, then also the quantiles of $\varepsilon$ are sorted $u_{\tau_1} < \cdots <  u_{\tau_K}$. Define $u_{\tau_{0}}=-\infty$ and
$u_{\tau_{K+1}}=\infty$.

The middle panel of Figure~\ref{Fig:quantileloss} depicts such a composite quantile loss function $\rho_C=\sum_{k=1}^Kw_k\rho_k$ for $K=3$ quantile levels 0.25, 0.5 and 0.75 with equal weights $w=(1/3,1/3,1/3)^\top$. The panel on the right in Figure~\ref{Fig:quantileloss} uses the same quantile levels but depicts the quantile loss function $\rho_C$ with weights $w=(0.15, 0.55, 0.3)^\top$.

In general, the composite quantile loss function can be rewritten in the following way,
\begin{equation}
\label{eq:CQR-loss}
\rho_{\rm{C}}(x) = \left\{\begin{array}{ll}
  \sum_{k = 1}^K w_k(1 - \tau_k) (u_{\tau_k}-x),  & x < u_{\tau_1}\\
  \sum_{k = 1}^K w_k \tau_k (x - u_{\tau_k}),  & x \ge u_{\tau_K} \\
  \sum_{k =1 }^\ell w_k \tau_k |x - u_{\tau_k}| + \sum_{k = \ell+1}^K w_k(1 - \tau_k)|x - u_{\tau _k}|,
   & x \in [u_{\tau_\ell}, u_{\tau_{\ell+1}})\\
    \qquad \mbox{ for } \ell  = 1, \ldots, K-1.
\end{array}\right.
\end{equation}

Note that a single quantile loss function can be seen as a particular case of a composite loss function: take  $K=1$ and the single weight $w_1=1$. Theoretical results regarding regularized estimation for a single quantile loss function can be found in \citet{bradic2016robustness}.
Henceforth, we concentrate on the example of the composite case.

One aim of this paper it to investigate the weight choice $w$ by minimizing the asymptotic mean squared error of the estimators $\widehat{\beta}_{{\rm MA}}(\bm{\lambda})$ and $\widehat{\beta}_{{\rm C}}(\lambda)$.

\section{Robust approximate message passing} \label{sec:RAMP}

The idea behind approximate message passing algorithms is to provide an iterative procedure that has as its fixed point the estimator of interest; in this case the minimizer (\ref{eq:MA}) of the regularized loss function in the case of model averaging, and the estimator (\ref{eq:Composite}) in the case of composite estimation.
Due to a convergence in the mean square between the solution of the approximate message passing algorithm and the estimator (\ref{eq:MA}), respectively (\ref{eq:Composite}), the asymptotic mean squared error that holds for the solution of the approximate message-passing algorithm, is also the asymptotic MSE of the other estimator. Studying effects of regularization while allowing $n/p \to \delta \in (0,1)$ is challenging. The AMP provides theoretical advantages in these cases as it enables a complete and tractable, albeit challenging, structure for obtaining AMSE.
This paper is the first to obtain and use the asymptotic mean square error of the regularized estimators to optimize the weight choice of both the model-averaged estimator and the composite estimator. We extend the theory of the RAMP to apply to the model-averaged estimator; see Theorem~\ref{thm:amse_mod_avrg}. Challenges arise with incorporating dependence into the AMSE expression; see Theorem~\ref{thm:cov_like_estimator}.
Theorem~\ref{thm:cov_like_estimator}, in turn, leads to a new Stein-type estimator of RAMPs asymptotic MSE. While we focus on the weight choice, the availability of an estimated AMSE may be used in other contexts, for instance, for the construction of confidence intervals.

\subsection{Notation}

When the composite loss function $\rho_C=\sum_{k=1}^K w_k\rho_k$ is used in the RAMP algorithm with tuning parameter $\alpha$ we denote the estimator at iteration number $t$ by $\widehat\beta_{\text{C},(t)}(\alpha)$.
When the value of the tuning parameter is clear from the context, we also denote the RAMP estimator by $\widehat\beta_{\text{C},(t)}$.

For constructing the model averaging estimator we denote the separate estimators from the RAMP algorithm using regularity parameters $\alpha_k$, $k=1,\ldots,K$ by $\widehat\beta_{k,(t)}(\alpha_k)$ and the model-averaged estimator is denoted by $\widehat\beta_{\text{MA},(t)}(\bm{\alpha})=\sum_{k = 1}^K w_k \widehat\beta_{k, (t)}(\alpha_k)$ with $\bm{\alpha}=(\alpha_1,\ldots,\alpha_K)^\top$.
When the value of the tuning parameters is clear from the context, we denote the model averaging RAMP estimator by $\widehat\beta_{\text{MA},(t)}$.

A generic estimator, without referring to a specific loss function or construction, is denoted by $\widehat{\beta}_{(t)}$, using tuning parameter $\alpha$; the subscript $(t)$ refers to the iteration number.

\subsection{The robust approximate message passing algorithm}\label{section:algo}

We first revise the (robust) approximate message passing algorithm, which consists of three steps  iterated until convergence. In comparison with the more straightforward AMP for the case with a differentiable convex loss function \citep{donoho2009message}, this procedure for robust high-dimensional parameter estimation \citep{donoho2016high,bradic2016robustness} adjusts the residuals to incorporate the valid score directly.
While more details are given in Algorithm~\ref{algo_single_quantile}, which is applied to the different loss functions $\rho_1,\ldots,\rho_K$ and to their weighted sum $\rho_C=\sum_{k=1}^K w_k\rho_k$, we here provide the main outline.
The used notation does not explicitly indicate a dependence on the number of coefficients $p$ to not overcomplicate the formulas.

\citet{donoho2016high} proposed to use the following proximal mapping operator to adjust the residuals. With $b>0$,
$$\mbox{Prox}(z, b) = \arg \min_{x \in \mathbbm R} \{ b\rho(x) + \frac{1}{2}(x-z)^2 \}$$
which minimizes the square loss regularized by the non-differentiable loss, $\rho$. The parameter $b$ controls how the proximal operator map points to the minimum of the non-differentiable loss, where small values correspond to a small movement towards the minimum of $\rho$. The fixed point solution of the proximal operator coincides with the minimum of the loss function $\rho$. For more information, see \citet{ParikhBoyd2014}.

We continue with the worked out example on quantile regression, see \eqref{eq:CQR-loss}.
For $\ell = 0, \ldots, K$, define
\begin{equation}\label{eq:hl}
h(\ell) = \sum_{k =1}^\ell w_k \tau_k - \sum_{k = \ell+1}^K w_k(1-\tau_k),
\end{equation} where we define a summation sign to be equal to zero in the case where the upper summation index is smaller than the lower one, that is, $\sum_{i=a}^{b} x_i = 0$ if $b<a$.
The proximal operator for the composite quantile case, see \eqref{eq:CQR-loss}, is
\begin{equation}\label{eq:proximal cqr}
\mbox{Prox}(z;b)=
\left\{\begin{array}{ll}
  z - b h(\ell), & z \in (u_{\tau_\ell} + b h(\ell), u_{\tau_{\ell+1}} + b h(\ell)),\  \ell = 0, \ldots, K \\
   u_{\tau_{\ell}} & z\in [u_{\tau_\ell} + b h(\ell-1), u_{\tau_\ell} + b h(\ell)],\  \ell = 1, \ldots,  K.
\end{array}\right.
\end{equation}
See Section~\ref{proof:proxcqr} for the derivation of the algorithm.

We now describe the three steps in more detail.

\textbf{Step 1: Create adjusted residuals.} \\
We use the estimates $\widehat\beta_{(t-1)}$ and $\widehat{\beta}_{(t)}$ from iteration steps $t-1$ and $t$ to compute the adjusted residuals
  \begin{eqnarray}\label{eq:adjust_resid}
    \lefteqn{z_{(t)} =  Y - {X}\widehat\beta_{(t)} + }
      \\& n^{-1} G(z_{(t-1)}; b_{(t-1)})  \nonumber
      \displaystyle\sum_{j=1}^p I\left\{ \eta\big(\widehat\beta_{(t-1),j} + {X}_{\cdot j} G(z_{(t-1)}; b_{(t-1)}); \theta_{t-1}\big) \not= 0\right\},
    \end{eqnarray}
where the soft-thresholding function $\eta(x; \theta ) = \mbox{sign}(x)\max(|x|-\theta,0)$ and the score function $G$ is defined in \eqref{eq:scaled_eff_score}.

In Algorithm~\ref{algo_single_quantile}, see Section~\ref{section:state_evolution}, we give details on how to set the soft-thresholding parameter $\theta$, which might change in each iteration, and we explain that
a proper choice of $\theta$ as a function of the regularity constant $\lambda$ leads to an equivalence of the RAMP estimator and the regularized estimator.

The effective score function used in \citet{donoho2016high} is
\begin{eqnarray}
\label{eq:Gtilde}
\widetilde G(z;b) = b \cdot \partial\rho(x)|_{x = {\rm Prox}(z;b)}, \mbox{ with } b>0;
\end{eqnarray}
a subgradient is used in case of nondifferentiability. That is, for a value $x$ where $\rho$ is non-differentiable
$$\partial\rho(x) = \{y: \rho(u)\ge \rho(x) + y(u-x), \forall u \}.$$
Throughout, we use $\partial_1$ as the notation for the partial derivative or partial subgradient of a function with respect to its first argument.
Functions (e.g.~$\widetilde G$) are applied componentwise to vectors.

For the example on composite quantile regression
the subgradient of $\rho_{\rm{C}}$ is computed as,
\begin{equation} \label{eq:subgradient_cqr}
 \partial \rho_{\rm{C}}(x)
\left\{\begin{array}{ll}
 = h(\ell), & x \in (u_{\tau_\ell} , \ u_{\tau_{\ell+1}}), \mbox{ for } \ell = 0, \ldots, K, \\
  \in [h(\ell-1), h(\ell)], & x = u_{\tau_\ell}, \mbox{ for } \ell = 1, \ldots, K,
 \end{array}\right.
\end{equation}
where $h(\ell)$ is defined in \eqref{eq:hl}.
The effective score function for composite quantile regression, see Section~\ref{Proof:effscore}, is
\begin{equation}\label{eq:CQR_effective_score}
\widetilde G(z;b) =
\left\{\begin{array}{ll}
  b h(\ell), & z \in (u_{\tau_\ell} + b h(\ell), u_{\tau_{\ell+1}} + b h(\ell)),\ \ell = 0, \ldots, K \\
  z - u_{\tau_\ell}, & z\in [u_{\tau_\ell} + b h(\ell-1), u_{\tau_\ell} + b h(\ell)],\ \ell = 1, \ldots,  K.
\end{array}\right.
\end{equation}

To incorporate the sparsity, \citet{bradic2016robustness}, see also \citet{bayati2011dynamics},  used the rescaled, min regularized effective score function,
\begin{eqnarray}\label{eq:scaled_eff_score}
G(z;b) = \delta\omega^{-1} \widetilde G(z;b)
\end{eqnarray}
where $\omega=E[\|B_0\|_0]$, see condition~\ref{cond:beta} in the Appendix, which corresponds to the limit of $s/p$, with $s=\|B_0\|_0$, the true number of nonzero components, as $p$ tends to infinity.

\textbf{Step  2: Use the effective score function to set $b$.}
\\
We choose the scalar $b_{(t)}$ such that the empirical average of the effective score function $G(z; b)$ has slope 1, thus
$n^{-1}\sum_{i = 1}^n\partial_1 G(z_{i,(t)}; b_{(t)})=1$.
In the case of a non-differentiable loss function, \citet{bradic2016robustness} proposed to solve $\widehat \nu(b_{(t)})=1$ with
\begin{eqnarray}\label{eq:estimator_nu}
    \widehat\nu(b_{(t)}) = \frac{b_{(t)}\delta}{\omega}\big(\frac{1}{n}\sum_{i=1}^n\sum_{j=1}^2 \partial v_{j} \{z_{(t), i}\} +
    \sum_{l=1}^{L-1}\gamma_l  \{\widehat
    f_{P}(r_{l+1}) - \widehat f_{P}(r_{l})\}\big).
\end{eqnarray}
Condition~\ref{cond:loss} in the Appendix defines $\gamma_l, r_l$ and the differentiable functions $v_1$ and $v_2$
\citep[See also Condition (R) of][]{bradic2016robustness}, $\widehat f_P$ is the estimated density of $\mbox{Prox}(z_{i,(t)};b_{(t)})$ for $i=1,\ldots,n$.

The derivation of the estimator $\widehat\nu(b_{(t)})$, see also Section~\ref{section:estimator_nu}, relies on the limiting behaviour of the system, see Section~\ref{section:state_evolution}.

For the composite quantile loss, see \eqref{eq:CQR-loss}, we clearly see the dependence on the quantiles. The estimator of $\nu$ in \eqref{eq:estimator_nu} uses $v_1(z) = 0$ and $v_2(z) = z - u_{\tau_\ell}$ $z\in [u_{\tau_\ell} + b h(\ell-1), u_{\tau_\ell} + b h(\ell)],\ \ell = 1, \ldots,  K$, corresponding to the differentiable pieces in \eqref{eq:CQR_effective_score}. The step functions $v_3(z)=bh(\ell)$ when $z \in (u_{\tau_\ell} + bh(\ell), u_{\tau_{\ell} + 1} + bh(\ell)),\ \ell = 0, \ldots, K$.
Solving for $b$ in the equation $\widehat\nu(b)=1$ is equivalent to solving for $b$ in the following equation,
\begin{eqnarray}\label{eq:eff_score_b}
 \frac{s}{n} &=& b\big[\sum_{k=0}^{K-1} h(k) f_z\{u_{\tau_{k+1}} + b h(k)\} -
 \sum_{k=1}^K h(k) f_z\{u_{\tau_k} + b h(k)\}\big]   \nonumber\\
 && +F_z\{b h(K)\} - F_z\{b h(0)\},
\end{eqnarray}
where $F_z$ is the cumulative distribution function and
$f_z$ the density function of the adjusted residuals. In practice, a grid search is performed to approximate the solution $\widehat b_t$. For each $b$ in the grid, we use the empirical cumulative distribution, that is, $\widehat F_z(bh(K)) = n^{-1}\sum_{i=1}^n I\{ z_{i;(t)} \leq b h(K)\}$. A kernel density estimator of $f_z$ with the Gaussian kernel estimates defined as
$\widehat f_z\{u_{\tau_k} +bh(k)\} = (nh)^{-1} \sum_{i =1}^n \phi\{ (z_{i;(t)} -u_{\tau_k} -bh(k))/h \}$ with $\phi$ being the standard normal density function.
 The solution $\widehat b_t$ is taken to be the average of the smallest $b$ in the grid that makes the righthand side of (\ref{eq:eff_score_b}) smaller than $\frac{s}{n}$ and the next value in the grid.

\textbf{Step  3: Update the estimator of $\beta$.}
\\
Use the estimated $b_{(t)}$ from the previous step to update the estimate of $\beta$ to
    \begin{eqnarray}\label{eq:pseudo_data and eq:esti}
     \widehat\beta_{(t+1)} = \eta (\widetilde\beta_{(t)}; \theta_{(t)}), \mbox{ where }   \widetilde\beta_{(t)} = \widehat\beta_{(t)} + {X}^\top G(z_{(t)}; b_{(t)}).
    \end{eqnarray}

The estimator $\widetilde\beta_{(t)}$, before applying the soft-thresholding function, is of interest too since it can be interpreted as a debiased estimator
\citep{javanmard2014confidence, javanmard2014hypothesis, van2014asymptotically, javanmard2018debiasing}; a thorough study of which, however, is beyond  the current work.

\allowdisplaybreaks
\begin{algorithm}[!t]
   \SetAlgoLined
   \SetKwFunction{FMain}{singleRAMP}
   \caption{RAMP algorithm for a single loss function with tuning parameter $\alpha$}\label{algo_single_quantile}
   \Fn{\FMain {$\alpha$}}{
          \Initialization{$\widehat\beta_{(0)} \gets 0 \in \mathbbm R^p$, \\
            iteration index $t \gets 0$, final iteration $t_{\rm final} \gets 0$,\\
            adjusted residuals $z_{(0)} \gets Y\in \mathbbm R^n$,\\
            empirical state evolution $\bar\zeta_{(0)}^2$ using (\ref{eq:state_evo_zeta_empirical}),\\
            tuning parameter of the soft-thresholding function $\theta_{(0)} = \alpha \bar\zeta_{(0)}$
         }
    \While{iteration $t \leq T$ and tolerance $tol > \varepsilon_{\rm tol}$ }
    {\begin{enumerate}
        \item \textit{Adjust residuals:}  adjust the residuals $z_{(t)} \in \mathbbm R^n$:
      \begin{eqnarray*} \hspace*{-10mm}
      z_{(t)} \gets Y - {X}\widehat\beta_{(t)}  +  \frac{1}{n} G(z_{(t-1)}; b_{(t-1)})
      \sum_{j=1}^p I\left\{ \eta \big(\widehat\beta_{j,(t-1)} + {X}_{\cdot j}^\top G(z_{(t-1)}; \right.
      \\  \left.  b_{(t-1)}); \theta_{(t-1)}\big) \neq 0 \right\}.
                   \end{eqnarray*}
        \item \textit{Effective score:}
                  \begin{enumerate}
                     \item choose the scalar $b_{(t)}$ satisfying \linebreak
                       \leIf{$G$ \rm differentiable}{
                         $1= \frac{1}{n}\sum_{i = 1}^n\partial_1 G(z_{i,(t)}; b_{(t)})$\linebreak}
                        {$1 = \widehat \nu(b_{(t)})$, see (\ref{eq:estimator_nu})}

                      \item update the state evolution parameter $\bar\zeta_{(t)}^2$ using (\ref{eq:state_evo_zeta_empirical})
                      \item update the tuning parameter $\theta_{(t)} \gets \alpha \bar\zeta_{(t)}$.
                      \end{enumerate}
        \item \textit{Estimation:}
                 Update the coefficient estimation
                        \begin{eqnarray*}
                          \widetilde\beta_{(t)} \gets \widehat\beta_{(t)} + {X}^\top G(z_{(t)}; b_{(t)}) \mbox{ and }
                                                   \widehat\beta_{(t+1)} \gets \eta (\widetilde\beta_{(t)}; \theta_{(t)}),
                        \end{eqnarray*}
        \item \textit{Adjust iteration index:}
                 $t \gets t+1$; $t_{\rm final} \gets t$.
        \item \textit{Calculate tolerance:}
                $tol = \|\widehat\beta_{(t)} - \widehat\beta_{(t - 1)}\|^2 /p $
    \end{enumerate}
   }
  \KwRet{$\widehat\beta \gets \widehat\beta_{(t_{\rm final})}$,
  $\widetilde\beta \gets \widetilde\beta_{t_{\rm final}}$,
  \text{the estimated $\mbox{AMSE}(\widehat\beta; \beta)$ for $\widehat\beta$, see Theorems \ref{thm:amse_mod_avrg} and \ref{thm:cov_like_estimator}.}}
  }
 \end{algorithm}

\section{State evolution}\label{section:state_evolution}

Within each iteration step $t$ of the approximate message passing algorithm, state evolution studies the limiting behaviour of the estimators when the sample size goes to infinity.
We now define the state evolution parameter $\bar\zeta_{(t)}^2$
which is critical for Algorithm~\ref{algo_single_quantile}. We start by defining the empirical version as follows
\begin{eqnarray}\label{eq:state_evo_zeta_empirical}
  \bar\zeta_{\rm{emp}, (t)}^2 = \frac{1}{n}\sum_{i = 1}^n G(z_{i, (t)} ; b_{(t)})^2.
 \end{eqnarray}

This quantity is linked to the state evolution recursion which describes the limiting behaviour of large systems, see Theorem 2 in \citet{bayati2011dynamics} and Lemma 1 in \citet{bradic2016robustness}.
It holds that, see Section~\ref{Proof:ofconv_eff_score} for details,
\begin{equation}\label{eq:converge_eff_score}
\lim_{n\to \infty} \frac{1}{n} \sum_{i = 1}^n G \big(z_{i, (t)}; b_{(t)}\big)^2
    \stackrel{a.s.}{=} E [G(\varepsilon  - \bar\sigma_{(t)}Z; b_{(t)})^2]=\bar\zeta_{(t)}^2,
\end{equation}
where $\bar\zeta_{(t)}$
  is the state evolution parameter for the large system, $Z$ is a random variable with standard normal distribution independent of everything else and $\bar \sigma(t)$ is defined in \eqref{eq:state_evolution_sigma}.

 Due to the symmetry of $Z$, the state evolution parameter is formally defined as
  \begin{equation}\label{eq:state_evolution_zeta}
      \bar\zeta_{(t)}^2 = E [G(\varepsilon  + \bar\sigma_{(t)}Z; b_{(t)})^2].
  \end{equation}
 This definition explicitly features the extra Gaussian component $\bar\sigma_{(t)}Z$ in the limiting version, with variance
  \begin{equation}\label{eq:state_evolution_sigma}
    \bar\sigma_{(t)}^2 = \delta^{-1} E[ (\eta(B_0 + \bar\zeta_{(t-1)}Z; \theta_{(t - 1)} ) -  B_0 )^2]
  \end{equation}
   with $B_0$ defined in \ref{cond:beta}. To connect the theoretical expression of $\bar\sigma_{(t)}^2$ to Algorithm~\ref{algo_single_quantile}, we apply Eq.(3.6) in \citet{bayati2011dynamics}, and Eqs.(7.10) and (7.19) in \citet{bradic2016robustness}. This leads to
\begin{eqnarray}\label{eq:state_evolution_as_sigma}
 \delta^{-1}\lim_{p \to \infty} \frac{1}{p} \sum_{j = 1}^p \{\eta(\widehat\beta_{(t),j} + X_{\cdot j}^\top G(z_{i, (t)};b_{(t)}); \theta_{(t)}) - \beta_{j}
    \}^2 \stackrel{a.s.}{=}\bar\sigma_{(t)}^2.
\end{eqnarray}
Note that \eqref{eq:state_evolution_as_sigma} features
the debiased estimator from \eqref{eq:pseudo_data and eq:esti}.

We now explain the connection between the estimators that explicitly use an $l_1$-regularization and the corresponding estimators from the RAMP algorithm.

By applying Theorem~2 of \citep{bradic2016robustness}, we get the immediate connection between the regularized estimators $\widehat{\beta}_k(\lambda_k)$ for $k=1,\ldots, K$ and the corresponding estimators obtained by applying the RAMP algorithm with a suitable choice of its regularity parameter $\alpha$. We explain this below.
Since the regularized estimators $\widehat{\beta}_k(\lambda_k)$ for $k=1,\ldots,K$ are used for $\widehat{{\beta}}_{{\rm MA}}({\bm\lambda})$, \eqref{eq:MA}, the connection between the model-averaged estimators from regularization and from application of the RAMP algorithm, follows immediately from the connections between the $K$ separate estimators. The composite estimator $\widehat{{\beta}}_{{\rm C}}(\lambda)$, \eqref{eq:Composite}, is a special case of a model-averaged estimator with $K = 1$, weight equal to one, and loss function $\rho_{\rm C}=\sum_{k=1}^Kw_k\rho_k$.

Denote $(\bar\zeta^2, b)$ as the fixed point solution when the iteration number $t \to \infty$ of the following equations,
\begin{eqnarray}
\label{eq:zeta_n>p}
    \bar\zeta_{(t)}^2 = E[G (\varepsilon + \bar\zeta_{(t)}Z ; b_{(t)})^2] &= &
    (\delta/\omega)^2 E[\widetilde G (\varepsilon + \bar\zeta_{(t)}Z ; b_{(t)})^2]
\\ \label{eq:effective_score_n>p}
    1 = E[\partial_1 G (\varepsilon + \bar\zeta_{(t)}Z ; b_{(t)})] &=& (\delta/\omega) E[\partial_1 \widetilde G (\varepsilon + \bar\zeta_{(t)}Z_k ; b_{(t)})].
\end{eqnarray}
Note that \eqref{eq:zeta_n>p} is the state evolution recursion for the large system while in  \eqref{eq:effective_score_n>p} the first equality is the population version of the requirement in step 2 in Algorithm~\ref{algo_single_quantile} which states that
$n^{-1}\sum_{i = 1}^n\partial_1 G(z_{i,(t)}; b_{(t)})=1$. The second equalities of both \eqref{eq:zeta_n>p} and \eqref{eq:effective_score_n>p} follow by using the definition of $G$ in \eqref{eq:scaled_eff_score}, with $\widetilde G$ being defined in \eqref{eq:Gtilde}.

Then, under assumptions \ref{Afirst}--\ref{Alast} (see the Appendix), for the RAMP algorithm with $\theta=\alpha \overline\zeta$, where the tuning parameter $\alpha >0$ (which motivates the definition of $\theta_{(t)}=\alpha \bar\zeta_{(t)}$ in Algorithm~\ref{algo_single_quantile}), and for the $l_1$-optimization with
 \begin{eqnarray}\label{eq:lambda-alpha}
 \lambda = \frac{\alpha \overline\zeta}{b\delta} P(|B_0+\overline\zeta Z|\ge \alpha\overline\zeta),
 \end{eqnarray}
it follows by Theorem 2 of \citet{bradic2016robustness} that
\begin{eqnarray}
\label{eq:equivalence}
\lim_{t\to\infty}\lim_{p\to\infty}\frac{1}{p}\sum_{j=1}^p \{\widehat{\beta}_{{\rm C},j}(\lambda)- \widehat{\beta}_{{\rm C},(t),j}(\alpha)\}^2 = 0 \mbox{ a.s.}
\end{eqnarray}

The convergence in \eqref{eq:equivalence} explicitly connects the two composite estimators: one estimator uses an explicit $l_1$-regularization as in \eqref{eq:Composite}, the other estimator is obtained via the RAMP algorithm.
 Similar results can be found in \citet[Theorem 2.2]{huang2020robust} for a generalized AMP algorithm with non-negative convex loss function, and in \citet[Theorem 1.8]{bayati2011lasso} for the AMP algorithm with least squares loss function.

For the model averaging estimator we use such an equivalence for estimation with each separate loss function $\rho_k$, $k=1,\ldots,K$.
When using explicit $l_1$-regularization as in \eqref{eq:MA_single} with the regularization constants $\lambda_k$ matching as in \eqref{eq:lambda-alpha} the values $\theta_k=\alpha_k\bar\zeta$, for $k=1,\ldots,K$ that are used in the RAMP algorithm, again Theorem 2 of \citet{bradic2016robustness} applies. It hence follows that
 $$
 \lim_{t\to\infty}\lim_{p\to\infty}\frac{1}{p}\sum_{j=1}^p \{\widehat{\beta}_{{\rm MA},j}(\bm\lambda)- \widehat{\beta}_{{\rm MA},(t),j}(\bm\alpha)\}^2 = 0, \mbox{a.s.} $$


\section{Theoretical contributions} \label{sec:AMSE}
This section contains detailed theoretical developments for the composite as well as the model-averaged AMP estimators in high-dimensions.

\subsection{Asymptotic mean squared error} \label{section:theoretical expression}
We first define the asymptotic mean squared error as
\begin{equation}\label{eq:single_amse}
    \AMSE(\widehat\beta_{(t)}, \beta) = \lim_{p\to\infty} \frac{1}{p}\sum_{j = 1}^p(\widehat\beta_{(t),j} - \beta_j)^2.
\end{equation}
Combining (\ref{eq:state_evolution_as_sigma}) and \eqref{eq:converge_eff_score}, we obtain
\begin{eqnarray}\label{eq:single_theoretical_amse}
\AMSE(\widehat\beta_{(t)}, \beta)
&=&
\lim_{p \to \infty} \frac{1}{p} \sum_{j = 1}^p \Big(\eta(\widetilde\beta_{(t - 1),j} ; \theta_{(t-1)}) - \beta_{j}\Big)^2 \nonumber\\
&&\stackrel{a.s.}{=} E[\{\eta(B_0 - \bar\zeta_{(t-1)}Z; \theta_{(t - 1}) ) -  B_0 \}^2],
\end{eqnarray}
which corresponds to Eq.(3.4) in \citet{bradic2016robustness} with $\widetilde\beta_{(t),j}$ the debiased estimator in (\ref{eq:pseudo_data and eq:esti}).

In Section~\ref{section:state_evolution}, we defined the empirical state evolution parameter $\bar\zeta_{\rm{emp}, (t)}^2$, and we described the connections between the empirical updates in Algorithm~\ref{algo_single_quantile} and the theoretical state evolution recursion, which connects to the theoretical expression of the AMSE. While Algorithm~\ref{algo_single_quantile} and the theoretical state evolution recursion involve only a single estimator, the model-averaged estimator, on the other hand, is the weighted sum of $K$ such estimators $\widehat{\beta}_k$, $k=1,\ldots,K$, each obtained by Algorithm~\ref{algo_single_quantile}. Consequently, the estimators $\widehat{\beta}_k$, $k=1,\ldots,K$ are correlated.

Lemma~\ref{lemma:converge_cov} extends Theorem 2 in \citet{bayati2011dynamics} and (3.16) in Lemma 1(b) in \citet{bayati2011dynamics} to the almost sure convergence of the product for any two recursions among $K$ paralleled recursions.
All proofs are contained in Appendix~\ref{sec:proofs}.

\begin{lemma}
\label{lemma:converge_cov}
Let the sequences of design matrices $\{X(p)\}$, coefficient vectors $\{\beta(p)\}$, error vectors $\{\varepsilon(p)\}$, initial condition vectors $\{q_0(p)\}$ be the common sequences for $K$ recursions satisfying assumptions \ref{cond:design}--\ref{A4} in the Appendix. Let $\{\bar\sigma_{k, (t)}^2, \bar\zeta_{k, (t)}^2\}$ be defined uniquely by the recursions in (\ref{eq:state_evolution_zeta}) and (\ref{eq:state_evolution_sigma}). These are the  state evolution parameters for the $k$th estimation with initialization $\bar\sigma^2_{k, (0)} = \lim_{n\to \infty}\frac{1}{n} \sum_{i = 1}^n q_{(0), i}^2 /\delta$. Then Lemma 1 in \citet{bayati2011dynamics} holds individually for each of the $K$ recursions; additionally, for all pseudo-Lipschitz functions $\tilde\psi_{\rm c}: \mathbbm R^{t+2} \to \mathbbm R$ of order $\kappa_{\rm c}$ for some $1 \leq \kappa_{\rm c} \leq \kappa/2$ with $\kappa$ as in \ref{cond:lemma_convergence} and $t$ a natural number larger than or equal to 0,
\begin{eqnarray*}
{\lim_{p\to \infty} \frac{1}{p}\sum_{j = 1}^p \tilde\psi_{\rm c}(h_{k_1, (1), j},\ldots, h_{k_1, (t+1), j}, \beta_{j})
 \tilde\psi_{\rm c}(h_{k_2, (1), j}, \ldots, h_{k_2, (t+1), j}, \beta_{j})} \stackrel{a.s.}{=}\\
 E[\tilde\psi_{\rm c}(\bar\zeta_{k_1, (0)}Z_{k_1, (0)},\ldots, \bar\zeta_{k_1, (t)}Z_{k_1, (t)}, B_0)
 \tilde\psi_{\rm c}(\bar\zeta_{k_2, (0)}Z_{k_2, (0)},\ldots, \bar\zeta_{k_2, (t)}Z_{k_2, (t)}, B_0)]
\end{eqnarray*}
where $(Z_{k, (0)}, \ldots, Z_{k, (t)}) \sim \mathcal{N}(0, I_{t+1})$, $k = k_1, k_2$, is a $(t+1)$-dimensional zero-mean multivariate standard normal vector independent of $B_0$, $\varepsilon$; at iteration $t$, $(Z_{k_1, (t)}, Z_{k_2, (t)})$ is a bivariate standard normal vector with covariance not necessarily equal to zero.
\end{lemma}

Note that Algorithm~\ref{algo_single_quantile} belongs to the general recursion in \citet{bayati2011dynamics}, the initial condition takes $q_{(0)} = - \beta$ and the $k$th estimator calculated by Algorithm \ref{algo_single_quantile} takes $h_{k, (t+1)} = \beta - X^\top G(z_{k, (t)}; b_{k, (t)}) - \beta_{k, (t)}$.

We obtain at iteration $t$, for $k_1, k_2\in\{1,\ldots,K\}$,
\begin{eqnarray*}\label{eq:cov_like_convergence}
\lefteqn{ \lim_{p \to \infty} \frac{1}{p}\sum_{j = 1}^p(\widehat\beta_{k_1, (t), j}- \beta_j) (\widehat\beta_{k_2,(t), j}- \beta_j)  } \\
&& \stackrel{a.s.}{=} E\Big[\prod_{r=1}^2 \big\{\eta(B_0 + \bar\zeta_{k_r, (t-1)} Z_{k_r}; \theta_{k_r, (t - 1)}) - B_0 \big\} \Big],
\end{eqnarray*}
where $Z_{k_1}$ and $Z_{k_2}$ are possibly dependent standard normal random variables.

Since the estimators $\widehat\beta_{k_r}$, $r=1,2$ use the same design matrix, a correlation between $Z_{k_1}$ and $Z_{k_2}$ exists (see Corollary~\ref{cor:cov_Zs}) and contributes to the correlation between $\widehat \beta_{k_1}$ and $\widehat \beta_{k_2}$.
Using Lemma~\ref{lemma:converge_cov},
 we obtain the theoretical AMSE for the regularized model-averaged estimator.

\begin{theorem} \label{thm:amse_mod_avrg}
Assume conditions \ref{Afirst}--\ref{Alast} in the Appendix.
At Algorithm~\ref{algo_single_quantile}'s iteration step $t$ for the estimator $\widehat\beta_{k,(t)}$, for each $k=1,\ldots,K$,
and for a weight vector $w = (w_1, \ldots, w_K)^\top$,
the model-averaged estimator $\widehat\beta_{\MA,(t)} = \sum_{k = 1}^K w_k \widehat\beta_{k,(t)}$
has asymptotic mean squared error
\begin{eqnarray}\label{eq:ma_amse_theoretical}
 \AMSE(\widehat\beta_{\MA, (t)}, \beta)
   &=& \lim_{p \to \infty} \frac{1}{p} \sum_{j = 1}^p (\widehat\beta_{\MA, (t), j} - \beta_{j})^2 \nonumber \\
   &=& \lim_{p \to \infty} w^\top \Sigma_{0, (t)}(p) w
   \stackrel{a.s.}{=} w^\top \Sigma_{(t)} w
\end{eqnarray}
 where $\Sigma_{0, (t)}(p)$ is a $K\times K$ matrix with $(k_1, k_2)$th component
 \begin{equation}\label{eq:empirical_sigma}
 (\Sigma_{0, (t)})_{(k_1, k_2)}(p) = p^{-1}\sum_{j=1}^p (\widehat\beta_{k_1, (t), j} - \beta_j)(\widehat\beta_{k_2, (t), j} - \beta_j);
 \end{equation}
 similarly, $\Sigma_{(t)}$ is a $K\times K$ matrix with the $(k_1, k_2)$th component
$$
 (\Sigma_{(t)})_{(k_1, k_2)} = E\Big[\prod_{r=1}^2\{\eta(B_0 + \bar\zeta_{k_r, (t-1)} Z_{k_r}; \theta_{k_r, (t-1)}) - B_0 \} \Big].
$$
\end{theorem}
Since the AMSE expression of the regularized model-averaged estimator is a quadratic function of the weight vector $w$, Corollary~\ref{thm:ma_lower_bound} readily provides the lower bound of the AMSE as well as the weight vector reaching this lower bound.
The $K$-vector $\mathbf 1_{K}$ consists of ones only.

\begin{corollary}\label{thm:ma_lower_bound}
Constraining the weights to sum to one, the lower bound of the AMSE at iteration $t$ for the model-averaged estimator as in (\ref{eq:ma_amse_theoretical}) is equal to
$\big( \mathbf 1_{K}^\top (\Sigma_{(t)})^{-1} \mathbf 1_{K} \big)^{-1}$.
This lower bound is attained for the theoretical optimal weights $w_{\rm MA} = \big(\Sigma_{(t)}\big)^{-1} \mathbf 1_{K} \big(\mathbf 1_{K}^\top (\Sigma_{(t)})^{-1} \mathbf 1_{K} \big)^{-1}$.
\end{corollary}

\subsection{Estimating optimal weights}\label{ssec:Sigma_estimator}

The expression of the core matrix $\Sigma_{(t)}$, which is the limit matrix for $n, p \to \infty$, contains the random variable $B_0$ which satisfies assumption \ref{cond:beta} in the Appendix. Likewise,  $\Sigma_{0, (t)}$ which is the limit matrix for fixed $p$ while $n \to \infty$, contains the true coefficient $\beta$ (see (\ref{eq:empirical_sigma})).
In practice, neither the true coefficient vector $\beta$ nor the random variable $B_0$ is known.
To make practical use of the expressions of the AMSE, we derive an estimator of the matrix $\Sigma_{0,(t)}$ relying only on sequences generated in Algorithm~\ref{algo_single_quantile}.

\subsubsection{Model-averaged estimator}
Before deriving the estimator of the AMSE for the model-averaged estimator, we first define
$\bar\zeta_{\textrm{emp}, (k_1, k_2), (t)}$ which is an estimator of the parameter $\bar\zeta_{(k_1, k_2), (t)}$, a quantity similar to the state evolution parameter $\bar\zeta_{k, (t)}^2$, which records the covariance between the unbiased sequences $\widetilde \beta_{k_1, (t)}$ and $\widetilde \beta_{k_2, (t)}$ generated in (\ref{eq:pseudo_data and eq:esti}) in Algorithm~\ref{algo_single_quantile} when $p \to \infty$.
Since model-averaged estimators combine estimators constructed from the same data into one weighted average, the correlation between
$\widehat\beta_{k_1}$ and $\widehat\beta_{k_2}$ is needed to understand the AMSE of the model-averaged estimator.

Notice that the unbiasedness of the sequence $\widetilde\beta_{k, (t)}$ follows from the argument that $\widetilde\beta_{k, j, (t)}$ converges weakly to $B_0 + \bar\zeta_{k, (t)} Z_k$ when $p \to \infty$, while assigning $1/p$ point mass to each entry of the vector. Then, $\widetilde \beta_{k, j, (t)} | (B_0 = \beta_j) \sim N(\beta_j, \bar\zeta^2_{k, (t)})$ for large $p$, indicating  that $\widetilde \beta_{k, j, (t)}$ centers at $\beta_j$ ensuring the unbiasedness. Moreover, the vector $\widetilde \beta_{k, (t)}$ has  Gaussian distribution. By applying the soft-thresholding function $\eta$ on $\widetilde \beta_{k, j, (t)}$ in Lemma~\ref{lemma:stein's lemma}, we avoid the usage of the true coefficient vector $\beta$ in $\Sigma_{0, (t)}$ resulting in a Stein-type risk estimator requiring only observables from Algorithm~\ref{algo_single_quantile}.
A Gaussianity argument has also been used in \citet{bayati2011lasso, bayati2013estimating, mousavi2013parameterless, mousavi2018consistent} to derive a similar Stein-type risk estimator for the Lasso.
Details can be found in Section~\ref{section:proof of thm:cov_like_estimator}.
The bias of the estimator $\widehat\beta_{k, (t)}$ is introduced in Algorithm~\ref{algo_single_quantile} by applying the soft-thresholding function componentwise to the unbiased sequence $\widetilde\beta_{k, (t)}$.

\begin{corollary}\label{cor:cov_Zs}
Assume conditions \ref{Afirst}--\ref{Alast} in the Appendix. For any $k_1, k_2 = 1, \ldots, K$, at iteration $t$,
\begin{eqnarray*}\label{eq:converge_cov_Zs}
  \lim_{p\to \infty}\frac{1}{p}\sum_{j = 1}^p(\widetilde\beta_{k_1, (t), j} - \beta)(\widetilde\beta_{k_2, (t), j} - \beta) \stackrel{a.s.}{=}
  \bar\zeta_{k_1, (t)}\bar\zeta_{k_2, (t)} \Cov(Z_{k_1}, Z_{k_2}),
\end{eqnarray*}
 where $\bar\zeta_{k, (t)}, k = k_1, k_2$ is the state evolution parameter corresponding to $\widehat\beta_k$.
\end{corollary}

Corollary~\ref{cor:cov_Zs} indicates both the existence and a feasible estimation of the covariance between $Z_{k_1}$ and $Z_{k_2}$.
As an estimator for
$$\bar\zeta_{(k_1, k_2), (t)} = \bar\zeta_{k_1, (t)}\bar\zeta_{k_2, (t)} \Cov(Z_{k_1}, Z_{k_2})$$
 we define
\begin{eqnarray}\label{eq:zeta.emp}
\bar\zeta_{\textrm{emp}, (k_1, k_2), (t)} = \frac{1}{p-1}\sum_{j = 1}^p \big(\widetilde{\beta}_{k_1, (t), j} - \frac{1}{p}\sum_{j = 1}^p \widetilde{\beta}_{k_1, (t), j}\big) \big(\widetilde{\beta}_{k_2, (t), j} - \frac{1}{p}\sum_{j = 1}^p\widetilde{\beta}_{k_2, (t), j} \big).
\end{eqnarray}

We now state an unbiased estimator for the matrix $\Sigma_{0, (t)}$, and a consistent estimator for the matrix $\Sigma_{(t)}$ upon convergence of Algorithm~\ref{algo_single_quantile}.
\begin{theorem}\label{thm:cov_like_estimator}
Assume conditions \ref{Afirst}--\ref{Alast} in the Appendix, and that the state evolution parameter in (\ref{eq:state_evo_zeta_empirical}) satisfies $\bar\zeta_{{\rm emp}, (t)}^2 - \bar\zeta_{{\rm emp}, (t-1)}^2 = o(1)$.
For any $k_1, k_2 = 1, \ldots, K$, define
\begin{align*}
(\widehat\Sigma_0)_{(k_1, k_2), (t)}&=
  - {\bar\zeta}_{{\rm emp}, (k_1, k_2), (t - 1)}
 + \frac{1}{p} \sum_{j = 1}^p \prod_{r=1}^2\big\{
  \eta(\widetilde{\beta}_{k_r, (t - 1 ), j} ; \theta_{k_r, (t - 1)}) - \widetilde{\beta}_{k_r, (t - 1), j} \big\}
\\
  & \qquad + {\bar\zeta}_{{\rm emp}, (k_1, k_2), (t-1)}\cdot\frac{1}{p}
 \sum_{j = 1}^p \sum_{r=1}^2 I\{|\widetilde{\beta}_{k_r, (t - 1), j}| \ge \theta_{k_r, (t - 1)}\},
\end{align*}
 with $\widetilde{\beta}_{k_1, (t-1)}$, $\widetilde{\beta}_{k_2, (t-1)}$
 in (\ref{eq:pseudo_data and eq:esti})
 Then, $(\widehat\Sigma_0)_{(k_1, k_2), (t)}$ is an unbiased estimator of component $(k_1, k_2)$ of the matrix $\Sigma_{0, (t)}$ at iteration $t$. Further, $(\widehat\Sigma_0)_{(k_1, k_2), (t)}$ is a consistent estimator of the matrix $\Sigma_{(t)}$ in Theorem~\ref{thm:ma_lower_bound}.
\end{theorem}

This new estimator can be compared to the estimator used in \citet[][Def.~2]{bayati2013estimating} and \citet[][Eq.(9)]{mousavi2018consistent} for the case of a single estimator ($K=1$).
The proof of Theorem~\ref{thm:cov_like_estimator}, see
Section~\ref{section:proof of thm:cov_like_estimator} uses Stein's lemma (see Lemma~\ref{lemma:stein's lemma}) to estimate the covariances that appear in the matrix $\Sigma_{0,(t)}$.
 The soft-thresholding function $\eta(\cdot; \theta)$ that appears in the estimator $\widehat\Sigma_{0,(t)}$ links the estimator $\widehat{\beta}_k$ to the estimator $\widetilde\beta_k$.
 The proof also uses the joint asymptotic normality of the $j$th components of the vectors $\widetilde\beta_{k_1}$ and $\widetilde\beta_{k_2}$.
The obtained estimator for $\Sigma_{0,(t)}$ in the case $K>1$ is nontrivial and new to the literature.

Estimated AMSE-optimal weights for the model-averaged estimator are obtained by using the estimator $\widehat\Sigma_{0,(t)}$ at the final iteration in Theorem~\ref{thm:cov_like_estimator}. In combination with the sum-to-one constrained weights this gives the estimated weights that minimize the estimated AMSE for the model-averaged estimator
$$\widehat{w}_{\rm MA} = \big(\widehat\Sigma_{(t)}\big)^{-1} \mathbf 1_{K} \big(\mathbf 1_{K}^\top (\widehat\Sigma_{(t)})^{-1} \mathbf 1_{K} \big)^{-1}.$$
When additional constraints such as positivity are needed, the optimal weights no longer have an explicit formula, but they are straightforward to compute, see \eqref{eq:maqr opt+ weight}.

\subsubsection{Composite estimator}
The AMSE of a composite estimator can be
obtained from Theorem~\ref{thm:amse_mod_avrg}
as a special case, treating the composite loss function as a single loss function with weight one, thus $\rho_{\rm C}=\sum_{k=1}^Kw_k\rho_k$ as in \eqref{eq:Composite}.
At iteration $t$,
$$
\Sigma_{(t)} =  E[\{\eta(B_0 + \bar\zeta_{(t-1)} Z; \theta_{(t-1)}) - B_0 \}]^2, \mbox{ and } \Sigma_{0, (t)} = p^{-1}\sum_{j = 1}^p (\widehat\beta_{(t), j} - \beta_j)^2.
$$
The matrices $\Sigma_{(t)}, \Sigma_{0, (t)}$ are now real numbers and coincide with the AMSE of the estimator in (\ref{eq:single_theoretical_amse}). We obtain the corresponding estimator for the AMSE
\begin{align}\label{eq:composite AMSE estimator}
 \widehat\Sigma_{\rm{C}, 0} &= \widehat{\AMSE}_{\rm{C}}(w) \\
 & =  - {\bar\zeta}_{\textrm{emp}}^2(w) + \frac{1}{p} \sum_{j = 1}^p \Big[ \big\{\eta (\widetilde{\beta}_{j}(w) ; \theta) - \widetilde{\beta}_{j}(w) \big\}^2
  + 2 {\bar\zeta}_{\textrm{emp}}^2(w)I\{|\widetilde{\beta}_{j}(w)| \ge \theta\} \Big].
 \nonumber
\end{align}

For the single loss function, $\rho_{\rm C}$, the estimator of AMSE in \eqref{eq:composite AMSE estimator} can be compared to the Stein-type estimator that has been obtained in Definition 2 in \citet{bayati2013estimating} for the AMP algorithm using the least squares loss, which is a particular case of Algorithm~\ref{algo_single_quantile}.

Finding optimal weights for the composite estimator is complicated. Indeed, while the model-averaged estimator
has an AMSE, which is a quadratic function in the weights, see \eqref{eq:ma_amse_theoretical}, the composite estimator and its AMSE depend on the weights in a highly nonlinear fashion; e.g., observe that the soft-thresholding function in \eqref{eq:composite AMSE estimator} depends on $w$.

Therefore, optimization of the estimated AMSE with respect to the weights proceeds numerically;
$$
w_{\rm C, 1} = \arg\min_{w} \widehat{\mbox{AMSE}}_{{\rm C}}(w).$$
See Section~\ref{ssec:optw} for more details.

\subsection{ The case of dense (non-sparse)  linear models with $n / p \to \delta \geq 1$ : asymptotic variance optimality }

\citet{donoho2016high} and \citet{el2013robust} showed that the asymptotic variance of the M-estimators in the case where $p, n \to \infty$ and $n/p \to \delta \in [1, \infty)$ contains an extra Gaussian component. Recently, \citet{lei2018asymptotics} obtained  the coordinate-wise asymptotic normality of regression M-estimators in the moderate $p/n$ regime for a fixed design matrix. In the sparse high-dimensional linear model setting where $\delta \in (0, 1)$, it was shown that the sequence $\widetilde \beta_{(t)}$ in (\ref{eq:pseudo_data and eq:esti}) follows for the Lasso estimator \citep{bayati2013estimating} a similar normal distribution with the variance containing an extra Gaussian component.
The above-mentioned literature focuses on the asymptotics for a single M-estimator; we extend the asymptotic result to the model-averaged estimator. In this section, we only characterize the asymptotic variance of the model-averaged estimator for dense linear models with $n/p \to \delta \geq 1 $, following \citet{donoho2016high}.

Under the dense linear model with $n \geq p$, the soft-thresholding function $\eta(\cdot; \theta)$ is replaced by the identity function and the ratio $\omega = E[\| B_0\|_0]=1$. Consequently, Algorithm~\ref{algo_single_quantile} is adjusted to estimate
$$\widehat\beta_k = \arg\min_{\beta \in \mathbbm R^p}\Big\{ \sum_{i = 1}^n \rho_k(Y_i - X_i\beta) \Big\},$$
where $\beta$ is dense.
It is trivial to show that Algorithm~\ref{algo_single_quantile} still belongs to the general recursion in \citet{bayati2011dynamics}.
For a single estimator at iteration $t$ denoted as $\widehat\beta_{k, (t)}$, the two state evolution parameters $\bar\zeta_{k, (t)}^2$ and $\bar\sigma_{k, (t)}^2$ coincide and Theorem 4.1 in \citet{donoho2016high} holds.
\begin{theorem}
\label{prop:asymp_variance}
Assume conditions \ref{Afirst}--\ref{Alast} in the Appendix. Let $n/p \to \delta \geq 1$ when $n, p \to \infty$. For the asymptotic variance of the model-averaged estimator $\widehat \beta_{\rm MA}$ holds that
\begin{eqnarray}\label{eq:varMADonoho}
    \lim_{n, p\to \infty} \frac{1}{p}\sum_{j = 1}^p \Var(\widehat \beta_{\rm MA, \it j})
  \stackrel{a.s}{=}
   \sum_{k_1 = 1}^K \sum_{k_2 = 1}^K \Cov(Z_{k_1}, Z_{k_2})
    \prod_{r=1}^2 \{w_{k_r} V^{1/2}(\widetilde G_{k_r}; \widetilde F_{k_r}) \}
\end{eqnarray}
for differentiable $\widetilde G$, where $V(\widetilde G_k; F_{k}) = (\int \widetilde G_k^2 dF_{k}) / (\int \partial_1 \widetilde G_k dF_{k})^2$ denotes the Huber asymptotic variance formula for M-estimators.
For non-differentiable $\widetilde G$,
we replace $V$ in \eqref{eq:varMADonoho} by the consistent estimator $\widehat V(\widetilde G_k; F_{k}) = (\int \widetilde G_k^2 dF_{k}) / \widehat\nu({b_k})^2$.
The extra Gaussian component is identified in the convolution of the regression noise distribution and a Gaussian distribution: $\widetilde F_{k} = F_\varepsilon \star N(0, \bar\zeta_{k}^2)$.
\end{theorem}
Recall that the componentwise empirical distribution of $\widehat\beta_k(p)$, when $p\to \infty$, converges weakly to $B_0 + \bar\zeta_{k} Z_k$ following \citet{bayati2011dynamics} and \citet{donoho2016high}. Then for large $p$, while the iteration $t \to \infty$, $\widehat\beta_k(p) \sim N(\beta, \bar\zeta_{k}^2 I_p)$ \citep{donoho2016high, mousavi2013parameterless} with $I_p$ the $p\times p$ identity matrix.
The $(k_1, k_2)$th component of the empirical variance matrix is denoted by $\big(\Sigma_{\rm emp}(p)\big)_{(k_1, k_2)} = p^{-1} \sum_{j = 1}^p (\widehat\beta_{k_1, j} - \beta_j)(\widehat\beta_{k_2, j} - \beta_j)$, which is unbiasedly estimated by
$$
\big(\widehat\Sigma_{\textrm{emp}}(p)\big)_{(k_1, k_2)} = \sum_{j = 1}^p (\widehat\beta_{k_1, j} - \frac{1}{p} \sum_{j = 1}^p \widehat\beta_{k_1, j})(\widehat \beta_{k_2, j} - \frac{1}{p}\sum_{j = 1}^p \widehat\beta_{k_2, j}) / (p - 1).
$$
Note that this estimator coincides with \eqref{eq:zeta.emp} for the special case that $n \geq p$ and the soft-thresholding function is replaced by the identity function.

\section{Computational details} \label{sec:examples}

\subsection{Regularized model-averaged quantile estimation}\label{section:model_avg_practicals}

The estimation  of the quantile $u_{\tau_k}= F_\varepsilon^{-1}(\tau_k)$ follows a two-step procedure.
\begin{enumerate}
    \item Obtain an initial slope estimate $\widehat\beta_{\rm init}$ and calculate the residuals. Example initial slope estimates are the Lasso or regularized quantile estimation with a single quantile level.
    \item For $k=1,\ldots,K$, estimate the quantile intercepts $\widehat u_{\tau_k}$ by taking the corresponding $\tau_k \times 100\%$ quantile of the residuals from the previous step.
    \end{enumerate}

The regularized model-averaged estimator is obtained by averaging over $K$ paralleled estimators. See Algorithm~\ref{algo_K_quantile} for the pseudo-code, of which the core is
Algorithm~\ref{algo_single_quantile}; there the effective score function $G$ is that of a single quantile loss function with $K=1$, see also Example 2 in \citet{bradic2016robustness}. In our numerical work, the upper bound for the number of iteration steps $T$ is set to be 50 in both the simulation and the data analysis sections.
With $K=1$, this algorithm applies to the regularized composite estimator too.

 \begin{algorithm}
    \SetAlgoLined
    \SetKwFunction{FMain}{KparallelRAMP}
   \caption{RAMP algorithm for $K$ paralleled estimations with tuned $\alpha$'s}\label{algo_K_quantile}
   \algorithmfootnote{AMSE refers to the estimated version. The $\widetilde \beta_k$s are recorded for calculating the weights in  Corollary~\ref{thm:ma_lower_bound}.}
    \Fn{\FMain{$K$}}{
         \For{$k$ in $\{1, \ldots, K\}$}
         {\Initialization{$\widehat\beta(\alpha_{k, \rm opt}) \gets 0 \in \mathbbm R^p$, $\widetilde\beta_k(\alpha_{k, \rm opt}) \gets 0 \in \mathbbm R^p$ and $\mbox{AMSE}(\widehat\beta_k(\alpha_{k, \rm opt}); \beta) \gets 0$}
         \For{$\alpha$ \rm{in candidate set} $\mathcal A$}{
           \texttt{singleRAMP}($\alpha$) in Algorithm \ref{algo_single_quantile} \linebreak
           \If{\rm{AMSE}$(\widehat\beta_k(\alpha); \beta) \leq$ \rm{AMSE}$(\widehat\beta_k(\alpha_{k, \rm opt}); \beta)$}{
               $\widehat\beta_k(\alpha_{k, \rm opt}) \gets \widehat\beta_k(\alpha)$,
               $\widetilde\beta_k(\alpha_{k, \rm opt}) \gets \widetilde\beta_k(\alpha)$,
               $\mbox{AMSE}(\widehat\beta_k(\alpha_{k, \rm opt}); \beta) \gets \mbox{AMSE}(\widehat\beta_k(\alpha); \beta)$
           }
         }
         }
         {\KwRet{$(\widehat\beta_1(\alpha_{1,\rm opt}), \ldots,
                    \widehat\beta_{K}(\alpha_{K,\rm opt}))$,
                    $(\widetilde\beta_1(\alpha_{1,\rm opt}), \ldots,
                   \widetilde\beta_{K}(\alpha_{K,\rm opt}))$, and
                    $(\mathrm{AMSE}(\widehat\beta_1(\alpha_{1, \rm opt}); \beta), \ldots,
                    \mathrm{AMSE}(\widehat\beta_K(\alpha_{K, \rm opt}); \beta))$
         }
         }
         }
 \end{algorithm}

The tuning parameter $\alpha$ of Algorithm \ref{algo_K_quantile} controls the sparsity of the estimators and requires a tuning procedure to choose it in practice. In Section~\ref{sec:numerical}, we consider the one dimensional Golden-section search algorithm \citep{kiefer1953sequential} for tuning the value $\alpha$ in the range $[\alpha_{\rm min}, \alpha_{\rm max}]$ that minimize the estimated MSE of $\widehat\beta$ using the estimator derived in Section~\ref{ssec:Sigma_estimator}. The upper bound $\alpha_{\rm max}$ is chosen to be 2.3 for the simulations and data analysis. The lower bound $\alpha_{\rm min}$ in the data analysis follows the lower bound in Proposition 9.2 in \citet{eldar2012compressed} and is chosen to be the unique non-negative solution to the equation
$(1 + \alpha^2)\Phi(-\alpha) - \alpha \phi(\alpha) = \delta/2$,
where $\phi(x)$ and $\Phi(x)$
denote the p.d.f and c.d.f of the standard normal distribution respectively.
In the simulation study, the lower bound $\alpha_{\rm min}$ is chosen to be 1.3 for computational efficiency purposes, since the optimal tuning parameter for those settings was rarely less than 1.3.

\subsection{Optimization of the weights}\label{ssec:optw}

To obtain the regularized model-averaged quantile estimations with the AMSE-type weight derived in Corollary \ref{thm:ma_lower_bound}, we follow the following procedure:
\begin{enumerate}
    \item Obtain optimally tuned paralleled regularized quantile estimates, see \eqref{eq:pseudo_data and eq:esti},
    $(\widehat\beta_{\tau_1}(\alpha_{1,\rm opt})$, $\ldots, \widehat\beta_{\tau_K}(\alpha_{K,\rm opt}))$, and the additional $K$ estimates $(\widetilde\beta_{\tau_1}(\alpha_{1,\rm opt}),$ $\ldots, \widetilde\beta_{\tau_K}(\alpha_{K,\rm opt}))$ from the converged iterations using Algorithm \ref{algo_K_quantile}.
    \item Estimate the AMSE-type optimal weight $\widehat w_{\rm MA, 1}$ with constraints by
    \begin{equation}\label{eq:maqr opt+ weight}
        \widehat w_{\rm MA, 1} = \arg\min_{w \geq 0, \mathbf 1_{K}^\top w = 1} w^\top \widehat\Sigma_0w
    \end{equation}
    where the $K \times K$ matrix $\widehat\Sigma_0$ is the consistent estimator of Theorem \ref{thm:cov_like_estimator}.
   \item Obtain the regularized model-averaged estimate (\ref{eq:MA}) with the estimated AMSE-type optimal weight.
\end{enumerate}
It is worth mentioning that $\widehat w_{\rm MA, 1}$ is a constrained version of $w_{\rm MA}$ attaining the lower bound of the AMSE in Corollary~\ref{thm:ma_lower_bound}. $\widehat w_{\rm MA, 1}$ focuses on approximating the lower bound of the AMSE of the sparse coefficient vector $\beta$ without assuming that the nonzero entries are selected perfectly; whereas another type of weight choice derived in \citet{bradic2011penalized, BloznelisClaeskensZhou2019} aims at the lower bound of the variance of the nonzero part of $\beta$ by imposing the perfect selection assumption. A numerical comparison of these two types of weight choices is presented in Section~\ref{sec:numerical}.

To equip the regularized composite quantile estimator with the weight minimizing the estimated AMSE, we cannot make use of an analytical solution to the weight minimization problem. Instead, a numerical search for a better weight choice in the neighbourhood of an initial weight proposal is employed. The basic idea is that the estimator $\widehat\beta_{\rm C}(w_{\rm C})$ is treated as a function of the weights.
We propose a collection of candidate weight vectors in the neighborhood of the weight chosen in the previous step. The weight for $\widehat\beta_{\rm C}(w_{\rm C})$ is updated in each step by the one having the lowest estimated AMSE, i.e.,
$$w_{\rm C, 1} = \arg\min_{w_{\rm cand}} \widehat{\mbox{AMSE}}(\widehat\beta_{\rm C}(\alpha_{\rm opt}; w_{\rm cand}); \beta).
$$

A more detailed search procedure is as follows.
 \begin{enumerate}
     \item Propose a reasonable initial weight vector $w_{\rm C, init}$, e.g.~the vector of equal weights;  estimate $\widehat\beta_{\rm C}$ at the initial weight $w_{\rm C, init}$ and obtain the  estimate of $\mbox{AMSE}\big(\widehat\beta_{\rm C}(\alpha_{ \rm opt} ;w_{\rm C, init}); \beta\big)$.
     \item Initiate the searching step calculator $s_{\mathcal D} = 0$, the candidate optimal weight $w_{\rm C, 1} = w_{\rm C, init}$, and the corresponding candidate minimum MSE
         $$ \mbox{AMSE}(w_{\rm C, 1}) = \mbox{AMSE}\big(\widehat\beta_{\rm C}(\alpha_{ \rm opt} ;w_{\rm C,init}); \beta\big) $$ estimated by the AMSE estimator in Theorem \ref{thm:cov_like_estimator} for $K = 1$, the collection of the used weight vectors $\mathcal{V}_{w} = \{w_{\rm C, 1} \}$.
     \item Propose a set of candidate weight vectors $\mathcal{V}_{w_{\rm{cand}}}$. This is to exclude  those recorded in the collection of the used weight vectors $\mathcal{V}_{w}$. In addition,  $\mathcal{V}_{w_{\rm{cand}}}$ should be in the neighborhood of the current optimal weight $w_{\rm C, 1}$. Rules of proposing candidate weight vectors are user-decided; here, we consider a $(K - 1)$-dimensional grid search centering at $w_{\rm C, 1}$.
     \item Obtain the regularized composite quantile estimates at all candidate weight vectors in $\mathcal{V}_{w_{\rm{cand}}}$ with Algorithm~\ref{algo_K_quantile}. Update the used weight vector collection $\mathcal{V}_{w}$, increase the counter $s_{\mathcal V} = s_{\mathcal V} + 1$, update the candidate optimal weight $w_{\rm C, 1}$ by the weight with the lowest estimated AMSE in $\mathcal{V}_{w} = \{w_{\rm C, 1} \}$, and update the candidate minimum AMSE value $\mbox{AMSE}(w_{\rm C, 1})$.
     \item Stop the iteration if the searching step calculator $s_{\mathcal V} > S_{\mathcal V}$ or the candidate weight vector collection $\mathcal{V}_{w_{\rm{cand}}} = \emptyset$; otherwise repeat steps 3 and 4.
 \end{enumerate}
The pseudocode of the search procedure is stated in Algorithm~\ref{algo_weight_search_composite}.

\begin{algorithm}
    \SetAlgoLined
    \SetKwFunction{FMain}{Weight Search}
   \caption{Weight search for regularized composite estimator}\label{algo_weight_search_composite}
   \algorithmfootnote{A possible initial weight vector $w_{\rm C, init}$ is the vector of equal weights or the weight proposed in \citet{bradic2011penalized}; $\widehat\beta_{C}$ is estimated by Algorithm~\ref{algo_K_quantile}, and $\mbox{AMSE}(\widehat\beta_{\rm C}; \beta)$ is estimated by (\ref{eq:composite AMSE estimator}).}
    \Fn{\FMain}{
    \Initialization{Better weight recorder $w_{\rm C, 1} \gets w_{\rm C, init}$, step calculator $s_{\mathcal V} \gets 0$, MSE recorder $\mbox{AMSE}(w_{\rm C, 1}) \gets \widehat{\mbox{AMSE}}\big(\widehat\beta_{\rm C}(\alpha_{ \rm opt} ;w_{\rm C, init}); \beta\big)$, and the collection of the used weight vectors $\mathcal{V}_{w} = \{w_{\rm C, 1} \}$.}
         \While{searching step $s_{\mathcal V} \leq S_{\mathcal V}$ or candidate weight collection $\mathcal{V}_{w_{\rm{cand}}} = \emptyset$}{
         \begin{enumerate}
             \item Propose a new $\mathcal{V}_{w_{\rm{cand}}}$ in the neighbourhood of $w_{\rm C, 1}$. Rules of proposing candidate weight vectors are user-decided; here, we consider a $(K - 1)$-dimensional grid search centering at $w_{\rm C, 1}$.
             \item \For{$w_{\rm cand}$ in $\mathcal{V}_{w_{\rm{cand}}} \cap{\mathcal{V}_{w}^\complement }$}{
             Estimate $\widehat\beta_{\rm C}(\alpha_{ \rm opt} ;w_{\rm cand})$ and $\widehat{\mbox{AMSE}}\big(\widehat\beta_{\rm C}(\alpha_{ \rm opt} ;w_{\rm cand}); \beta\big)$\\
             \If{$\widehat{\mbox{AMSE}}\big(\widehat\beta_{\rm C}(\alpha_{ \rm opt} ;w_{\rm cand}); \beta\big) < \mbox{AMSE}(w_{\rm C, 1})$}{
         $w_{\rm C, 1} \gets w_{\rm cand}$, $\mbox{AMSE}(w_{\rm C, 1}) \gets \widehat{\mbox{AMSE}}\big(\widehat\beta_{\rm C}(\alpha_{ \rm opt} ;w_{\rm cand}); \beta\big)$
             }
         }
         \item Update $s_{\mathcal V} = s_{\mathcal V} + 1$.
         \end{enumerate}
         }
         {\KwRet{$w_{C, 1}$, $\widehat\beta(\alpha_{\rm opt}, w_{\rm C, 1})$, and  $\widehat{\mbox{AMSE}}(\widehat\beta_{\rm C}(\alpha_{\rm opt}; w_{\rm cand}); \beta)$
         }
         }
         }
 \end{algorithm}

\section{Numerical results}\label{sec:numerical}

\subsection{Simulation study}
In this section, we consider the following setup under the high-dimensional linear model setting.

\begin{enumerate}
    \item Fix the dimension $p = 500$, the sample size $n = 250$, the ratio $\delta = 0.5$. The number of non-zero components $s$ is taken to be 5 for the high-sparsity setting and $50$ for the medium-sparsity setting; the non-zero part is generated from the Dirac distribution with a point mass equally distributed on -1 and 1, or a standard normal distribution.
    \item In each repetition, we generate a new dataset by randomly generating a sensing matrix $X$, a coefficient vector $\beta$, and an error vector $\bm{\varepsilon}$. The components of the sensing matrix $X$ are independent and generated from $N(0, 1/250)$.
    \item  As error distributions, we take the standard normal $N(0, 1)$, student-$t$ with degrees of freedom 3, and the mixture of normal distributions $0.5N(0, 1) + 0.5 N(5, 9)$; errors generated in Step 2 are centered and rescaled to have standard deviation 0.2.
\end{enumerate}

The objective is to compare the performance of the regularized model-averaged estimator and the composite estimator with different weights, with emphasis on the weights where the selection uncertainty is taken into account.
The simulation is repeated to get 500 estimates for each setup.
For both the regularized model-averaged and composite quantile estimator, the weights considered are (1) the estimated AMSE-type weights (i.e. $w_{\rm MA, 1}$ for the model-averaged quantile estimator and $w_{\rm C, 1}$ for the composite quantile estimator), (2) the estimated weights based on minimising the asymptotic variance of the estimators of only the active set of coefficients, denoted by $w_{\rm MA, 2}$ \citep{BloznelisClaeskensZhou2019} and $w_{\rm C, 2}$ \citep{bradic2011penalized} where, with the $(k_1, k_2)$th component of $A$ equal to $A_{k_1, k_2} = \min(\tau_{k_1},\tau_{k_2})\{1 - \max(\tau_{k_1},\tau_{k_2})\}$,
$A_\varepsilon= \textrm{diag}( f_\varepsilon(u_{\tau_1}), \ldots, f_\varepsilon(u_{\tau_K}))$, and $a_\varepsilon  = ( f_\varepsilon(u_{\tau_1}), \ldots, f_\varepsilon(u_{\tau_K}))^\top$
\begin{equation}\label{eq:blonznelis_opt_weight}
    w_{\rm MA, 2} = \arg\min_{w, \mathbf 1_K^\top w = 1, w_k \geq 0} \Big\{w^\top A_\varepsilon^{-1} A A_\varepsilon^{-1} w \Big\}
\end{equation}
and
\begin{equation*}\label{eq:bradic2011_opt_weight}
    w_{\rm C, 2} = \arg\min_{w, a_\varepsilon^\top w = 1, w_k \geq 0} \Big[w^\top A w \Big].
\end{equation*}
Only considering the variance has been the standard practice so far.
(3) Equal weights $1/K$ for each component.

The number of quantiles $K$ for both estimators is taken to be 3, with quantile levels $25\%, 50 \%, 75\%$.

We present the empirical MSEs of the abovementioned estimators for estimation of three vectors of coefficients.
First, we consider the estimator of the subvector of the full coefficient   that consists of only the non-zero true coefficients, we refer to this as the ``non-zero part".
Second, we consider the estimator of the subvector of the coefficients that are truly zero. This is referred to as the ``zero part".
Third, we consider the full vector of estimated coefficients. Note that some truly zero coefficients might have a non-zero estimate, while some truly non-zero coefficients might be estimated as zero.
For each of these three vectors, ``parts", we compare the estimated values with the true values to get
$$
\mbox{MSE}(\widehat\beta_{\rm part}) =  \sum_{j_{\rm part} = 1}^{p_{\rm part}} (\widehat\beta_{j_{\rm part}} - \beta_{j_{\rm part}})^2 / p_{\rm part}
$$
for the appropriate part of the full vectors.
Results for the regularized model-averaged quantile estimator with different weights are presented in Table~\ref{table:mse_simu_mod_avrg_comb}.
We observe that the model-averaged quantile estimator using the weight in (\ref{eq:maqr opt+ weight}) has lower MSEs for estimating the non-zero part of $\beta$ and for the full vector $\beta$, and this for $t_3$ and the mixture of normally distributed errors in the high-sparse case where the number of non-zero components $s = 5$. Using equal weights leads to a fair performance of the model-averaged quantile estimator, especially for estimating the all-zero part of $\beta$.
The Lasso estimator is considered as the baseline comparison, which from Table~\ref{table:mse_simu_mod_avrg_comb} seems to have a competitive performance, especially in the medium sparsity settings. However, the Lasso mostly gives over-sparse estimations, which can be observed in the top half of Table~\ref{table:rate_simu_comb} summarizing the averaged true positive (TP) and true negative (TN) recovery rates which are defined as
$$\mbox{TP (TN)} = \frac{\mbox{number of correctly identified as non-zeros (zeros)}}{\mbox{number of true non-zeros (zeros)}}$$
The Lasso has the highest TN rate consistently and mostly the lowest TP rate. Further, while increasing the standard deviation of the errors, the Lasso's overly-sparse estimation becomes clearer, i.e., Lasso gives sparser estimations and becomes all-zeros eventually. The regularized model-averaged estimator with equal weights mostly has the highest TP rate, except for the medium sparsity settings where the non-zero part of the true regression coefficient is sampled from a Dirac distribution at -1 and 1, and the errors are sampled from $N(0,1)$ or $0.5N(0,1) + 0.5N(5, 9)$. The model-averaged estimator with the weight in (\ref{eq:maqr opt+ weight}) has the second-highest TN rate consistently.

\begin{table} 
	\centering
	\bgroup
	\def\arraystretch{0.85}
	\begin{adjustbox}{max width=\textwidth}
	\begin{tabular}{c c c c c c} \Xhline{.8pt}
	$f_\varepsilon$&  part & MSE($\widehat\beta_{\widehat w_{\rm MA,1}}$) & MSE($\widehat\beta_{\widehat w_{\rm MA, 2}}$) &
    MSE($\widehat\beta_{\widehat w_{\rm eq}^{\rm MA}}$) & MSE($\widehat \beta_{\rm Lasso}$) \\ \Xhline{.8pt}
    \rowcolor{lightgray}
	  \multicolumn{6}{c}{Non-zero part of $\beta$: \textbf{Dirac distribution at -1 and 1} ($*$: $\times 10^{-2}$, $\dagger$: $\times 10^{-3}$, $\ddagger$: $\times 10^{-4}$)} \\ \Xhline{.8pt}
  \multicolumn{1}{l}{$\mathbf{s = 5}$}  & Non-zero & 0.312 &   0.299 & 0.306 & 0.480   \\
   $N(0, 1)$& Zero ($\ddagger$) & 6.812 & 6.436 & 5.276 &   0.526 \\
    & Full vec ($\dagger$) & 3.790 & 3.630 &   3.585 & 4.854   \\   \cline{1-6}
     \multirow{3}*{$t_3$} & Non-zero &   0.167 & 0.168 & 0.182 & 0.681  \\
    & Zero ($\ddagger$) & 4.051 & 3.579 & 3.041 &   0.106\\
    & Full vec ($\dagger$)& 2.078 &   2.039 & 2.121 & 6.816 \\  \cline{1-6}
     \multirow{3}*{$0.5N(0, 1) + $} & Non-zero &   0.247 & 0.355 & 0.314 & 0.412  \\
   & Zero ($\ddagger$)& 4.593 & 7.516 & 5.418 &   0.791 \\
    $0.5 N(5, 9)$ & Full vec ($\dagger$) &   2.920 & 4.294 & 3.680  & 4.207 \\ \Xhline{.8pt}
     \multicolumn{1}{l}{$\mathbf{s = 50}$}
      & Non-zero & 0.487 & 0.502 & 0.526 &   0.376   \\
    $N(0, 1)$ & Zero ($\dagger$)& 5.498  & 4.438  &   3.675  & 5.710 \\
    & Full vec ($*$)& 5.364  & 5.419 & 5.590 &   4.275  \\   \cline{1-6}
     \multirow{3}*{$t_3$} & Non-zero & 0.399 & 0.427 & 0.452 &   0.384   \\
    & Zero ($\dagger$)& 4.976  & 3.945  &   3.412  & 5.317 \\
    & Full vec ($*$)& 4.436  & 4.630 & 4.832 &   4.318  \\   \cline{1-6}
     \multirow{3}*{$0.5N(0, 1) + $} & Non-zero & 0.504 & 0.517 & 0.540 &   0.371  \\
   & Zero ($\dagger$)& 5.303 & 4.386  &   3.635  & 5.913\\
    $0.5 N(5, 9)$ & Full vec ($*$)& 5.514 & 5.566 & 5.724  &   4.241 \\  \Xhline{.8pt}
    \rowcolor{lightgray}
    \multicolumn{6}{c}{Non-zero part of $\beta$: $\mathbf{N(0,1)}$ ($*$: $\times 10^{-2}$, $\dagger$: $\times 10^{-3}$, $\ddagger$: $\times 10^{-4}$)} \\ \Xhline{.8pt}
  \multicolumn{1}{l}{$\mathbf{s = 5}$} & Non-zero & 0.206 &   0.197 & 0.203 & 0.378 \\
   $N(0, 1)$ & Zero ($\ddagger$) & 5.624 & 5.683 & 4.439 &   0.158 \\
    & Full vec ($\dagger$) & 2.613 & 2.537 &   2.465 & 3.800   \\   \cline{1-6}
     \multirow{3}*{$t_3$} & Non-zero &   0.123 & 0.126 & 0.132 & 0.540  \\
    & Zero ($\ddagger$) & 3.727 & 3.153 & 2.752 &   0.017\\
    & Full vec ($\dagger$)& 1.601 &   1.574 & 1.590 & 5.403 \\  \cline{1-6}
    \multirow{3}*{$0.5N(0, 1) + $} & Non-zero &   0.159 & 0.230 & 0.204 & 0.313  \\
   & Zero ($\ddagger$)& 3.788 & 6.723 & 4.720 &   0.348 \\
    $0.5 N(5, 9)$ & Full vec ($\dagger$) &   1.969 & 2.970 & 2.511  & 3.162 \\  \hline
    \multicolumn{1}{l}{$\mathbf{s = 50}$} & Non-zero & 0.257 & 0.256 & 0.265 &   0.216   \\
    $N(0, 1)$ & Zero ($\dagger$)& 3.377  & 2.835  &   2.401  & 2.445 \\
    & Full vec ($*$)& 2.870  & 2.819 & 2.870 &   2.376  \\   \cline{1-6}
     \multirow{3}*{$t_3$} & Non-zero &   0.201 & 0.207 & 0.216 & 0.244   \\
    & Zero ($\dagger$)& 2.831  & 2.275  & 2.009  &   1.859 \\
    & Full vec ($*$)&   2.264  & 2.278 & 2.336 & 2.611  \\   \cline{1-6}
     \multirow{3}*{$0.5N(0, 1) + $} & Non-zero & 0.275 & 0.278 & 0.285 &   0.220  \\
   & Zero ($\dagger$)& 3.571 & 2.945  &2.536  &    2.530\\
    $0.5 N(5, 9)$ & Full vec ($*$)& 3.076 & 3.049 & 3.077  &   2.423 \\  \Xhline{.8pt}
 		\end{tabular}
	\end{adjustbox}
	\caption{The mean, over 500 simulation repetitions, of the empirical MSE of the regularized \emph{model-averaged} quantile estimator with $K=3$ for three error distributions. Empirical MSEs are calculated for the non-zero parts, all-zero parts, and the full vector of the true coefficient $\beta$. The non-zero part of the true coefficient vector is generated from Dirac distribution with point mass equally distributed on -1 and 1 (top half), or standard normal distribution (bottom half). Smaller values of MSE among competitors indicate more accurate estimations.
	}
	\label{table:mse_simu_mod_avrg_comb}
	\egroup
 \end{table}

\begin{table} 
	\centering	
	\bgroup
	\def\arraystretch{1.2}
	\begin{adjustbox}{angle=90, max width=0.87\linewidth}
	\begin{tabular}{c c| c c c c c | c c c c c}\Xhline{.8pt}
	  \multicolumn{2}{c|}{Non-zero part of $\beta$:}&  \multicolumn{5}{c|}{\bf{Dirac distribution at -1 and 1}} &
	  \multicolumn{5}{c}{$\mathbf{N(0, 1)}$} \\ \Xhline{.8pt}
	  \rowcolor{lightgray}
     $f_\varepsilon$&  rate & $\widehat\beta_{\widehat w_{\rm MA,1}}$ & $\widehat\beta_{\widehat w_{\rm MA, 2}}$ &
    $\widehat\beta_{\widehat w_{\rm eq}^{\rm MA}}$ & $\widehat \beta_{\rm Lasso}$ & $\widehat\beta_{0.5}$ &$\widehat\beta_{\widehat w_{\rm MA,1}}$ & $\widehat\beta_{\widehat w_{\rm MA, 2}}$ &
    $\widehat\beta_{\widehat w_{\rm eq}^{\rm MA}}$ & $\widehat \beta_{\rm Lasso}$ & $\widehat\beta_{0.5}$ \\ \Xhline{.8pt}
     $\mathbf{s = 5}$ \multirow{2}*{$N(0, 1)$}& TP &    0.992 & 0.991 &   0.993 & 0.906 & 0.982& 0.677 &    0.683 &   0.688 & 0.419 & 0.660 \\
     & TN & 0.904 & 0.903 & 0.896 &   0.995 &   0.940 &  0.916 & 0.912 & 0.907 &   0.998 &    0.945\\   \cline{1-10}
     \multirow{2}*{$t_3$} & TP &    0.999 & 0.999 &    0.999 & 0.663 &   1.000 & 0.754 &    0.762 &   0.765 & 0.294 & 0.739\\
    & TN & 0.910 & 0.903 & 0.896 &   0.999 &    0.941& 0.913 & 0.905 & 0.899 &   1.000 &    0.943 \\ \cline{1-10}
    $0.5N(0, 1) + $ & TP &    0.992 & 0.984 &   0.992 & 0.942 & 0.820&    0.719 & 0.711&   0.724 & 0.486 & 0.482\\
   $0.5 N(5, 9)$ & TN & 0.922 & 0.912 & 0.906 &   0.992 &     0.942& 0.927 & 0.916 & 0.911 &   0.997 &    0.946\\ \Xhline{.8pt}
    $\mathbf{s = 50}$\multirow{2}*{$N(0, 1)$} & TP & 0.836 & 0.847 &    0.854 &   0.889 & 0.548 & 0.647&    0.658 &   0.666 & 0.619 & 0.600\\
    & TN &    0.843  & 0.830  & 0.823  &   0.868 & 0.606 & 0.843 & 0.832 & 0.807 &    0.883 &   0.894\\   \cline{1-10}
   \multirow{2}*{$t_3$} & TP & 0.892 &    0.899 &   0.904 & 0.882 & 0.453 & 0.696 &     0.706 &    0.715 & 0.590 & 0.622\\
    & TN &    0.833  & 0.816  & 0.807  &   0.873 & 0.707 & 0.839 & 0.822 & 0.811 &    0.931 &    0.842\\  \cline{1-10}
     $0.5N(0, 1) + $ & TP & 0.822 & 0.837 &    0.843 &   0.892 &0.531 & 0.633 &    0.643 &   0.650 & 0.621 & 0.539\\
    $0.5 N(5, 9)$ & TN &    0.845  & 0.833 & 0.826  &   0.864 & 0.601&    0.846 & 0.834 & 0.827 &   0.911 & 0.834\\ \Xhline{.8pt}
    \rowcolor{lightgray}
     $f_\varepsilon$&  rate & $\widehat\beta_{\widehat w_{\rm C,1}}$ & $\widehat\beta_{\widehat w_{\rm C, 2}}$ &
    $\widehat\beta_{\widehat w_{\rm eq}^{\rm C}}$ & $\widehat \beta_{\rm Lasso}$ & $\widehat\beta_{0.5}$ & $\widehat\beta_{\widehat w_{\rm C,1}}$ & $\widehat\beta_{\widehat w_{\rm C, 2}}$ &
    $\widehat\beta_{\widehat w_{\rm eq}^{\rm C}}$ & $\widehat \beta_{\rm Lasso}$ & $\widehat\beta_{0.5}$\\ \Xhline{.8pt}
   $\mathbf{s = 5}$ \multirow{2}*{N(0, 1)} & TP &   0.993 &     0.991 & 0.991 & 0.911 & 0.982&   0.664 & 0.656 &    0.657 & 0.444 & 0.660\\
   & TN &    0.946 & 0.946 & 0.946 &   0.994 & 0.940&    0.963 & 0.963 & 0.963 &   0.998 & 0.945\\   \cline{1-10}
     \multirow{2}*{$t_3$} & TP &    1.000 &    1.000 & 1.000 & 0.675 & 1.000 &   0.740  & 0.729 &    0.736 & 0.303 & 0.739\\
    & TN &    0.945 & 0.944 & 0.944 &   0.998 & 0.941&    0.957 & 0.957 & 0.956 &   1.000 & 0.943\\  \cline{1-10}
    $0.5N(0, 1) + $ & TP &   0.990 &    0.982 & 0.968 & 0.939 & 0.820 &   0.699 &    0.650& 0.626 & 0.485& 0.482\\
   $0.5 N(5, 9)$ & TN &    0.955 & 0.953 & 0.951 &   0.991 & 0.942& 0.965 & 0.966 &    0.967 &   0.993 & 0.946 \\  \Xhline{.8pt}
   $\mathbf{s = 50}$ \multirow{2}*{$N(0, 1)$} & TP &   0.903 &     0.895 & 0.891 & 0.883 & 0.548 &  0.669 &    0.668 & 0.665 & 0.613 & 0.600\\
    & TN &    0.824  & 0.823  & 0.824  &   0.872 & 0.606 & 0.848 & 0.846 & 0.846 &   0.916 &    0.894\\  \cline{1-10}
     \multirow{2}*{$t_3$} & TP &   0.939 & 0.933 &    0.934 & 0.871& 0.453&   0.724 &     0.722 &  0.720 & 0.590 & 0.622 \\
    & TN &0.819  & 0.817  &     0.820  &   0.879 &0.707 & 0.835 & 0.834 & 0.835 &    0.929 &    0.842\\   \cline{1-10}
     $0.5N(0, 1) + $ & TP &    0.886 & 0.881 & 0.875 &   0.891 & 0.531  &   0.651 &    0.648 & 0.642 & 0.615 & 0.539\\
    $0.5 N(5, 9)$ & TN &    0.827  & 0.824 & 0.825  &   0.867 & 0.601 &    0.855 & 0.852 & 0.852 &   0.914 & 0.834 \\  \Xhline{.8pt}
 		\end{tabular}
 	\end{adjustbox}
	\caption{The mean, over 500 simulation repetitions, of the true positive (TP) and true negative (TN) rate of the regularized \emph{model-averaged} (top half) and \emph{composite} (bottom half) quantile estimator with $K=3$ for three error distributions. The TP and TN rates of the regularized single quantile estimator at quantile level 0.5 are presented in the 7th and 12th columns. The non-zero part of the true coefficient vector is generated from Dirac distribution with point mass equally distributed on -1 and 1 (left), or standard normal distribution (right). Larger values of TP and TN indicate a better identification power; the largest values among competitors are highlighted in green, whereas the second largest values are highlighted in yellow.}
	\label{table:rate_simu_comb}
	\egroup
 \end{table}

Since there is no analytical expression for the selection incorporated weight of the regularized composite quantile estimator $w_{\rm C, 1}$, the choice of weights can only be determined numerically by an exhaustive search. To reduce the searching time of the composite quantile estimator, we set the stopping criterion $S_{\mathcal V}$ to be five and only randomly select 4 points in the neighborhood $\mathcal{V}_{w_{\rm cand}}$; the tuning parameter $\alpha$ of the soft-thresholding function is tuned once for the regularized composite quantile estimator with the weight $w_{\rm C, 2}$, then fixed after that.

Table~\ref{table:mse_simu_comp_comb} summarizes the empirical MSEs of the regularized composite quantile estimator with different weights. Since the tuning parameter, $\alpha$ is selected for $w_{\rm C, 2}$ and a fixed tuning parameter is used for obtaining the regularized composite quantile estimates with other weights, it is not surprising that using $w_{\rm C, 2}$ leads to lower MSEs in most cases. However, it is worth noticing that using equal weights, while $\alpha$ is not optimally tuned, leads to the regularized composite quantile estimator's fair performances. The Lasso estimator consistently has the lowest empirical MSEs recovering the all-zero parts, through the largest empirical MSEs recovering the non-zero parts. This is caused by overly sparse estimations of the Lasso, which is indicated in the bottom half of Table~\ref{table:rate_simu_comb}. The regularized composite estimator with locally optimized $w_{\rm C, 1}$ consistently has the highest TP rate, and second-highest TN rate among all competitors, except the TN rate for $t_3$ distributed errors and TP rate for $0.5N(0,1) + 0.5 N(5, 9)$ distributed errors. At the same time, the non-zero parts of $\beta$ are generated from Dirac distribution at -1 and 1.

\begin{table} 
	\centering
	\bgroup
	\def\arraystretch{0.85}
	\begin{adjustbox}{max width=\textwidth}
	\begin{tabular}{c c c c c c c} \Xhline{.8pt}
	 $f_\varepsilon$&  part & MSE($\widehat\beta_{\widehat w_{\rm C,1}}$) & MSE($\widehat\beta_{\widehat w_{\rm C, 2}}$) &
    MSE($\widehat\beta_{\widehat w_{\rm eq}^{\rm C}}$) & MSE($\widehat \beta_{\rm Lasso}$) & MSE($\widehat \beta_{0.5}$)\\ \Xhline{.8pt}
    \rowcolor{lightgray}
    \multicolumn{7}{c}{Non-zero part of $\beta$: \textbf{Dirac distribution at -1 and 1} ($*$: $\times 10^{-2}$, $\dagger$: $\times 10^{-3}$, $\ddagger$: $\times 10^{-4}$)} \\ \Xhline{.8pt}
   \multicolumn{1}{l}{$\mathbf{s = 5}$}
    & Non-zero &   0.226 & 0.246 & 0.249 & 0.479 & 0.272  \\
   $N(0, 1)$ & Zero ($\ddagger$) & 6.566 & 7.119 & 7.199 &   0.571 & 11.641 \\
    & Full vec ($\dagger$) &   2.906 & 3.163 & 3.202 & 4.847 & 3.752 \\   \cline{1-7}
    \multirow{3}*{$t_3$} & Non-zero &   0.122 & 0.135 & 0.133 & 0.674 & 0.142\\
    & Zero ($\ddagger$) & 3.782 & 4.148 & 4.109 &   0.112 & 5.650 \\
    & Full vec ($\dagger$)&   1.593 & 1.756 & 1.740 & 6.747 & 1.911\\  \cline{1-7}
     \multirow{3}*{$0.5N(0, 1) + $} & Non-zero &   0.184 & 0.246 & 0.311 & 0.420 & 0.461\\
   & Zero ($\ddagger$)& 4.635 & 6.165 & 7.353 &   1.015 & 20.016\\
    $0.5 N(5, 9)$ & Full vec ($\dagger$) &   2.301 & 3.068 & 3.839  & 4.303 &6.011 \\  \Xhline{.8pt}
   \multicolumn{1}{l}{$\mathbf{s = 50}$} & Non-zero & 0.342 & 0.359 & 0.367 & 0.384 &   0.310  \\
    $N(0, 1)$ & Zero ($\dagger$)& 8.722  & 9.273  & 9.536  & 5.572 &   4.368\\
    & Full vec ($*$)&   4.203  & 4.423 & 4.524 & 4.339 & 4.308 \\   \cline{1-7}
     \multirow{3}*{$t_3$} & Non-zero &   0.280 & 0.294 & 0.298 & 0.398 & 0.598  \\
    & Zero ($\dagger$)& 7.026  & 7.511  & 7.618  &   5.159 & 19.415 \\
    & Full vec ($*$)&   3.429  & 3.617 & 3.663 & 4.443 & 5.709\\   \cline{1-7}
     \multirow{3}*{$0.5N(0, 1) + $} & Non-zero&  0.358 & 0.376 & 0.385 & 0.375  &   0.318\\
  & Zero ($\dagger$)& 9.099 & 9.644  & 10.007  & 5.750 &   4.301\\
    $0.5 N(5, 9)$ & Full vec ($*$)& 4.403 & 4.631 & 4.751  &   4.268 & 4.360 \\  \Xhline{.8pt}
    \rowcolor{lightgray}
    \multicolumn{7}{c}{Non-zero part of $\beta$: $\mathbf{N(0,1)}$ ($*$: $\times 10^{-2}$, $\dagger$: $\times 10^{-3}$, $\ddagger$: $\times 10^{-4}$)} \\ \Xhline{.8pt}
  \multicolumn{1}{l}{$\mathbf{s=5}$} & Non-zero &   0.157 & 0.173 & 0.175 & 0.363 & 0.177\\
   $N(0, 1)$ & Zero ($\ddagger$) & 3.782 & 4.311 & 4.327 &   0.220 & 9.528 \\
    & Full vec ($\dagger$) &   1.946 & 2.153 & 2.178 & 3.655  & 2.555 \\   \cline{1-7}
     \multirow{3}*{$t_3$} & Non-zero &   0.099 & 0.110 & 0.108 & 0.527 & 0.107 \\
    & Zero ($\ddagger$) & 2.575 & 2.861 & 2.826 &   0.043 & 5.169\\
    & Full vec ($\dagger$)&   1.245 & 1.382 & 1.360 & 5.273 & 1.492\\  \cline{1-7}
    \multirow{3}*{$0.5N(0, 1) + $} & Non-zero &   0.134 & 0.176 & 0.212 & 0.317 & 0.406 \\
   & Zero ($\ddagger$)& 2.787 & 3.772 & 4.550 &   0.476 & 23.379 \\
    $0.5 N(5, 9)$ & Full vec ($\dagger$) &   1.615 & 2.135 & 2.568  & 3.213 &4.469 \\  \Xhline{.8pt}
     \multicolumn{1}{l}{$\mathbf{s = 50}$} & Non-zero &   0.174 & 0.182 & 0.186 & 0.216  & 0.230  \\
    $N(0, 1)$ & Zero ($\dagger$)& 4.925  & 5.220  & 5.369  &   2.299 & 3.970\\
    & Full vec($*$)&   2.181 & 2.289 & 2.342 & 2.371 &2.571  \\   \cline{1-7}
     \multirow{3}*{$t_3$} & Non-zero &   0.130 & 0.138 & 0.139 & 0.236 & 0.168    \\
    & Zero ($\dagger$)& 3.849  & 4.077  & 4.152  &   1.888 & 2.887\\
    & Full vec ($*$)&   1.645  & 1.744 & 1.765 & 2.531 & 1.925 \\   \cline{1-7}
     \multirow{3}*{$0.5N(0, 1) + $} & Non-zero &   0.189 & 0.198 & 0.205 & 0.214  & 0.236\\
   & Zero ($\dagger$)& 5.149 & 5.507  & 5.738  &    2.386 & 3.984\\
    $0.5 N(5, 9)$ & Full vec ($*$)& 2.354 & 2.476 & 2.562  &   2.353 & 2.776 \\  \Xhline{.8pt}
 		\end{tabular}
	\end{adjustbox}
	\caption{The mean, over 500 simulation repetitions, of the empirical MSE of the regularized \emph{composite} quantile estimator with $K=3$ and the regularized single quantile estimator at quantile level 0.5 for three error distributions. Empirical MSEs are calculated for the non-zero parts, all-zero parts, and the full vector of the true coefficient $\beta$. The non-zero part of the true coefficient vector is generated from Dirac distribution with point mass equally distributed on -1 and 1 (top half), or standard normal distribution (bottom half). Smaller values of MSE among competitors indicate more accurate estimations.
	}
	\label{table:mse_simu_comp_comb}
	\egroup
 \end{table}

Tables~\ref{table:mse_simu_comp_comb} and \ref{table:rate_simu_comb} illustrate that the regularized composite quantile estimator mostly improves the performance of regularized single quantile estimator.
For the same simulations settings, we compare the averaged empirical MSEs, true positive and true negative rates of the regularized composite quantile estimator, see Table~\ref{table:rate_simu_comb}, column 7 and 12, and Table \ref{table:mse_simu_comp_comb}, column 7, with the single regularized quantile estimator at the median $\tau = 0.5$.
For settings where $s = 5$, the composite quantile estimator clearly dominates the single quantile estimator for all three error distributions. For settings where $s =50$, the composite estimator still mostly outperforms the single quantile estimator, except for the following cases: (1) the MSE for the non-zero and zero estimated subvector of $\beta$ in settings where errors are generated from $N(0, 1)$ and $0.5N(0,1) + 0.5 N(5, 9)$ distribution and the true non-zero subvector of $\beta$ is generated from a Dirac distribution; (2) TN rates in settings where errors are generated from $N(0, 1)$ and $t_3$ distribution and the true non-zero subvector of $\beta$ is generated from $N(0, 1)$.

The percentage of converged cases for the model-averaged and composite estimator, while setting the tolerance $\varepsilon_{{\rm tol}}$ to be $10^{-6}$ for different error distributions, are included in Table~\ref{table:converge_table_full}, where we define a estimator to have converged when the needed number of iterations was less than 50.

\begin{table} 
	\centering
	\def\arraystretch{1.2}
	\begin{tabular}{c c c c c} \Xhline{.8pt}
   $(\%)$ &  \multicolumn{2}{c}{$s = 10$}  & \multicolumn{2}{c}{$s = 50$}  \\ \Xhline{.8pt}
   $f_\varepsilon$ & model-averaged & composite & model-averaged  & composite  \\ \Xhline{.8pt}
   $N(0, 1)$ & 76 & 90 & 69 & 86 \\
   $t_3$ & 78 & 78 & 71 & 82 \\
   $0.5N(0, 1) + 0.5 N(5, 9)$ & 77 & 86 & 71 & 85 \\ \Xhline{.8pt}
 		\end{tabular}
	\caption{Percentage of converged cases of both regularized model-averaged and composite quantile estimators with the convergence tolerance $\varepsilon_{\rm tol}=10^{-6}$. The convergence percentage of the regularized model-averaged estimator is calculated by including only those cases of which all single quantile component estimates converge in less than 50 iterations.	
	\label{table:converge_table_full}}
 \end{table}

 Condition~\ref{cond:design} restricts  Algorithm~\ref{algo_single_quantile} to a special design matrix that does not allow correlations between the $X_{\cdot j}$'s.
However, since such   correlation might be present in reality, it is of interest to see if Algorithm~\ref{algo_single_quantile} is still numerically robust while Condition~\ref{cond:design} is relaxed in practice.
We consider a similar simulation setup as used before with $p=500$, the sample size $n = 250$, and $\delta = 0.5$. The number of non-zero components $s$ is taken to be 5 or 50; the non-zero components are generated from the Dirac distribution with point mass equally distributed on -1 or 1, or a standard normal distribution. In each simulation replication, a design matrix is first generated from a multivariate Gaussian distribution $N(0, \Sigma_X)$, then the components $X_{i,j}$ are centered and scaled such that the components of the rescaled matrix $X$ have sample variance $1/n$. Here, we allow for a Toeplitz covariance matrix $\Sigma_{X}$  of which its $(i, j)$th component $(\Sigma_{X})_{i,j} = \sigma_{X}^{|i - j|}, i, j = 1, \ldots, p$. We consider $\sigma_{X} = 0, 0.1, 0.3$.
\\
To investigate the effect of the correlation on the RAMP algorithm we consider the regularized single quantile estimator at quantile level 0.5. The error distribution considered is $t_3$. Table~\ref{table:correlated settings_single} records the performance of Algorithm~\ref{algo_single_quantile} with tolerance $\varepsilon_{\rm tol}=10^{-6}$ for such a correlated design matrix; the performance is evaluated by the empirical MSEs, the TP and TN rates, and the percentage of convergence.

 \begin{table} 
	\centering
	\bgroup
	\def\arraystretch{1}
	\begin{adjustbox}{max width=\textwidth}
	\begin{tabular}{c c c c c c c c} \Xhline{.8pt}
	$f_\varepsilon: t_3$& $\widehat\beta_{\widehat{w}_{\rm MA, 1}}$ & \multicolumn{3}{c}{$\mbox{MSE}{(\widehat\beta_{\rm vec})}$} & TP & TN & Convergence \%\\ \hline
 & & Non-zero & Zero & Full vec & \\ \hline
   \rowcolor{lightgray}\multicolumn{8}{c}{Non-zero part of $\beta$: \textbf{Dirac distribution at -1 and 1} ($*$: $\times 10^{-2}$, $\dagger$: $\times 10^{-3}$, $\ddagger$: $\times 10^{-4}$)} \\ \Xhline{.8pt}
     \multirow{2}*{$\sigma_X = 0$} & $s = 5$ & 0.142 & 5.650 ($\ddagger$) & 1.911 ($\dagger$) & 1.000 & 0.941 & 98 \\ 
     & $s = 50 $ & 0.598 & 19.415 ($\dagger$) & 5.709 ($*$) & 0.453 & 0.707 & 87 \\ \hline 
	\multirow{2}*{$\sigma_X = 0.1$} & $s = 5$ & 0.143 & 5.455($\ddagger$) & 1.969 ($\dagger$) & 1.000 & 0.944 & 97 \\ 
	& $s = 50$ & 0.403 & 5.170 ($\dagger$) & 4.492 ($*$) & 0.843 & 0.881 & 85  \\ \hline 
   \multirow{2}*{$\sigma_X = 0.3$} & $s = 5$ & 0.145 & 5.923($\ddagger$) & 2.037 ($\dagger$) & 1.000 & 0.940 & 87 \\ 
	& $s = 50$ & 0.461 & 5.159 ($\dagger$) & 5.074 ($*$) & 0.793 & 0.895 & 41  \\ \Xhline{.8pt} 
 \rowcolor{lightgray}\multicolumn{8}{c}{Non-zero part of $\beta$: \textbf{$\mathbf{N(0,1)}$} ($*$: $\times 10^{-2}$, $\dagger$: $\times 10^{-3}$, $\ddagger$: $\times 10^{-4}$)} \\
 \Xhline{.8pt}
   \multirow{2}*{$\sigma_X = 0$} & $s = 5$ & 0.107 &  5.169 ($\ddagger$) &  1.492 ($\dagger$) & 0.943 & 0.482 & 98 \\ 
   & $s = 50$ & 0.168 & 2.887 ($\dagger$) & 1.925 ($*$) & 0.622 & 0.842 &87  \\ \hline 
	\multirow{2}*{$\sigma_X = 0.1$} & $s = 5$ & 0.106 & 4.842 ($\ddagger$) & 1.534 ($\dagger$) & 0.741 & 0.948 & 97 \\ 
	& $s = 50$ & 0.182 & 2.875 ($\dagger$) & 2.079 ($*$) & 0.653 & 0.892 & 86  \\ \hline 
  	\multirow{2}*{$\sigma_X = 0.3$} & $s = 5$ & 0.109 & 4.636 ($\ddagger$) & 1.546 ($\dagger$) & 0.738 & 0.948 & 87 \\
  	& $s = 50$ & 0.198 & 2.857 ($\dagger$) & 2.236 ($*$) & 0.638 & 0.900 & 47  \\  \Xhline{.8pt} 
    \end{tabular}
	\end{adjustbox}
	\caption{The mean, over 500 simulation repetitions, of the empirical MSEs, the true positive (TP), the true negative (TN), and convergence percentages of the regularized model-averaged quantile estimators for $t_3$ distributed errors. Empirical MSEs are calculated for the non-zero parts, all-zero parts, and the full vector of the true coefficient $\beta$. The non-zero part of the true coefficient vector is generated from Dirac distribution with point mass equally distributed on -1 and 1 (top), or standard normal distribution (bottom).
	}
	\label{table:correlated settings_single}
	\egroup
 \end{table}

 We see from Table~\ref{table:correlated settings_single} that parameter estimation using Algorithm~\ref{algo_single_quantile} remains accurate and stable when weak correlations such as with $\sigma_X = 0.1$ exist between the $X_{\cdot j}$'s; the accuracy drops when we further increase the correlations as with $\sigma_X= 0.3$; it is worth mentioning that the convergence percentages decrease when the correlation increases. Further research concerning correlated data is worth considering.

\subsection{Data analysis}
We consider  the audio wave file of a waveshape from Octave in the R package \texttt{signal}. The dataset is a list of 3 elements; the audio wave sample is a vector of 17380 entries stored in the element ``sound", the sample rate is 22050 Hz stored in the element ``rate", and the resolution of the wave file is 16 bits recorded in the element ``bits". To alleviate the computational burden of the signal compression and reconstruction, we only consider the signal from the 6145th entry to the 8192th entry of the original sound wave signal.

\subsubsection{The preprocessing -- discrete wavelet transform}
Originated from the compressed sensing problem, the sparse linear model $Y = X\beta + \bm{\varepsilon}$ describes the image or signal compression. The $s$-sparse $p$-dimensional input signal $\beta$ is first compressed by a known sensing matrix $X \in \mathbbm R^{n\times p}$ with $n < p$; the compressed signal vector $X \beta \in \mathbbm R^{n}$ can be corrupted by the noise $\varepsilon$ with $\varepsilon_i$'s i.i.d. via transmission.
Notice that the $p$-dimensional input signal vector $\beta$ is assumed to be $s$-sparse which is usually unsatisfied by signals expressed in the standard basis.
To obtain the sparse representation of $\beta$ in practice, an intermediate stage of expressing the natural non-sparse vector $\beta^*$ in a proper orthonormal basis $ \Psi^* = (\psi_1^*, \ldots, \psi_p^*)$ is required. Examples of such an orthonormal basis include the orthonormal wavelet basis, the Fourier basis, and so forth.
To perform the discrete wavelet transform, we use the R package \texttt{wavethresh}. The collection of the coefficients at all resolution levels is used for further compression.

\subsubsection{The artificially corrupted compression}
To imitate the compressed sensing process, we process the audio wave signal vector as follows:
\begin{enumerate}
    \item Perform the Daubechies' least asymmetric wavelet transform with 8 vanishing moments using the \texttt{wd} function in the R package \texttt{wavethresh} on the original signal $\beta^*\in \mathbbm R^{2048}$ and obtain the corresponding wavelet coefficient vector $\beta \in \mathbbm R^{2047}$ with $p = 2047$.
    \item Randomly generate the sensing matrix $X$ with i.i.d components $X_{ij} \sim N(0, 1/n)$, where $n = \lfloor \delta^\prime p \rfloor$ and $\delta^\prime$ is the undersampling ratio chosen to be 0.5 here; compress the corresponding wavelet coefficients $\beta$ by computing $X \beta$.
    \item Corrupt the compressed wavelet coefficients by the error vector $\varepsilon$ with i.i.d.~components $\varepsilon_i$ having p.d.f $f_\varepsilon$;
    obtain the artificial observed signal vector $Y = X\beta + \varepsilon$. Additionally, the standard normal $N(0,1)$, student-$t$ with 3 degrees of freedom, and the bimodal mixed normal $0.5 N(0, 1 ) + 0.5  N(5, 9)$ are used as the corruption error distributions; the errors are sampled according to the distributions first, then centered and rescaled to have standard deviation 0.03.
\end{enumerate}
 In practice, the artificial vector $Y$ and the sensing matrix $X$ are observed. The accurate recovery of the original wavelet coefficient vector $\beta$ is of practical interest. To obtain an impression on the performance of the AMSE-type optimal weight, we generate the sensing matrix $X$ under a fixed seed number, which is set to be 1 in our case, then generate the error vector $\bm\varepsilon$ under various seed numbers. However, we only present the reconstructions under one seed for each setting in Section~\ref{section:audio_recovery} due to limited space.

 \subsubsection{Signal recovery}\label{section:audio_recovery}

\begin{figure}[!t]
    \centering
    \includegraphics[width = \linewidth]{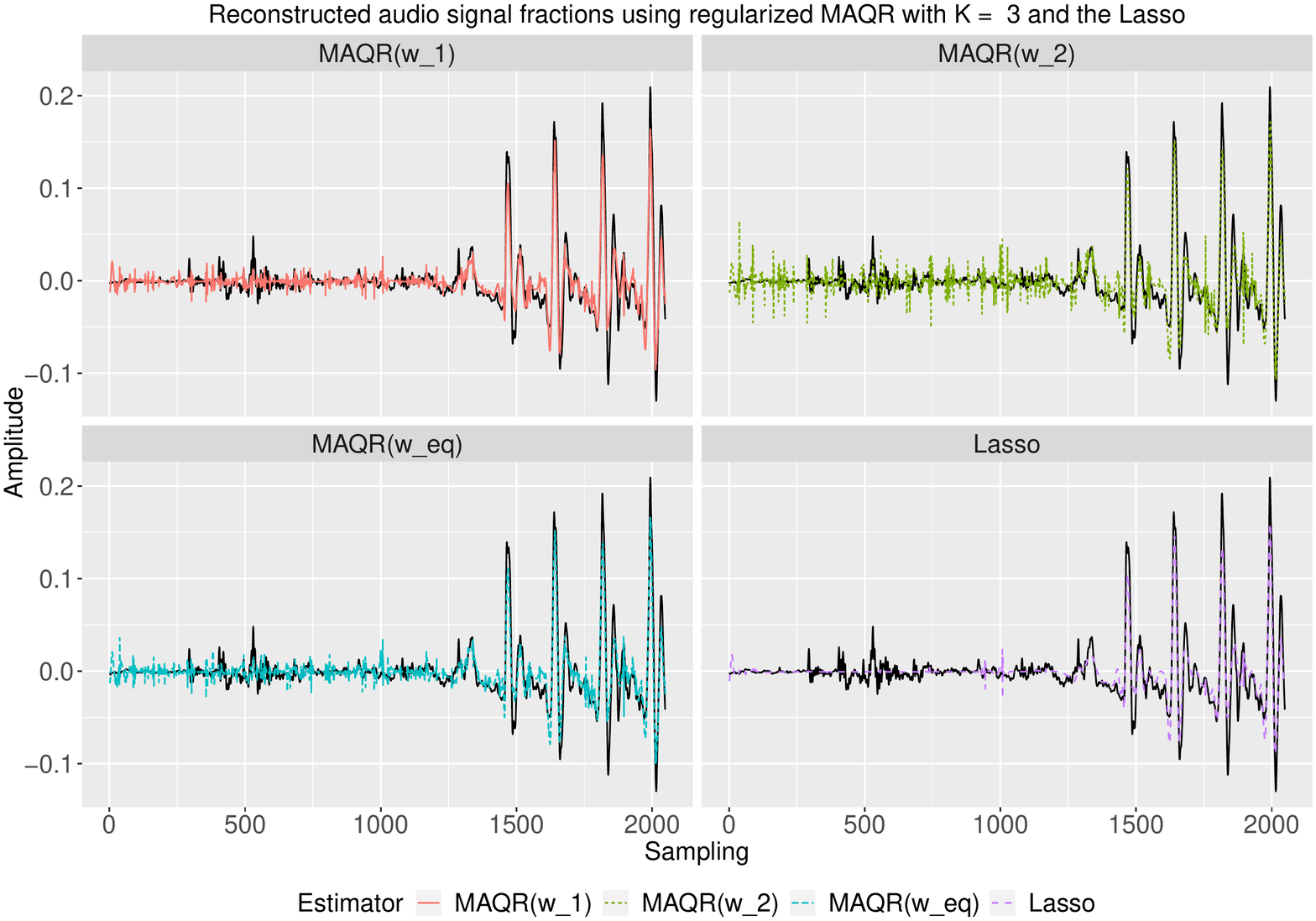}
    \caption{Reconstructed audio signal from  using the regularized \emph{model-averaged} estimator with the estimated AMSE-type weights in (\ref{eq:maqr opt+ weight}), oracle-type optimal weights in (\ref{eq:blonznelis_opt_weight}), and equal weights. The original audio curve is depicted in black. The Lasso reconstruction is presented at bottom-right. The error used for corruption follows the mixture of normals distribution $0.5 N(0,1) + 0.5N(5, 9)$.}
    \label{fig:audio_maqr_K3_mixnorm}
\end{figure}
\begin{figure}[!t]
    \centering
    \includegraphics[width = \linewidth]{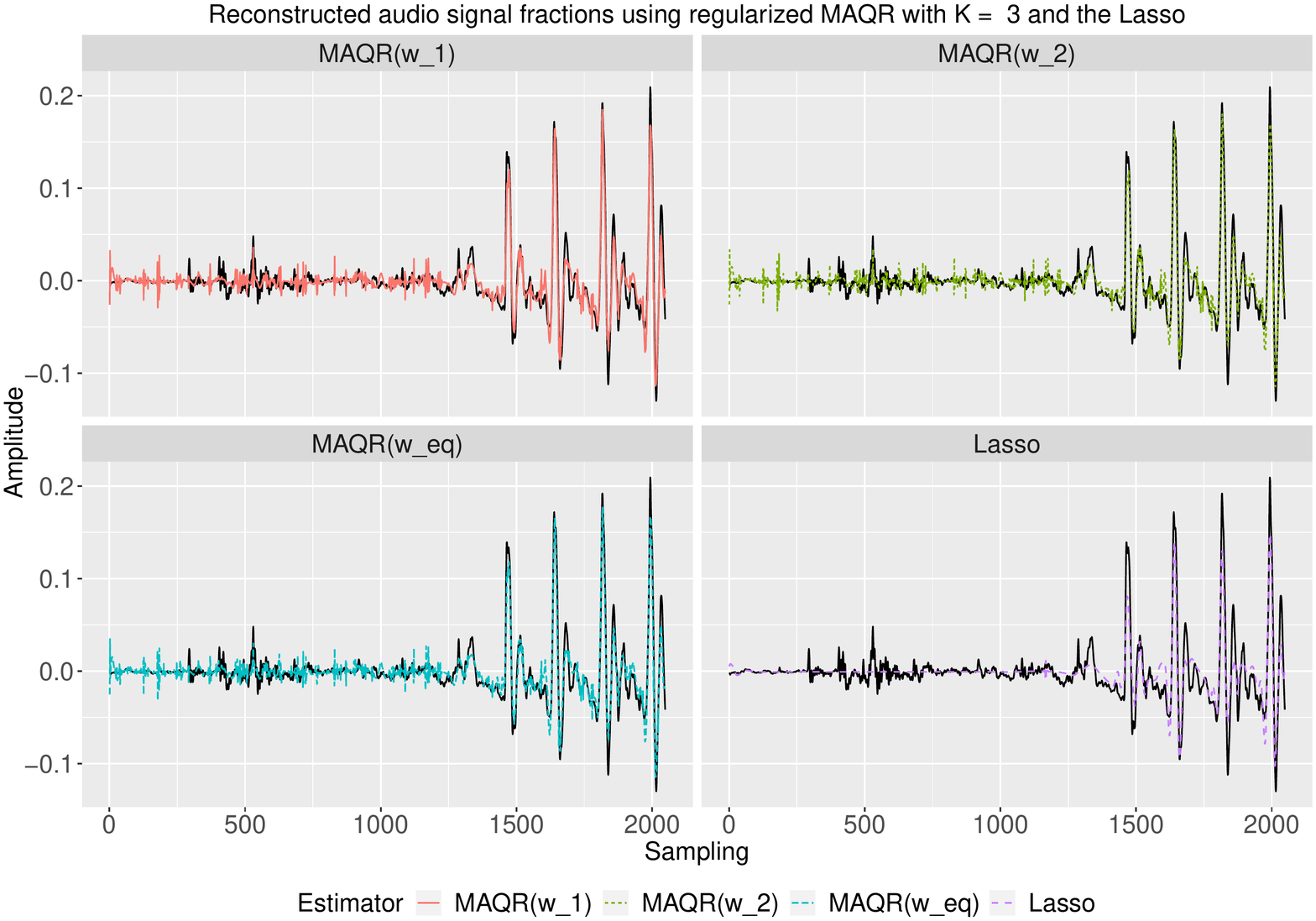}
    \caption{Reconstructed audio signal using the regularized \emph{model-averaged} estimator with the estimated AMSE-type weights in (\ref{eq:maqr opt+ weight}), oracle-type optimal weights in (\ref{eq:blonznelis_opt_weight}), and equal weights. The original audio curve is depicted in black. The Lasso reconstruction is presented at bottom-right. The error used for corruption is $t_3$ distributed.}
    \label{fig:audio_maqr_K3_t3}
\end{figure}

To reconstruct the signal vector $\beta$ expressed in the wavelet basis from the sensing matrix $X$ and the observed compressed signal vector $Y$ corrupted by potentially non-Gaussian distributed error $\bm \varepsilon$, we consider the regularized model-averaged and the composite quantile estimator weighting over three equally-spaced quantiles ($25\%, 50 \%, 75\%$) using equal weights, the oracle-type weights and the new AMSE-type weights. The tolerance in the RAMP algorithm is set as $\varepsilon_{\rm tol}=10^{-8}$. The Lasso estimator is considered as the baseline comparison.
Notice that the regularized estimates $\widehat\beta_{\rm MA}$ and $\widehat\beta_{\rm C}$ after reconstruction are the representations in the wavelet domain. To compare the accuracy of the reconstruction, we perform a back-transform on the estimates and obtain the corresponding signal vectors $\widehat{\beta}_{\rm MA}^*$ and $\widehat{\beta}_{\rm C}^*$ with representations in the natural basis.

Example reconstructions of the audio signal for $K = 3$ using the regularized model-averaged estimator equipped with different weights, with the baseline recovery from the Lasso represented in the natural basis are presented in
Figure~\ref{fig:audio_maqr_K3_mixnorm} for the mixture of normals distributed error, and in Figure~\ref{fig:audio_maqr_K3_t3} for the $t_3$ distributed error. We observe that the strong signals corresponding to large values located at the end of the sound signal are well captured by the model-averaged quantile estimator using different weights for both error distributions. For the weak signals clustering at the front of the signal, the model-averaged estimators using $\widehat{w}_{{\rm MA, 1}}$ and equal weights outperform the counterpart with $\widehat{w}_{{\rm MA, 2}}$ for $0.5 N(0,1) + 0.5 N(5, 9)$ distributed errors; recovery differences for the weak signals of the model-averaged estimator using different weights are hardly observable for the $t_3$ distributed errors. Recovery using the Lasso is competitive to the model-averaged estimator using $w_{\rm MA, 1}$ for strong signals. However, the Lasso estimates the signals in an over-sparse way with too many zeros entries; one can observe the almost flat recovery for the weak signals for both error distributions.

\begin{figure}[!t]
    \centering
    \includegraphics[width = \linewidth]{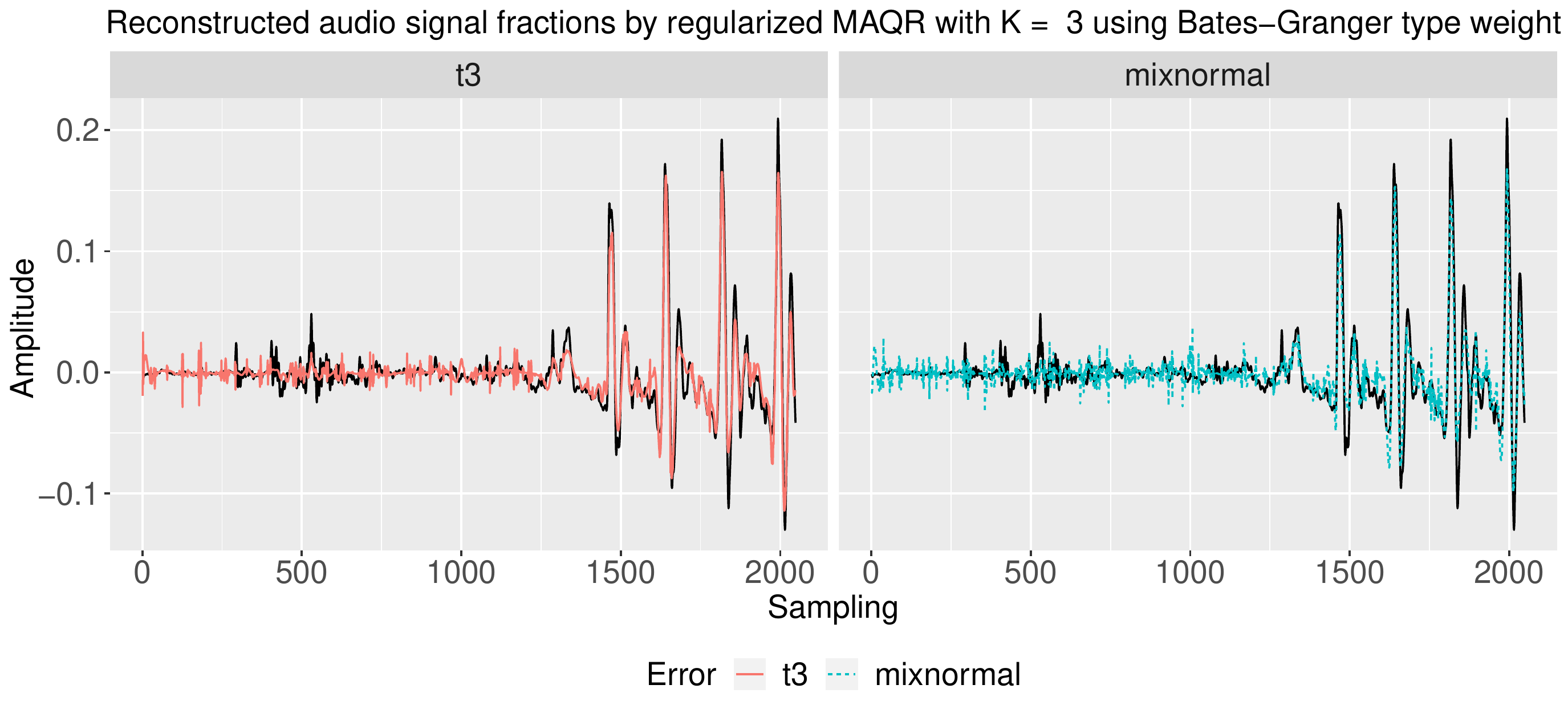}
    \caption{Reconstructed audio signal using the regularized \emph{model-averaged} estimator with Bates-Granger type weight in \eqref{eq:bates_granger weight}. The original audio curve is depicted in black. The left figure uses $t_3$ distributed corruption error, while the figure on the right used $0.5 N(0,1) + 0.5N(5, 9)$ distributed corruption error.}
    \label{fig:audio_maqr_K3_Bates}
\end{figure}

\citet{BatesGranger1969} provide an alternative weight choice for the model-averaged estimator obtained by considering only the variances of $\widehat\beta_k$'s and ignoring the covariances. This leads to
\begin{equation}\label{eq:bates_granger weight}
    \widehat{w}_{\rm MA, 3} = \arg\min_{w\geq0, \mathbf{1}_K^\top w = 1}w^\top \mbox{diag}(\widehat\Sigma_{0,(t)}) w,
\end{equation}
where $\mbox{diag}(\widehat\Sigma_{0,(t)})$ denotes the diagonal matrix obtained from $\widehat\Sigma_{0,(t)}$ which keeps the diagonal and has zeros in all off-diagonal entries.
Figure~\ref{fig:audio_maqr_K3_Bates} contains the recovery of the audio signal using the model-averaged estimator using this weight.

For the composite quantile estimator $\widehat\beta_{\rm C}$,
we performed the same weight searching method as for the simulation study. This is, $S_{\mathcal V}=5$ and randomly select 4 candidate weights in the neighbourhood of the previous value. We select the tuning parameter $\alpha$ once for the starting weight $\widehat w_{\rm C, 2}$, it remains unchanged thereafter. The recovered signals by the composite estimator with different weights are very similar in all cases.

To compare the recovery of the regularized model-averaged and composite estimator combined with different weights, as well as the Lasso estimator, we present the mean absolute percentage error (MAPE) in Table~\ref{table:recovery accuracy} where the MAPE is defined as
\begin{equation}\label{eq:mape}
    \mbox{MAPE}(\widehat\beta, \beta) = \frac{1}{p} \sum_{j = 1}^p \Big|{\big(\widehat\beta_j - \beta_j \big)} \Big/ {\beta_j} \Big|
\end{equation}
Table~\ref{table:recovery accuracy_mse} reports the MSE.
\begin{table} 
	\centering
	\bgroup
	\def\arraystretch{1.2}
	\begin{adjustbox}{max width=\textwidth}
	\begin{tabular}{c c c c c c c c c} \Xhline{.8pt}
  $f_\varepsilon$ & \multicolumn{4}{c}{$t_3$} & \multicolumn{4}{c}{0.5N(0,1) + 0.5 N(5, 9)} \\\Xhline{.8pt}
   est: MA / C & $w_{\rm est, 1}$ &  $w_{\rm est, 2}$ & $w_{\rm eq}$ & $w_{\rm est, 3}$ & $w_{\rm est, 1}$ &  $w_{\rm est, 2}$ & $w_{\rm eq}$ & $w_{\rm est, 3}$ \\\Xhline{.8pt}
    MAQR & 3.177 & 3.339 & 3.341 & 3.005 & 2.934 & 5.746 & 3.722 & 4.152 \\
    CQR & 2.798 & 2.730& 2.671 & - & 6.100 & 6.090 & 6.462 & -\\
    Lasso & \multicolumn{4}{c}{1.346} & \multicolumn{4}{c}{1.536} \\ \Xhline{.8pt}
 		\end{tabular}
	\end{adjustbox}
	\caption{The MAPE defined in (\ref{eq:mape}) of the audio signal recovered by the regularized model-averaged and composite estimators with different weights, and the Lasso estimator. The seed number used to generated the errors for corrupting the compressed signal vector is 37 for both $t_3$ and mixed normal distributed errors.
	}
	\label{table:recovery accuracy}
	\egroup
 \end{table}

 \begin{table} 
	\centering
	\bgroup
	\def\arraystretch{1.2}
	\begin{adjustbox}{max width=\textwidth}
	\begin{tabular}{c c c c c c c c c} \Xhline{.8pt}
  $ f_\varepsilon $ & \multicolumn{4}{c}{$t_3$} & \multicolumn{4}{c}{0.5N(0,1) + 0.5 N(5, 9)} \\\Xhline{.8pt}
   est: MA / C ($\times 10^{-4}$) & $w_{\rm est, 1}$ &  $w_{\rm est, 2}$ & $w_{\rm eq}$ & $w_{\rm est, 3}$ & $w_{\rm est, 1}$ &  $w_{\rm est, 2}$ & $w_{\rm eq}$ & $w_{\rm est,3}$\\\Xhline{.8pt}
    MAQR & 1.286 & 1.288 & 1.274  & 1.304 & 2.044 & 2.417 & 2.051 & 2.006 \\
    CQR & 1.279 & 1.273 & 1.271 & - & 2.363 & 2.068 & 2.120 & -\\
    Lasso & \multicolumn{4}{c}{2.566} & \multicolumn{4}{c}{2.003} \\ \Xhline{.8pt}
 		\end{tabular}
	\end{adjustbox}
	\caption{The MSE of the audio signal recovered by the regularized model-averaged and composite estimators with different weights, and the Lasso estimator. The seed number used to generated the errors for corrupting the compressed signal vector is 37 for both $t_3$ and mixed normal distributed errors.
	}
	\label{table:recovery accuracy_mse}
	\egroup
 \end{table}

 We see that the Lasso has the lowest MAPE for both $t_3$ and mixed normal distributed errors; at the same time, it estimates the weak signals in an over-sparse way and is not capable of capturing the weak signals. Comparing the effect of different weight choices on the regularized model-averaged quantile estimator with its composite quantile counterpart, we see that the MAPEs of the composite quantile estimators are relatively stable using different weights. The model-averaged estimator with the AMSE-type weight $\widehat w_{\rm MA, 1}$ has excellent performance compared to the composite estimator, especially for the mixed normal distributed error.
 The Bates-Granger weighting provides good results regarding MAPE for the $t_3$ error case, but not for the mixed normal. Regarding MSE, it performs well for the mixed normal case but is worst for the $t_3$ errors, wherein this example the equal weights perform best, although all results are close. Searching for the selection incorporated weight $\widehat w_{\rm C, 1}$ for the regularized composite quantile estimator is computationally infeasible for large $p$ (2047 in our case). Estimating the regularized model-averaged quantile estimator averaging three quantiles here takes approximately 4 -- 5 hours whereas estimating the regularized composite quantile estimator takes more than 16 hours with only five steps in a nearby search with four surrounding candidate weights, and the tuning parameter $\alpha$ tuned only once for the starting weight.

 Additionally, we present the estimated weights for both regularized model-averaged and composite estimators in Table~\ref{table:recovery opt weight}. An interesting observation is made by comparing the estimated weights $\widehat w_{\rm MA, 1}$ and $\widehat w_{\rm MA, 2}$ for the mixed normal distributed error. The weight $\widehat w_{\rm MA, 1}$ presented here is quite representative; it assigns weight 0 to the quantile estimate at 50\% quantile level suggesting the final model-averaged estimate is obtained by averaging estimates at 25\% and 75\% quantile levels. On the contrary, $\widehat w_{\rm MA, 2}$ assigns the largest weight to the estimate at a 50\% quantile level indicating the most significant contribution to the final model-averaged estimate.

  \begin{table} 
	\centering
	\bgroup
	\def\arraystretch{1.2}
	\begin{adjustbox}{max width=\textwidth}
	\begin{tabular}{c c c c } \Xhline{.8pt}
  $ f_\varepsilon $ & est: MA / C & MAQR & CQR \\\Xhline{.8pt}
  \multirow{2}*{$t_3$} & $w_{\rm est, 1}$ & (0.156, 0.725, 0.119) & (0.089, 0.492, 0.419) \\
  & $w_{\rm est, 2}$ & (0.077, 0.650, 0.273) &  (0.314, 0.267, 0.467) \\ \Xhline{.8pt}
  $0.5 N(0,1) +$ & $w_{\rm est, 1}$ & (0.548, 0, 0.452) & (0.469, 0.495, 0.036) \\
 $ 0.5N(5,9)$ &  $w_{\rm est, 2}$ & (0.147, 0.843, 0.010) & (0.369, 0.345, 0.286) \\ \Xhline{.8pt}
 	\end{tabular}
	\end{adjustbox}
	\caption{The estimated weights $\widehat w_{\rm MA, 1}$ and $\widehat w_{\rm MA, 2}$ for the model-averaged estimator, and $\widehat w_{\rm C, 1}$ and $\widehat w_{\rm C, 2}$ for the composite estimator. The seed number used to generated the errors for corrupting the compressed signal vector is 37 for both $t_3$ and mixed normal distributed errors.
	}
	\label{table:recovery opt weight}
	\egroup
 \end{table}

\section{Discussion} \label{sec:discussion}

This paper is the first to take the selection uncertainty due to regularization into account when computing the weights used in model-averaged and composite estimation.
While we have studied both composite estimation and model-averaged estimation, the flexibility of allowing for parallel computation and a component-specific choice of regularization, combined with an explicit expression of the optimal weights for model averaging, places this method in a preferred position from a computational point of view.

It would be interesting to investigate whether AMSE expressions for other types of regularization may be obtained similarly.
Going yet one step further would be incorporating the effect of data-driven values of the regularization parameters $\lambda$ (for composite estimation) and $\lambda_1,\ldots,\lambda_K$ (for model-averaged estimation) on the choice of the weights.
To further study the weight selection and the effect of using data-driven weights, one should study the joint distribution of the estimated weights and the estimators of interest. To simplify such matters, sample splitting could be used such that the weights are computed on a hold-out sample and the estimation using those weights proceeds on the rest of the sample.
In this paper, we used the same dataset for estimating both $\beta$ and $w$.

To avoid overly complicated mathematical expressions, we followed earlier literature in the use of a design matrix where $X_{ij}\sim N(0, 1/n)$. Other applications might require studying, for example, fixed designs, which are beyond the scope of the current paper.

\section*{Appendix} \label{sec:appendix}
\appendix

\section{Assumptions}
\begin{enumerate}[label=(A{{\arabic*}})]
\item   \label{cond:design}\label{Afirst}
    Design: The elements of the design matrix $X$, that is $X_{ij}$ for $i=1,\ldots,p$ and $j=1,\ldots,n$, are independent and identically distributed according to a $N(0,1/n)$ which is also called a standard Gaussian design.

\item \label{cond:beta}
    Coefficients: The $p$-vector ${\beta}$ is such that the sequence of uniform distributions that is placed on its components converges, for $p$ tending to infinity, to a distribution with a bounded $(2k-2)$th moment for $k\ge 2$. Denote by $B_0$ a random variable with this limiting distribution function $F_{B_0}$.

\item \label{cond:loss}
    Loss function: (i) The subgradient $\partial\rho(u)=\sum_{j=1}^3 v_j(u)$ where $v_1$ has an absolutely continuous derivative, $v_2$ is continuous and consists of piecewise linear parts and is constant outside a bounded interval, and $v_3$ is a non-decreasing step function. Denote $v_2'(u)=\alpha_l$ and $v_3(u)=\gamma_l$ when $u\in(r_l,r_{l+1}]$ where $\alpha_0=\alpha_L=0$,  $-\infty=r_0<r_1<\ldots<r_L<r_{L+1}=\infty$ and $-\infty=\gamma_0<\gamma_1<\ldots<\gamma_L<\gamma_{L+1}=\infty$.
(ii) The subgradient's absolute value $|\partial\rho(u)|$ is bounded for all $u\in\mathbbm R$.
(iii)
$h(t)=\int \rho(z-t)dF_\varepsilon(z)$ has a unique minimum at $t=0$.
(iv)
There exists a $\delta>0$ and $\eta>1$ such that $E[\{\sup_{|u|\le \delta}|v_1''(z+u)|\}^\eta]$ is finite.

\item \label{cond:lemma_convergence}  \label{A4}
     We assume that for some $\kappa>1$,
     \begin{enumerate}
         \item $\lim_{p \to \infty} E_{\widehat{f}_{\beta}}(B_0^{2\kappa - 2}) = E_{f_{B_0}}(B_0^{2\kappa - 2}) < \infty$
         \item $\lim_{p \to \infty} E_{\widehat{f}_{\varepsilon}}(\varepsilon^{2\kappa - 2}) = E_{f_{\varepsilon}}(\varepsilon^{2\kappa - 2}) < \infty$
         \item $\lim_{p \to \infty} E_{\widehat{f}_{q_0}}(B_0^{2\kappa - 2}) < \infty$.
     \end{enumerate}

\item \label{cond:error} \label{Alast}
The regression errors $\varepsilon_1, \ldots, \varepsilon_n$ and $\varepsilon$ are i.i.d. random variables with mean zero and finite 2nd moment.
Assume $\varepsilon$ has cumulative distribution function $F_\varepsilon$ and probability density function $f_\varepsilon$. Let $F_\varepsilon$ have bounded derivatives $f_\varepsilon$ and $\partial f_\varepsilon$; further, let $f_\varepsilon > 0$ in the neighbourhood of $r_1, \ldots, r_L$ in \ref{cond:loss}.
\end{enumerate}
Assumption \ref{cond:design} has been used by \citet{bayati2011dynamics, donoho2016high, bradic2016robustness}, assumption \ref{cond:beta} has been used by \citet{bayati2011dynamics, bradic2016robustness}; while conditions \ref{cond:loss} and \ref{cond:error} correspond to conditions R and D of \citep{bradic2016robustness}.
Condition \ref{A4} is used in Lemma~\ref{lemma:converge_cov}, in addition to the moment condition stated in \ref{cond:beta} and \ref{cond:error}.
We take $\kappa = 2$ for Algorithm~\ref{algo_single_quantile}.

\section{Lemmas and Proofs}

\subsection{Auxiliary definitions and lemmas}
\theoremstyle{definition}
\newtheorem{definition}{Definition}
\begin{definition}{(Pseudo Lipschitz function)}\label{def:pseudo_lipschitz}
A function $\phi:\mathbbm R^m \to \mathbbm R$ is pseudo-Lipschitz of order $\kappa\geq 1$, if there exists a constant $L > 0$, such that $\forall x, y \in \mathbbm R^m$
$$
|\phi(x) - \phi(y)| \leq L(1 + \| x\|^{\kappa - 1} + \|y\|^{\kappa -1})\|x - y\|.
$$
\end{definition}

It follows that if $\phi$ is a pseudo-Lipschitz function of order $\kappa$, then there exists a constant $L'$ such that $\forall x\in \mathbbm R^m: |\phi(x)| \leq L' (1 + \|x\|^{\kappa})$.

\begin{lemma}[Theorem 1 in \citet{jameson2014some}]\label{eq:polynomial_ineq}
    If $x_{i} \ge 0$ where $i = 1, \ldots, n$ and $p \ge 1$,  then
    \begin{eqnarray*}
       \sum_{i = 1}^n x_{i}^p \leq (\sum_{i = 1}^n x_{i})^p \leq n^{p-1} \sum_{i = 1}^n x_{i}^p.
    \end{eqnarray*}
    The reversed inequality holds for $p \in (0, 1)$
\end{lemma}

\begin{lemma}[Extrema of quadratic forms in \citet{rao1973linear}]\label{lemma:extrema_rao}
    Let $A$ be a $m \times m$ matrix, $B$ be a $m\times k$ matrix, and $U$ be a $k$-vector. Denote by $S^-$ any generalized inverse of $B^\top A^{-1} B$, then
    \begin{eqnarray*}
       \inf_{B^\top X =U} X^\top A X = U^\top S^- U
    \end{eqnarray*}
where $X$ is a column vector and the infimum is attained at $A^{-1}B S^- U$.
\end{lemma}

\begin{lemma}[Stein's lemma in \citet{stein1981estimation}]\label{lemma:stein's lemma}
    Let $X_{1}, X_{2}$ jointly Gaussian distributed. Let $g:\mathbbm R \to \mathbbm R$ be absolutely continuous with derivative $\partial g$ and $E|\partial g(X_{1})| \allowbreak < \infty$. Then
    $$\Cov(g(X_1), X_2) = \Cov(X_1, X_2)E[\partial g(X_1)].$$
\end{lemma}

\begin{lemma}[Lemma 4 in \citet{bayati2011dynamics}]\label{lemma:lemma_1_coef_converge}
    Let $\kappa \ge 2$ and a sequence of vectors $\{\beta(p)\}_{p \ge 0}$ whose empirical distribution converges weakly to probability measure $f_{B_0}$ on $\mathbbm R$ with bounded $\kappa$th moment; additionally, assume that $\lim_{p \to \infty} E_{\widehat f_{\beta}}(B_0^\kappa)  = E_{f_{B_0}}(B_0^\kappa)$. Then for any pseudo-Lipschitz function $\psi: \mathbbm R \to \mathbbm R$ of order $\kappa$:
    $$\lim_{p\to \infty} \frac{1}{p} \sum_{j = 1}^p \psi(\beta_j) \stackrel{\rm a.s.}{=} E[\psi(B_0)].$$
\end{lemma}

\subsection{Proofs} \label{sec:proofs}

\subsubsection{Proof of \eqref{eq:proximal cqr}}
\label{proof:proxcqr}
\begin{proof}
By definition, the proximal mapping operator is the minimizer of the function
$
b\rho_{\rm C}(x) + 0.5(x -z)^2
$
which is non-differentiable but subdifferentiable, with subgradient
$
b \cdot \partial \rho_{\rm C}(x) + x - z.
$
$\mbox{Prox}(z; b)$ is the minimizer if and only if
$
0 \in \{b \cdot \partial \rho_{\rm C}(x)|_{x = \rm{Prox}(z; b)}+ \mbox{Prox}(z; b) - z\}.
$
We distinguish between intervals where $\rho_{\rm C}$ is differentiable and non-differentiable points.
For $x \in (u_{\tau_{\ell}}, u_{\tau_{\ell + 1}}), \ell = 0, \ldots, K$ the function $\rho_C$ is differentiable.
Using the expression of the subgradient in \eqref{eq:subgradient_cqr}, we obtain
     $ 0 =b h(\ell) + x - z$, which is solved for $x$ to get that $\mbox{Prox}(z; b) = z - bh(\ell)$.
    From $\mbox{Prox}(z; b) \in (u_{\tau_\ell}, u_{\tau_{\ell+1}})$ it follows that $
     z \in (u_{\tau_\ell} + bh(\ell), u_{\tau_{\ell+1}} + bh(\ell))$.
For the non-differentiable points, that is $x = u_{\tau_\ell}, \ell = 1, \ldots, K$, having
     $ 0 \in \{ b [h(\ell - 1), h(\ell)]  + x - z\}$  leads to $u_{\tau_\ell} =\mbox{Prox}(z; b) \in [z - b h(\ell), z - b h(\ell-1)].$
     This implies that $
     z \in [u_{\tau_\ell} + b h(\ell - 1), b_{\tau_\ell} + b h(\ell)]$.
\end{proof}

\subsubsection{Proof of \eqref{eq:CQR_effective_score}}
\label{Proof:effscore}
\begin{proof}
By definition, $\widetilde G(z;b) = b \cdot \partial\rho(x)|_{x = \text{Prox}(z;b)}$, and, see the Proof of \eqref{eq:proximal cqr} in Section~\ref{proof:proxcqr}, $0 \in \{b \cdot \partial\rho(x)|_{x = \text{Prox}(z;b)}+ \mbox{Prox}(z, b) -z.$ Without loss of generality, we show the calculation for the cases where $z < u_{\tau_1} + b h(0)$ and where $z \in [u_{\tau_1} + bh(0), u_{\tau_1} + b h(1)]$.

For $z < u_{\tau_1} + bh(0)$ it holds that $\mbox{Prox}(z;b) = z -bh(0) < u_{\tau_1}$, which leads to
 $\partial\rho(x)|_{x = \rm{Prox}(z;b)} = h(0)$.
  Hence,
   $\widetilde G(z;b) = b \cdot \partial\rho(x)|_{x = \rm{Prox}(z;b)}  =  b h(0)$.

Having $z \in [u_{\tau_1} + bh(0), u_{\tau_1} + b h(1)]$ corresponds to taking the nondifferentiable point $u_{\tau_1}=\mbox{Prox}(z; b)$, see the proof of \eqref{eq:proximal cqr}.
We have  $\partial\rho(x)|_{x = \text{Prox}(z;b)}\in [h(0), h(1)]$. The subgradient $\partial\rho_{\rm C}$ is non-decreasing (Condition \ref{cond:loss}) and linear. From \eqref{eq:proximal cqr} the proximal operator is also a linear function. An intuitive choice for $\widetilde G(z;b) = b \cdot \partial\rho(x)|_{x = Prox(z;b)}$ with $z \in [u_{\tau_1}+bh(0), u_{\tau_1} + b h(1)]$ is $\widetilde G(z;b) = z - b_{\tau_1}$ which keeps the linearity of the composition of the two functions $\partial \rho$ and $\mbox{Prox}(\cdot; b)$.
\end{proof}

\subsubsection{Proof of \eqref{eq:converge_eff_score}}
\label{Proof:ofconv_eff_score}

\begin{proof}
Theorem 2, Eq.(3.7) of \citet{bayati2011dynamics} states in our notation that
 \begin{equation}\label{eq:eq3.7 in bayati and Montanari}
  \lim_{n\to \infty} \frac{1}{n} \sum_{i = 1}^n \psi(\varepsilon_i-z_{(t), i}, \varepsilon_i) \stackrel{a.s.}{=}
E[\psi(\bar\sigma_{(t)} Z, \varepsilon)],
\end{equation}
where $\psi: \mathbbm{R}^2 \to \mathbbm{R}$ is any pseudo-Lipschitz function, $Z \sim N(0,1)$,
$\bar\sigma_{(t)}$ from \eqref{eq:state_evolution_as_sigma},
and $\varepsilon$ as in \ref{Alast}.
Motivated by Eqs. (7.16) and (7.18) in \citet{bradic2016robustness} we take
$
    \psi(d, \varepsilon) = \{G(\varepsilon - d; b_{(t)})\}^2,
$
with $b_{(t)}$ as in Algorithm~\ref{algo_single_quantile}, step 2.
Applying \eqref{eq:eq3.7 in bayati and Montanari}  we obtain
that as $n\to\infty$
\begin{eqnarray*}
 \frac{1}{n} \sum_{i=1}^n G(\varepsilon_i - (\varepsilon_i-z_{i, (t)}); b_{(t)})^2
 =\frac{1}{n}\sum_{i = 1}^n G(z_i; b_{(t)})^2
 &\stackrel{a.s.}{\to}& E[G(\varepsilon - \bar\sigma_{(t)} Z; b_{(t)})^2].
\end{eqnarray*}
\end{proof}

\subsubsection{Estimation of $\nu(b)$} \label{section:estimator_nu}
 The effective score step in Section~\ref{section:algo}, in cases where $G(\cdot; b_{(t)})$ is non-differentiable, requires a solution $b_{(t)}$ to the equation $1 = \widehat\nu(b_{(t)})$ where $\widehat\nu(b_{(t)})$ is a consistent estimator of a population parameter $\nu(b_{(t)})$ defined as
 $$\nu(b_{(t)}) =  E[\partial_1\widetilde G(C_{(t)}; b_{(t)})] = b_{(t)}(\delta/\omega) E( \partial[ \partial \rho\{ \mbox{Prox}(C_{(t)};b_{(t)})\} ])$$
 with $C_{(t)} = \varepsilon - \bar\sigma_{(t)}Z$ the random variable characterizing the limit distribution of the adjusted residuals $z_{(t)}$ when $p \to \infty$.
\\
Using Condition~\ref{cond:loss} and Lemma~3 of \citet{bradic2016robustness}, $\partial\rho$ can be written as a sum of three functions of which $v_1$ and $v_2$ are differentiable. For the step function $v_3$, we use Condition~\ref{cond:loss} on $\rho$, where $\gamma_l$ is the step height on the interval $(r_{l}, r_{l+1}]$. Let $f_{C_{(t)}-\widetilde G(C_{(t)};b)}$ denote the density of the variable $C_{(t)} - \widetilde G(C_{(t)};b)$ which is equivalent to $\mbox{Prox}(C_{(t)}; b)$. The equivalence is obtained by setting the derivative of the $b\rho(x) + \frac{1}{2}(x - C_{(t)})^2$ w.r.t. $x$ to zero and evaluate at $\mbox{Prox}(C_{(t)}; b)$, due to the fact that the proximal operator is the minimizer of the function $b\rho(x) + \frac{1}{2}(x - C_{(t)})^2$.
Then we arrive at
\begin{eqnarray*}
\frac{\omega \nu(b_{(t)})}{\delta b_{(t)}} = \sum_{j=1}^2
E[\partial v_j (C_{(t)})] +
\sum_{l=1}^{L-1}\gamma_l  \{f_{C_{(t)}-\widetilde G(C_{(t)};b_{(t)})}(r_{l+1}) -
f_{C_{(t)}-\widetilde G(C_{(t)};b)}(r_{l})\}.
\end{eqnarray*}
The consistent estimator in (\ref{eq:estimator_nu}) is obtained by replacing the expectation above with the empirical mean and replacing the density of the proximal operator $\mbox{Prox}(C_{(t)}; b)$ with its kernel density estimator.

\subsubsection{Proof of Lemma~\ref{lemma:converge_cov}}
Since this proof is based on the general recursion and Lemma 1 in \citet{bayati2011dynamics}, we first restate the general recursion to which Algorithm~\ref{algo_single_quantile} belongs with slight changes in the notations.
Given the noise $\bm{\varepsilon} \in \mathbbm R^n$ and the coefficient vector $\beta \in \mathbbm R^p$, the general recursion is defined
$$h_{(t+1)} = X^\top m_{(t)} - \xi_{1, (t)}q_{(t)},\ \ \ m_{(t)} = g_{1, t}(d_{(t)}, \varepsilon)$$
$$d_{(t)} = X q_{(t)} - \xi_{2, (t)} m_{(t-1)},\ \ \  q_{(t)} = g_{2, (t)}(h_{(t)}, \beta) $$
where $\xi_{1, (t)} = n^{-1} \sum_{i = 1}^n \partial_1 g_{1, (t)}(d_{(t), i}, \varepsilon_i)$,
$\xi_{2, (t)} = (\delta p)^{-1}\sum_{j = 1}^p\partial_1 g_{2, (t)}(h_{(t), j}, \beta_j)$.
Further, to connect the general recursion to Algorithm~\ref{algo_single_quantile}, we also state the exact form of $h_{(t+1)}, m_{(t)}, d_{(t)}, q_{(t)}$ taken in Algorithm~\ref{algo_single_quantile}.
Lemma 1 in \citet{bradic2016robustness} states that Algorithm~\ref{algo_single_quantile} takes
$h_{(t+1)} = \beta - X^\top G(z_{(t)}; b_{(t)}) - \beta_{(t)}, \ \ \  q_{(t)} = \beta_{(t)} - \beta,$  from \eqref{eq:adjust_resid} $z_{(t)} = \varepsilon - d_{(t)}$, which defines $d_{(t)}$, $m_{(t)} = - G(z_{(t)}; b_{(t)})$
with the functions $g_{1, (t)}(x_1,x_2) = - G(x_2-x_1; b_{(t)})$, and $g_{2, (t)}(x_1) = \eta(\beta - x_1; \theta) - \beta$.
To proceed with the proof of Lemma~\ref{lemma:converge_cov}, we first recall the technique used for proving Lemma 1 in \citet{bayati2011dynamics}, which uses induction on the iteration $t$. To not fully repeat the long proof and all notations we only give details about where our proof differs from theirs.
\begin{enumerate}
    \item $\mathcal B_{(0)}$: show properties (3.15), (3.17), (3.19), (3.21), (3.23) and (3.23) of \citet{bayati2011dynamics} which are related to the vectors $b_{(0)}$ and $m_{(0)}$, by conditioning on the $\sigma$-algebra $\mathcal D_{(0),(0)}$ generated by $\{\beta, \varepsilon, q_{(0)}\}$; obtain the $\sigma$-algebra $\mathcal D_{(1),(0)}$ by adding $b_{(0)}$ and $m_{(0)}$ to the set $S_{(0),(0)} = \{\beta, \varepsilon, q_{(0)}\}$.
    \item $\mathcal H_1$: show that the properties (3.14), (3.16), (3.18), (3.20), (3.22), (3.24) and (3.25), which are related to the vectors $h_{(1)}$ and $q_{(1)}$, hold by conditioning on the $\sigma$-algebra $\mathcal D_{(1),(0)}$; obtain the $\sigma$-algebra $\mathcal D_{(1),(1)}$ by adding $h_{(1)}$ and $m_{(1)}$ to the set $S_{(1), (0)} = \{\beta, \varepsilon, q_{(0)}, d_{(0)}, m_{(0)}\}$
    \item $\mathcal B_{(t)}$: Similar to $\mathcal B_{(0)}$; the proof is conditioning on the $\sigma$-algebra $\mathcal D_{(t),(t)}$ for the set containing $\beta, \varepsilon, q_{(0)}$ and all previous obtained vectors; obtain the new $\sigma$-algebra $\mathcal D_{(t+1), (t)}$ by adding $b_{(t+1)}$ and $m_{(t+1)}$ to the set.
    \item $\mathcal{H}_{(t+1)}$: Similar to $\mathcal H_{(1)}$; conditioning on the $\sigma$-algebra $\mathcal D_{(t+1),(t)}$ for the set containing $\beta, \varepsilon, q_{(0)}$ and all previous obtained vectors.
\end{enumerate}
Assuming Lemma 1 in \citet{bayati2011dynamics} holds for all $K$ estimators $\widehat\beta_{k}, k=1,\ldots,K$ in (\ref{eq:MA}), we add an additional step considering the correlations between the estimators. The main technique is conditioning on the $\sigma$-algebra generated by $\cup_{k=1}^K \mathcal S_{k, (1), (0)}$ and $\cup_{k=1}^K \mathcal S_{k, (t+1), (t)}$, where $\mathcal S_{k, (1), (0)}$ and $\mathcal S_{k, (t+1),(t)}$ are the sets described in step 2 and 4 above for the $k$th estimator. The proof is similar to that of (3.16) in Lemma 1(b) of \citet{bayati2011dynamics}, with different mathematical techniques in order to adjust the original proof from a single sequence of iterations to $K$ paralleled sequences of iterations.

\begin{proof}
Idea of the construction: The construction of $\mathcal B_{(0)}$, $\mathcal H_{(0)}$, $\mathcal B_{(t+1)}$ and $\mathcal{H}_{(t+1)}$ depends on the space $\mathcal D_{(t+1), (t)}$ which is the space generated by the true coefficient $\beta$, the noise $\varepsilon$, the initial condition $q_{(0)}$, and the subsequent terms generated from Algorithm~\ref{algo_single_quantile}.
The proof by induction is similar to the proof of Lemma 1(b) in \citet{bayati2011dynamics}. We prove that $\mathcal H_{(1)}$ holds and if $\mathcal B_{(r)}, \mathcal H_{(s)}$ holds for all $r \leq t$ and $s \leq t$, then $\mathcal H_{(t+1)}$ holds.  Let $o_{k, (t)}(1)$ denote a vector in $\mathbbm R^{t}$ for the $k$th estimator such that all of its entries converge to 0 almost surely for $p \to \infty$.

Step 2 from \citet{bayati2011dynamics}: $\mathcal H_{(1)}$: We know from Eq.(3.35) in \citet{bayati2011dynamics} that for each $k$ and a Gaussian matrix $\tilde{X}_k$ with the same distribution as the design matrix $X$, see also \citet{bayati2011dynamics} Lemma 2 (1),
$$
h_{k, (1)}|_{\mathcal D_{k, (1),(0)}} \stackrel{d}{=}  (\tilde{X}_{k})^\top m_{k, (0)} + o_{k, (1)}(1) q_{k, (0)}.
$$
Let $a_{k, j} = ([(\tilde{X}_{k})^\top m_{k, (0)}]_j + {o}_{1,{k}}(1) q_{k,(0), j} , \beta_{j})$ and $c_{k, j} = ([(\tilde{X}_{k})^\top m_{k, (0)}]_j, \beta_{j})$ where $k = k_1, k_2$. We first show that for any two $k_1$, $k_2$ $\in \{1, \ldots, K\}$.
\begin{eqnarray}\label{eq:ignore o}
 \lim_{p\to\infty}\frac{1}{p} \sum_{j=1}^p \Big[ \tilde\psi_{\rm c}(a_{k_1, j})\tilde\psi_{\rm c}(a_{k_2, j}) - \tilde\psi_{\rm c}(c_{k_1, j})\tilde\psi_{\rm c}(c_{k_2, j})\Big] = 0.
\end{eqnarray}
Since $\tilde\psi_{\rm c}$ is $\kappa_{\rm c}$ order pseudo-Lipschitz, hence we have
 \begin{eqnarray*}
     |\tilde\psi_{\rm c}(a_{k, j}) - \tilde\psi_{\rm c}(c_{k, j})|
     \leq L \{1 + \max(\| a_{k, j} \|^{\kappa_{\rm c}-1}, \| c_{k, j} \|^{\kappa_{\rm c}-1})\}
     | q^0_{k, j}| {o}_{1,{k}}(1) ;
\\  |\tilde\psi_{\rm c}(a_{k, j})| \leq L^\prime(1 + \|a_{k, j}\|^{\kappa_{\rm c}}), \ \ \ \
  |\tilde\psi_{\rm c}(c_{k, j})| \leq L^{\prime\prime}(1 + \|c_{k, j}\|^{\kappa_{\rm c}}) ;
\end{eqnarray*}
meanwhile, from the proof in $\mathcal H_0$ in Lemma 1 in \citet{bayati2011dynamics}, we have for an arbitrary $\kappa_{\rm c}$ order pseudo-Lipschitz function $\tilde\psi_{\rm c}$
\begin{eqnarray}\label{eq:drop_error_single}
 \lim_{p\to \infty}\frac{1}{p} \sum_{j=1}^p |\tilde\psi_{\rm c}(a_{k, j}) - \tilde\psi_{\rm c}(c_{k, j}) | = 0.
\end{eqnarray}
Notice that
\begin{eqnarray*}
\lefteqn{| \tilde\psi_{\rm c}(a_{k_1, j})\tilde\psi_{\rm c}(a_{k_2, j}) - \tilde\psi_{\rm c}(c_{k_1, j})\tilde\psi_{\rm c}(c_{k_2, j})|} \\
&=& |\tilde\psi_{\rm c}(a_{k_1, j})\tilde\psi_{\rm c}(a_{k_2, j}) -\tilde\psi_{\rm c}(a_{k_2, j})\tilde\psi_{\rm c}(c_{k_1, j}) +\tilde\psi_{\rm c}(a_{k_2, j})\tilde\psi_{\rm c}(c_{k_1, j}) - \tilde\psi_{\rm c}(c_{k_1, j})\tilde\psi_{\rm c}(c_{k_2, j})|\\
&\leq& |\tilde\psi_{\rm c}(a_{k_2, j})||\tilde\psi_{\rm c}(a_{k_1, j}) - \tilde\psi_{\rm c}(c_{k_1, j})| +
    |\tilde\psi_{\rm c}(c_{k_1, j})||\tilde\psi_{\rm c}(a_{k_2, j}) - \tilde\psi_{\rm c}(c_{k_2, j})|.
\end{eqnarray*}
Then we have
\begin{eqnarray}
\lefteqn{ \frac{1}{p} \sum_{j=1}^p | \tilde\psi_{\rm c}(a_{k_1, j})\tilde\psi_{\rm c}(a_{k_2, j}) - \tilde\psi_{\rm c}(c_{k_1, j})\tilde\psi_{\rm c}(c_{k_2, j})| }  \nonumber\\
  &\leq& \frac{1}{p} \sum_{j=1}^p |\tilde\psi_{\rm c}(a_{k_2, j})||\tilde\psi_{\rm c}(a_{k_1, j}) - \tilde\psi_{\rm c}(c_{k_1, j})| +
    |\tilde\psi_{\rm c}(c_{k_1, j})||\tilde\psi_{\rm c}(a_{k_2, j}) - \tilde\psi_{\rm c}(c_{k_2, j})|  \nonumber \\
 &{\leq}& \max_j |\tilde\psi_{\rm c}(a_{k_2, j})|\cdot \frac{1}{p}\sum_{j =1}^p |\tilde\psi_{\rm c}(a_{k_1, j}) - \tilde\psi_{\rm c}(c_{k_1, j})| \nonumber \\
 && +\max_j |\tilde\psi_{\rm c}(c_{k_1, j})| \cdot\frac{1}{p}\sum_{j =1}^p |\tilde\psi_{\rm c}(a_{k_2, j}) - \tilde\psi_{\rm c}(c_{k_2, j})|  \nonumber \\
 &\leq& L^\prime_2\{1 + \max_j(\|a_{k_2, j}\|^{\kappa_{\rm c}})\}
  \frac{1}{p} \sum_{j=1}^p|\tilde\psi_{\rm c}(a_{k_1, j}) - \tilde\psi_{\rm c}(c_{k_1, j})|    \nonumber \\
&& +L^{\prime\prime}_1 \{1 + \max_j(\|c_{k_1, j}\|^{\kappa_{\rm c}})\}
  \frac{1}{p} \sum_{j=1}^p|\tilde\psi_{\rm c}(a_{k_2, j}) - \tilde\psi_{\rm c}(c_{k_2, j})|.  \label{eq:tildepsi-a-c}
\end{eqnarray}
By \eqref{eq:drop_error_single}, for $k=k_1,k_2$, $p^{-1} \sum_{j=1}^p|\tilde\psi_{\rm c}(a_{k, j}) - \tilde\psi_{\rm c}(c_{k, j})|$  tends to 0 as $p\to +\infty$. The remaining two factors are finite almost surely: $[(\tilde{X}^{k})^\top m_{k, (0)}]_j$ is a Gaussian random variable which is finite almost surely; $\beta_{0, j}$ is finite almost surely since its limiting distribution
has bounded moments up to $(2 \kappa - 2)$ by assumption \ref{cond:beta}.
Hence, for any pairs $k_1, k_2 \in \{1, \ldots, K\}$
(\ref{eq:ignore o}) holds.

From here, we consider
$\tilde h_{k, (1)}|_{\mathcal D_{k, (1), (0)}} \stackrel{d}{=} (\tilde{X}_{k})^\top m_{k, (0)}$ of which the components have the same distribution as $\|m_{k, (0)} \| Z_{k}/\sqrt{n}$ for $Z_k\sim N(0,1)$.
Conditioning on $\mathcal D_{k_1, (1), (0)}$ and $\mathcal D_{k_2, (1), (0)}$, we use the strong law of large numbers for triangular arrays in Theorem 3 of \citet{bayati2011dynamics} to obtain that
\begin{eqnarray} \label{eq:SLLN}
\lim_{p\to \infty}\frac{1}{p}\sum_{j = 1}^p \Big\{ \prod_{r=1}^2 \tilde\psi_{\rm c}(\tilde h_{k_r, (1),j}, \beta_{j}) -
E_{(\tilde{X}_{k_1}, \tilde{X}_{k_2})}[\prod_{r=1}^2\tilde\psi_{\rm c}(\tilde h_{k_r,(1),j}, \beta_{j})]\Big\} \stackrel{\rm a.s.}{=} 0.
\end{eqnarray}
We first prove \eqref{eq:SLLN}. For $k_1 \neq k_2$, we show that the condition in Theorem 3 of \citet{bayati2011dynamics} holds.
To simplify the notation, we denote the independent copies of the matrices $\tilde{X}_{k_1}, \tilde{X}_{k_2}$ to be $X_{k_1}, X_{k_2}$.
We take the random variables in the triangular array to be
\begin{equation}\label{eq:triangular array def}
\tilde\psi_{\rm c}(\tilde h_{k_1, (1),j}, \beta_{j})\tilde\psi_{\rm c}(\tilde h_{k_2, (1),j}, \beta_{j}) -
E_{(\tilde{X}_{k_1}, \tilde{X}_{k_2})}[\tilde\psi_{\rm c}(\tilde h_{k_1, (1),j}, \beta_{j})\tilde\psi_{\rm c}(\tilde h_{k_2, (1),j}, \beta_{j})]
\end{equation}
and let $0<\rho<1$ then
\begin{eqnarray*}
\lefteqn{\frac{1}{p} \sum_{j = 1}^p E \Big| \prod_{r=1}^2\tilde\psi_{\rm c}(\tilde h_{k_r,(1),j}, \beta_{j}) -
E_{(\tilde{X}_{k_1}, \tilde{X}_{k_2})}[\prod_{r=1}^2\tilde\psi_{\rm c}(\tilde h_{k_r,(1),j}, \beta_{j})] \Big|^{2 + \rho} } \\
&=& \frac{1}{p}  \sum_{j = 1}^p E_{(X_{k_1}, X_{k_2}, \tilde{X}_{k_1}, \tilde{X}_{k_2})}
\Big[ \Big|\tilde\psi_{\rm c}([X_{k_1}^\top m_{k_1, (0)}]_j, \beta_{j})\tilde\psi_{\rm c}([X_{k_2}^\top m_{k_2, (0)}]_j, \beta_{j})  \\
&&- \tilde\psi_{\rm c}([\tilde{X}_{k_1}^\top m_{k_1, (0)}]_j, \beta_{j})\tilde\psi_{\rm c}([\tilde{X}_{k_2}^\top m_{k_2, (0)}]_j, \beta_{j}) \Big|^{2+\rho} \Big]\\
&=& \frac{1}{p}  \sum_{j = 1}^p E \Big[
\Big|\tilde\psi_{\rm c}([X_{k_1}^\top m_{k_1, (0)}]_j, \beta_{j})\tilde\psi_{\rm c}([X_{k_2}^\top m_{k_2, (0)}]_j, \beta_{j}) \\
&& \qquad - \tilde\psi_{\rm c}([X_{k_1}^\top m_{k_1, (0)}]_j, \beta_{j})\tilde\psi_{\rm c}([\tilde{X}_{k_2}^\top m_{k_2, (0)}]_j, \beta_{j})  \\
&& \qquad + \tilde\psi_{\rm c}([X_{k_1}^\top m_{k_1, (0)}]_j, \beta_{j})\tilde\psi_{\rm c}([\tilde{X}_{k_2}^\top m_{k_2, (0)}]_j, \beta_{j})  \\
&& \qquad -\tilde\psi_{\rm c}([\tilde{X}_{k_1}^\top m_{k_1, (0)}]_j, \beta_{j})\tilde\psi_{\rm c}([\tilde{X}_{k_2}^\top m_{k_2, (0)}]_j, \beta_{j}) \Big|^{2+\rho} \Big]\\
& \leq&  \frac{1}{p}  \sum_{j = 1}^p E \Big[
\Big|\tilde\psi_{\rm c}([X_{k_1}^\top m_{k_1, (0)}]_j, \beta_{j})\tilde\psi_{\rm c}([X_{k_2}^\top m_{k_2, (0)}]_j, \beta_{j})  \\
&& \qquad -\tilde\psi_{\rm c}([X_{k_1}^\top m_{k_1, (0)}]_j, \beta_{j})\tilde\psi_{\rm c}([\tilde{X}_{k_2}^\top m_{k_2, (0)}]_j, \beta_{j}) \Big|^{2+\rho} \\
 &&+\frac{1}{p}  \sum_{j = 1}^p E \Big|
 \tilde\psi_{\rm c}([X_{k_1}^\top m_{k_1, (0)}]_j, \beta_{j})\tilde\psi_{\rm c}([\tilde{X}_{k_2}^\top m_{k_2, (0)}]_j, \beta_{j}) \\
&& \qquad -\tilde\psi_{\rm c}([\tilde{X}_{k_1}^\top m_{k_1, (0)}]_j, \beta_{j})\tilde\psi_{\rm c}([\tilde{X}_{k_2}^\top m_{k_2, (0)}]_j, \beta_{j}) \Big|^{2+\rho} \Big]\\
& \leq &\frac{1}{p}  \sum_{j = 1}^p E \Big[
\Big|\tilde\psi_{\rm c}([X_{k_1}^\top m_{k_1, (0)}]_j, \beta_{j}) \Big|^{2+\rho}
\\ && \times \qquad
\Big| \tilde\psi_{\rm c}([X_{k_2}^\top m_{k_2, (0)}]_j, \beta_{j}) -  \tilde\psi_{\rm c}([\tilde{X}_{k_2}^\top m_{k_2, (0)}]_j, \beta_{j}) \Big|^{2+\rho} \big] \\
 && + \frac{1}{p}  \sum_{j = 1}^p E \Big[
\Big|\tilde\psi_{\rm c}([\tilde{X}_{k_2}^\top m_{k_2, (0)}]_j, \beta_{j}) \Big|^{2+\rho}
\\ && \times \qquad
\Big| \tilde\psi_{\rm c}([X_{k_1}^\top m_{k_1, (0)}]_j, \beta_{j}) -  \tilde\psi_{\rm c}([\tilde{X}_{k_1}^\top m_{k_1, (0)}]_j, \beta_{j}) \Big|^{2+\rho} \big]\\
&\leq& \frac{1}{p}  \sum_{j = 1}^p E \Big[
\Big|L^\prime\Big(1 + |[X_{k_1}^\top m_{k_1, (0)}]_j |^{\kappa_c} + |\beta_{j}|^{\kappa_c}\Big) \Big|^{2+\rho}  \\
&& \times \qquad \Big| \tilde\psi_{\rm c}([X_{k_2}^\top m_{k_2, (0)}]_j, \beta_{j}) -  \tilde\psi_{\rm c}([\tilde{X}_{k_2}^\top m_{k_2, (0)}]_j, \beta_{j}) \Big|^{2+\rho} \big]  \\
 && +\frac{1}{p}  \sum_{j = 1}^p E\Big[
\Big|L^\prime\Big(1 + |[\tilde{X}_{k_2}^\top m_{k_2, (0)}]_j |^{\kappa_c} + |\beta_{j}|^{\kappa_c}\Big) \Big|^{2+\rho}\\
&& \times \qquad \Big| \tilde\psi_{\rm c}([X_{k_1}^\top m_{k_1, (0)}]_j, \beta_{j}) -  \tilde\psi_{\rm c}([\tilde{X}_{k_1}^\top m_{k_1, (0)}]_j, \beta_{j}) \Big|^{2+\rho} \big] \\
&\leq&
\max_{j=1,\ldots,p}E \Big[
\Big|L^\prime\Big(1 + |[X_{k_1}^\top m_{k_1, (0)}]_j |^{\kappa_c} + |\beta_{j}|^{\kappa_c}\Big) \Big|^{2+\rho} \Big]
 \\
&& \times\qquad \frac{1}{p}  \sum_{j = 1}^p  E \Big[\Big| \tilde\psi_{\rm c}([X_{k_2}^\top m_{k_2, (0)}]_j, \beta_{j}) -  \tilde\psi_{\rm c}([(\tilde{X}_{k_2})^\top m_{k_2, (0)}]_j, \beta_{j}) \Big|^{2+\rho} \Big] \\
 && + \max_{j=1,\ldots,p} E \Big[
\Big|L^\prime\Big(1 + |[\tilde{X}_{k_2}^\top m_{k_2, (0)}]_j |^{\kappa_c} + |\beta_{j}|^{\kappa_c}\Big) \Big|^{2+\rho} \Big]
 \\ && \times \qquad \frac{1}{p}  \sum_{j = 1}^p
 E \Big[\Big| \tilde\psi_{\rm c}([X_{k_1}^\top m_{k_1, (0)}]_j, \beta_{j}) -  \tilde\psi_{\rm c}([\tilde{X}_{k_1}^\top m_{k_1, (0)}]_j, \beta_{j}) \Big|^{2+\rho} \Big].
\end{eqnarray*}
For the first term in the last inequality above, we see that the expectation
$E[|L^\prime(1 + |[\tilde{X}_{k_2}^\top m_{k_2, (0)}]_j|^{\kappa_c} + |\beta_{j}|^{\kappa_c}) |^{2+\rho}]$ is bounded by some constant, since the expectation is with respect to the matrices $X_{k_1}, X_{k_2}, \tilde{X}_{k_1}, \tilde{X}_{k_2}$ of which the components are Gaussian distributed with mean 0 and variance $1/n$; the rest terms are bounded by a constant; the moments of Gaussian distributed r.v. are all finite. Let us denote the upper bound of this expectation by $L''$, then the first term of the inequality above is bounded by
$$
L''\frac{1}{p}  \sum_{j = 1}^p
 E \Big[\Big| \tilde\psi_{\rm c}([X_{k_2}^\top m_{k_2, (0)}]_j, \beta_{j}) -  \tilde\psi_{\rm c}([\tilde{X}_{k_2}^\top m_{k_2, (0)}]_j, \beta_{j}) \Big|^{2+\rho}\Big],
$$
 which can be shown to be bounded by $c p^{\rho/2}$ following a similar argument as in Lemma 1(b) in \citet{bayati2011dynamics}. The second term similarly can be shown to be bounded by $c' p^{\rho/2}$. Hence the variable defined in (\ref{eq:triangular array def}) satisfies the condition in Theorem 3 in \citet{bayati2011dynamics}; thus the a.s.~convergence holds.
 \\
 In the special case where $k_1 = k_2$, we show that the square of $\psi_c$ is still pseudo-Lipschitz of order $2\kappa_c \le\kappa$, then the almost sure convergence hold by directly applying the result in Lemma 1 in \citet{bayati2011dynamics}.\\
  To simplify the notation, we use $\psi$ to denote any pseudo-Lipschitz function here. For any pairs $x, y \in \mathbbm R^m$, we have
 \begin{eqnarray*}
\lefteqn{ |\psi^2(x) - \psi^2(y)|} \\
  &\leq& |\psi(x) + \psi(y)||\psi(x) - \psi(y)|
  \leq (|\psi(x)| +| \psi(y)|)|\psi(x) - \psi(y)|\\
  &\leq& L^{'}(1 + \|x\|^\kappa + 1 + \|y\|^\kappa)\cdot
  L(1 + \|x\|^{\kappa-1} + \|y\|^{\kappa - 1})\|x - y\| \\
  &\leq& LL^{''} (1 + \|x\|^\kappa + \|y\|^\kappa)(1 + \|x\|^{\kappa-1} + \|y\|^{\kappa - 1})\|x - y\|\\
  &\leq&  LL^{''} (1 + \|x\| + \|y\|)^{2\kappa -1}\|x - y\| \\
  &\leq& LL^{''} 3^{\kappa - 1}(1 + \|x\|^{2\kappa -1} + \|y\|^{\kappa - 1})\|x - y\|.
 \end{eqnarray*}
 Since $\kappa \geq 1$, $\|x\|, \|y\| \geq 0 $, the last two inequalities are obtained by applying the first and second inequality in Lemma 2, respectively. Hence, the square of any arbitrary pseudo-Lipschitz function of order $\kappa$ is still pseudo-Lipschitz with order $2\kappa$.
This proves \eqref{eq:SLLN}.

Using Lemma~\ref{lemma:lemma_1_coef_converge} for $v = \beta$ and
$$\psi(\beta_{j}) = E_{(\tilde{X}_{k_1}, \tilde{X}_{k_2})} \tilde\psi_{\rm c}(\tilde h_{k_1, (1), j}, \beta_{j})\tilde\psi_{\rm c}(\tilde h_{k_2, (1), j}, \beta_{j}),$$ the following convergence holds
\begin{eqnarray*}
\lefteqn{ \lim_{p\to \infty}\frac{1}{p}\sum_{j = 1}^p
E_{(\tilde{X}_{k_1}, \tilde{X}_{k_2})} \Big[ \tilde\psi_{\rm c}(\tilde h_{k_1, (1), j}, \beta_{j})\tilde\psi_{\rm c}(\tilde h_{k_2, (1), j}, \beta_{j}) \Big] }\\
&\stackrel{\rm a.s.}{=}& E_{B_0}\bigg[E_{
(Z_{k_1, (0)}, Z_{k_2, (0)})}
\Big[
\prod_{r=1}^{2}\tilde\psi_{\rm c}(\|\frac{m_{k_r, (0)}}{\sqrt{n}} \| Z_{k_r, (0)}, B_0)
 \Big] \bigg]\\
&\stackrel{\rm a.s.}{=}& E\Big[
\prod_{r=1}^{2}
\tilde\psi_{\rm c}(\bar\zeta_{k_r, (0)} Z_{k_r, (0)}, B_0) \Big].
\end{eqnarray*}

Step 4 from \citet{bayati2011dynamics}: $\mathcal H_{t + 1}$:
Following the first expression in the proof of Lemma 1(b) in step 4 in \citet{bayati2011dynamics}, for any index $k=1,\ldots,K$
\begin{eqnarray*}
\lefteqn{\tilde\psi_{\rm c}(h_{k,(1), j}, \ldots, h_{k, (t+1), j}, \beta_{j})|_{\mathcal D_{k, (t+1), (t)}}
\stackrel{d}{=}} \\
& \hspace*{-3mm}  \tilde\psi_{\rm c} \Big(h_{k, (1), j}, \ldots, h_{k, (t), j}, \Big[\!\displaystyle \sum_{r= 0}^{t-1}\alpha_r h_{k, (r+1)} \!+\! (\tilde{X}_{k})^\top m_{k, (t)} \!+\!\tilde Q_{k, (t+1)} {o}_{k, (t+1)}(1) \Big]_j\!, \beta_{j}\!\Big).
\end{eqnarray*}
 The columns of $\tilde Q_{k, (t+1)}$ form an orthogonal basis for the column space of $Q_{k, (t+1)} = [q_{k, (0)} \ldots  q_{k, (t)}]$.
Define the matrix  $M_{k,(t)} = [m_{k, (0)} \ldots  m_{k, (t-1)}]$, the vector $(m_{k,(t)})_{\|}=\sum_{r=0}^{t-1}\delta_r m_{k,(r)}$ as the projection of $m_{k,(t)}$ on the column space of $M_{k,(t)}$ and
 the vector $(m_{k,(t)})_\perp = m_{k,(t)}-(m_{k,(t)})_{\|}$.
 Similar to the proof in $\mathcal H_1$, we first show that the error term $\tilde Q_{k, (t+1)} {o}_{k, (t+1)}(1)$ can be dropped.
 Let  $a_{k, j} = $
 $$
\Big(h_{k, (1), j}, \ldots, h_{k, (t), j}, \Big[ \sum_{r= 0}^{t-1}\delta_{r} h_{k, (r+1)} + (\tilde{X}_{k})^\top (m_{k, (t)})_\perp + \tilde Q_{k^{\prime}, (t+1)} {o}_{k, (t+1)}(1) \Big]_j , \beta_{j}\Big)$$
 and
$c_{k, j} = \Big(h_{k, (1), j}, \ldots, h_{k, (t), j}, \Big[ \sum_{r= 0}^{t-1}\delta_{r} h_{k, (r+1)} + (\tilde{X}_{k})^\top (m_{k, (t)})_\perp \Big]_j , \beta_{j}\Big).$
 \\To show that the left hand-side of \eqref{eq:tildepsi-a-c} is finite for the new $a_{k,j}$ and $c_{k,j}$, it suffices to show that both $\max_j(\|a_{k_2, j}\|^{\kappa_{\rm c}})$ and $\max_j(\|c_{k_1, j}\|^{\kappa_{\rm c}})$ are finite almost surely.
By Lemma~\ref{eq:polynomial_ineq}, we obtain the following inequality
\begin{eqnarray*}
 \max_j (\|a_{k_2, j}\|^{\kappa_{\rm c}})
& = & \max_j \Big(C (\sum_{r = 0}^t |h_{k_2, (r+1), j} |^{\kappa_c} + |\beta_{j}|^{\kappa_c}) \Big) \\
 & \leq &  C\Big(\sum_{r = 0}^t \max_j |h_{k_2,(r+1), j}|^{\kappa_c} + \max_j |\beta_{j}|^{\kappa_c} \Big)
\end{eqnarray*}
for some constant $C$.
The finiteness of $\max_j |\beta_{j}|^{\kappa_c}$ has been discussed in $\mathcal H_1$; $\max_j |h_{k_2, (r+1), j}|$ is finite almost surely since Lemma 1 in \citet{bayati2011dynamics} states that for a higher order $l = k - 1$, $\lim\sum_{p \to \infty} p^{-1}\sum_{j = 1}^p (h_{k_2, (t+1), j})^{2l} < \infty$. The almost-sure finiteness of $\max_j|h_{k_2, (r+1), j}|$ follows by a simple contradiction: assume $P(\max_j|h_{k_2,(r+1), j}| = \infty)  = P(|h_{k_2, (r+1), j_{\max}}| = \infty)> 0$, then
\begin{eqnarray*}
\lefteqn{P\Big(\sup_{p^\prime \ge p} \frac{1}{p^\prime}\sum_{j = 1}^{p^\prime}(h_{k_2, (t+1), j})^{2l} < \infty \Big) } \\
&=&\!\!\!\! P\Big(\sup_{p^\prime \ge p} \frac{p^\prime - 1}{p^\prime}\{\frac{1}{p^\prime -1}\!\sum_{j \neq j_{\max}}\!(h_{k_2, (t+1), j})^{2l}\}+  \frac{1}{p^\prime} (h_{k_2, (t+1), j_{\max}})^{2l} <  \infty \Big) < 1.
\end{eqnarray*}
The above equation contradicts the result in Lemma 1(e) in \citet{bayati2011dynamics}.
Follow similar arguments, we have $\max_j(\|c_{k_1, j}\|^{\kappa_{\rm c}})$ finite almost surely.
Now we consider the random variable
\begin{eqnarray*}
\lefteqn{\tilde A_{k, j} = \tilde\psi_{\rm c}\Big(h_{k, (1), j}, \ldots, h_{k, (t), j},}
\\
& \Big[ \sum_{r= 0}^{t-1}\delta_r h_{k, (r+1)} + (\tilde{X}_{k})^\top (m_{k, (t)})_\perp + \tilde Q_{k, (t+1)} {o}_{k, (t+1)}(1) \Big]_j , \beta_{j} \Big).
\end{eqnarray*}
Following arguments as in $\mathcal{H}_1$, it is easy to show that
\begin{eqnarray}
 \lim_{p \to \infty} \frac{1}{p} \sum_{j =1}^p \Big[ \tilde A_{k_1, j} \tilde A_{k_2, j} -
 E_{(\tilde{X}_{k_1}, \tilde{X}_{k_2})}\tilde A_{k_1, j} \tilde A_{k_2, j} \Big] \stackrel{\rm a.s.}{=} 0.
\end{eqnarray}
By Lemma~\ref{lemma:lemma_1_coef_converge} and arguments as in the proof of Lemma 1 (b) in \citet{bayati2011dynamics},
\begin{eqnarray*}
\lefteqn{  \lim_{p \to \infty} \frac{1}{p} \sum_{j =1}^p \tilde\psi_{\rm c}\Big(h_{k_1, (1), j}, \ldots, h_{k_1, (t), j}, \!\big[\!\sum_{r= 0}^{t-1}\!\delta_{k_1, (r)} h_{k_1, (r+1)} \!+\! (\tilde{X}_{k_1})^\top (m_{k_1,(t)})_\perp \big]_j , \beta_{j} \Big) } \\
&& \times\tilde\psi_{\rm c}\Big(h_{k_2, (1), j}, \ldots, h_{k_2, (t), j}, \big[ \sum_{r= 0}^{t-1}\delta_{k_2, (r)} h_{k_2, (r+1)} + (\tilde{X}_{k_2})^\top (m_{k_2,(t)})_\perp \big]_j , \beta_{j} \Big) \\
& \stackrel{\rm a.s.}{=} & E_{B_0}E_{(Z_{k_1,(0)},\ldots,Z_{k_1,(t)},Z_{k_2,(0)}, \ldots, Z_{k_2,(t)})}\\
&&\bigg[\prod_{r=1}^2\tilde\psi_{\rm c}\Big(\bar\zeta_{k_r,(0)}Z_{k_r,(0)}, \ldots,\bar\zeta_{k_r,(t)}Z_{k_r,(t)}, B_0 \Big)
 \bigg]\\
& = & E \bigg[ \prod_{r=1}^2
\tilde\psi_{\rm c}\Big(\bar\zeta_{k_r,(0)}Z_{k_r,(0)}, \ldots,\bar\zeta_{k_r, (t)}Z_{k_r,(t)}, B_0 \Big)
\bigg].
\end{eqnarray*}
\end{proof}

\subsubsection{Proof of Corollary~\ref{cor:cov_Zs}}
\begin{proof}
  The almost sure convergence holds by choosing $\tilde\psi_{\rm c}(y_{(0)},\ldots, y_{(t)}, \beta_j) = \psi_{\rm c}(y_{(t)}, \beta_j) = \Big(\beta_j - y_{(t)}\Big) - \beta_j$ in Lemma~\ref{lemma:converge_cov}.
\end{proof}

\subsubsection{Proof of Theorem~\ref{thm:amse_mod_avrg}}
\begin{proof}
By Lemma~\ref{lemma:converge_cov} and choosing $\tilde\psi_{\rm c}(y_{(0)},\ldots, y_{(t)}, \beta_j) = \psi_{\rm c}(y_{(t)}, \beta_j) = \eta(\beta_j - y_{(t)}; \theta_{(t)}) - \beta_j$ which is a pseudo-Lipschitz function of order $\kappa_c = 1$
the convergence in (\ref{eq:ma_amse_theoretical}) is obtained.
\end{proof}

\subsubsection{Proof of Theorem~\ref{thm:cov_like_estimator}}\label{section:proof of thm:cov_like_estimator}
\begin{proof}
Theorem 2 in \citet{bayati2011dynamics} showed that when assigning $1/p$ point mass to each entry of the vector, $\widetilde{\beta}_{k, j, (t-1)}(p)$ converges weakly to $B_0 +\bar\zeta_{k, (t-1)} Z_{k} $ for $p \to \infty$ where $Z_k \sim N(0,1)$ and $B_0$ has p.d.f. $f_{B_0}$.
When $p$ is large, $\widetilde{\beta}_{k, (t-1)}\mid (B_0 = \beta) \approx N(\beta, \bar\zeta_{k, (t-1)}^2 I_p)$;
the normality comes from $Z_{k} \sim N(0, 1)$. Similar results for the Lasso estimator can be found in \citet{bayati2013estimating} and \citet{donoho2016high}. The normality of $\widetilde\beta_{k, (t-1)}$ ensures that the Stein's unbiased risk estimate is applicable for constructing the AMSE estimator. We choose $\mu$, $x$, $\widehat \mu(x)$ and $g(x)$ in Lemma~\ref{lemma:stein's lemma} to be $\beta$, $\widetilde{\beta}_{k, (t-1)}$, $\eta(\widetilde{\beta}_{k, (t-1)}; \theta_{k, (t-1)})$ and $(\eta(\widetilde{\beta}_{k, (t-1)}; \theta_{k, (t-1)}) - \widetilde{\beta}_{k, (t-1)})$, respectively. Recall that $\eta(\widetilde{\beta}_{k, (t-1)}; \theta_{k, (t-1)})$ refers to applying the soft-thresholding function with parameter $\theta_t$ to each entry of the vector $\widetilde \beta_{(t-1)}$. Then the function $\eta(\cdot; \theta_{k, (t-1)})$ is weakly differentiable with the derivative defined almost everywhere on $\mathbbm R^p$ except at $-\theta_{k, (t-1)}$ and $\theta_{k, (t-1)}$ in each coordinate.
\\
Next, consider any pair $(k_1,k_2)$ with $k_1, k_2\in\{1,\ldots,K\}$.
The conditional normality holds for $\widetilde{\beta}_{k_r, (t-1)}$ ($r=1,2$).
Each component of the sequence $\widetilde{\beta}_{k_r, (t-1)}$  is independent of the remaining entries. Hence, the dependence between $\widetilde{\beta}_{k_1, (t-1)}$ and $\widetilde{\beta}_{k_2, (t-1)}$ comes from the entry-wise dependence of the two variables. In other words, there is only dependence between $\widetilde{\beta}_{k_1, (t-1), j_1}$ and $\widetilde{\beta}_{k_2, (t-1), j_2}$ when
$j_1=j_2$. The covariance between the two sequences is
$$
\bar\zeta_{(k_1, k_2), (t-1)} = \Cov(\widetilde{\beta}_{k_1, (t-1)}, \widetilde{\beta} _{k_2, (t-1)}).
$$
Notice that $\widetilde\beta_{k_1, j, (t-1)}$ and $\widetilde\beta_{k_2, j, (t-1)}$ are jointly Gaussian distributed; further, the univariate function $g: x\to \big(\eta(x; \theta) - x \big)$ satisfies the condition in Lemma~\ref{lemma:stein's lemma} \citep{stein1981estimation}. We apply Lemma~\ref{lemma:stein's lemma} to the jointly Gaussian distributed pairs $\widetilde\beta_{k_1, j, (t-1)}$ and $\widetilde\beta_{k_2, j, (t-1)}$, $j = 1, \ldots, p$ with the univariate function $g$.

Meanwhile, since $\bar\zeta_{\textrm{emp}, (t)}^2 = \bar\zeta_{\textrm{emp}, (t-1)}^2 + o(1)$ by assumption, $\theta_{k_r, (t)} = \alpha \bar\zeta_{k_r, (t)}$ where $\alpha$ is fixed for the different iterations, we obtain
\begin{eqnarray*}
 \lefteqn{\lim_{p\to \infty} \frac{1}{p}\sum_{j = 1}^p |\widehat\beta_{k_r, (t), j} - \widehat\beta_{k_r, (t-1), j} | } \\
 &\stackrel{a.s}{=}&\! E\Big|\eta(B_0 + \bar\zeta_{k_r, (t)} Z_{k_r, (t), j}; \theta_{k_r, (t)}) -
\eta(B_0 + \bar\zeta_{k_r, (t-1)} Z_{k_r, (t-1), j}; \theta_{k_r, (t-1)})\Big| \\
 &=&\! E\Big|\eta(B_0 + \bar\zeta_{k_r, (t)} Z_{k_r, (t), j}; \theta_{k_r, (t)}) -
\eta(B_0 + \bar\zeta_{k_r, (t)} Z_{k_r, (t-1), j}; \theta_{k_r, (t)}) + o(1)\Big| \\
& = & \! 0.
\end{eqnarray*}
The almost sure convergence holds by Lemma 1(b) \citep{bayati2011dynamics}. The next equality holds by $\bar\zeta_{(t)}^2 = \bar\zeta_{(t-1)}^2 + o(1)$ and the definition of $\theta_{k_r,(t)}$. The last equality holds because both $Z_{k_r, (t-1), j}$ and $Z_{k_r, (t), j}$ are standard Gaussian distributed.
Thus, $\widehat\beta_{k_r, (t), j} - \widehat\beta_{k_r, (t-1), j} \vert (B_0 = \beta_j)$ converges to 0 almost surely.
Further, by (\ref{eq:pseudo_data and eq:esti}), $\widehat{\beta}_{k_r,(t-1),j} - \widetilde{\beta}_{k_r,(t-1), j} \stackrel{d}{=} {\bar\zeta}_{k_r, (t-1)} Z_{k,j}$ where $Z_{k_r,j}\sim N(0,1)$.
Then, $\widehat{\beta}_{k_r, (t), j} - \widetilde{\beta}_{k_r,(t-1), j} = (\widehat\beta_{k_r, (t), j} - \widehat\beta_{k_r, (t-1), j}) + (\widehat{\beta}_{k_r,(t-1),j} - \widetilde{\beta}_{k_r,(t-1), j}) \stackrel{d}{=} {\bar\zeta}_{k_r, (t-1)} Z_{k,j}$, where $\widehat{\beta}_{k_r, (t), j} = \eta(\widetilde{\beta}_{k_r, (t - 1), j}; \theta_{k_r, (t - 1)})$, by Slutsky's theorem.
Next, Stein's lemma is applied. We denote by $A_j$, conditioning on $(B_0=\beta_j)$ and $\widetilde{\beta}_{k_r,(t - 1), -j} , r=1,2$. It holds that
\begin{eqnarray*}
\lefteqn{ E\Big[\{\eta(\widetilde{\beta}_{k_1, (t-1), j}; \theta_{k_1,(t-1)}) - \widetilde{\beta}_{k_1,(t-1),j}\}(\widetilde{\beta}_{k_2,(t-1),j} - \beta_{j})\vert A_j \Big] } \\
 &=& \Cov\Big(\eta(\widetilde{\beta}_{k_1,(t-1),j}; \theta_{k_1,(t-1)}) - \widetilde{\beta}_{k_1,(t-1),j} , \widetilde\beta_{k_2,(t-1),j} \vert A_j \Big)\\
 &=& \Cov(\widetilde{\beta}_{k_1,(t-1),j}, \widetilde\beta_{k_2,(t-1),j}\vert A_j )
E[\partial_1 \eta(\widetilde \beta_{k_1,(t-1),j}; \theta_{k_1, (t-1)}) - 1\vert A_j].
\end{eqnarray*}
Below we condition everywhere on $B$ which denotes the event that $B_{0,j}=\beta_j$ for $j=1,\ldots,p$ where $B_{0,j}$ are independent copies of $B_0$.
Taking expectation w.r.t.~$\widetilde{\beta}_{k_r, (t-1), -j}$, we obtain for the whole vector,
\begin{eqnarray*}\label{eq:stein_cal_1}
 \lefteqn{
 E\Big[\{\eta(\widetilde{\beta}_{k_1, (t-1)}; \theta_{k_1, (t-1)}) - \widetilde{\beta}_{k_1, (t-1)}\}(\widetilde{\beta}_{k_2, (t-1)} - \beta) \vert B\Big] } \\
 & = &  {\bar\zeta}_{(k_1, k_2), (t-1)}E[\partial_1 \eta(\widetilde{\beta}_{k_1, (t-1)}; \theta_{k_1, (t-1)}) - \mathbf{1}_{p}\vert B]
\\ \label{eq:stein_cal_2}
\lefteqn{
 E\Big[\{\eta(\widetilde{\beta}_{k_2, (t-1)}; \theta_{k_2, (t-1)}) - \widetilde\beta_{k_2, (t-1)}\} (\widetilde{\beta}_{k_1, (t-1)} - \beta) \vert B\Big] } \\
 & = & {\bar\zeta}_{(k_1, k_2), (t-1)}
E[\partial_1 \eta(\widetilde{\beta}_{k_2, (t-1)}; \theta_{k_2, (t-1)}) - \mathbf 1_{p}\vert B].
\end{eqnarray*}
Next, we show the construction of the estimator for $(\Sigma_0)_{(k_1, k_2), (t)}$ at iteration $t$.
The product-sign notation $\prod_{r=1}^2v_r=v_1^\top v_2$.
\begin{eqnarray*}
\lefteqn{E [(\widehat\beta_{k_1, (t)}- \beta)^\top (\widehat\beta_{k_2, (t)} - \beta)\vert B]
= E\big[
\prod_{r=1}^{2}
\{\eta(\widetilde{\beta}_{k_r, (t - 1)}; \theta_{k_r, (t - 1)}) - \beta \}
 \vert B \big] }  \\
 &=& E\big[
 \prod_{r=1}^{2}
 \{\eta(\widetilde{\beta}_{k_r, (t - 1)}; \theta_{k_r, (t-1)}) - \widetilde{\beta}_{k_r, (t - 1)} \} \vert B\big]
  + E\big[\prod_{r=1}^{2}(\widetilde{\beta}_{k_r, (t - 1)} - \beta)\vert B \big] \\
 && + E\big[\{\eta(\widetilde{\beta}_{k_1, (t - 1)}; \theta_{k_1, (t-1)}) -
 \widetilde{\beta}_{k_1, (t - 1)} \}^\top (\widetilde{\beta}_{k_2, (t - 1)} - \beta) \vert B\big]  \\
&& + E\big[(\widetilde{\beta}_{k_1, (t - 1)} - \beta)^\top \{\eta(\widetilde{\beta}_{k_2, (t - 1)}; \theta_{k_2,(t - 1)}) - \widetilde{\beta}_{k_2, (t - 1)} \} \vert B \big] \\
 &=& E\big[
 \prod_{r=1}^{2}
 \{\eta(\widetilde{\beta}_{k_r, (t - 1)}; \theta_{k_r,(t-1)}) -
 \widetilde{\beta}_{k_r, (t - 1)} \}\vert B \big] +
{\bar\zeta}_{(k_1, k_2), (t - 1)} \\
&& +  {\bar\zeta}_{(k_1, k_2), (t - 1)}
\sum_{r=1}^{2}E\big[\partial_1 \eta(\widetilde{\beta}_{k_r, (t - 1)}; \theta_{k_r, (t-1)}) - \mathbf 1_{p} \vert B\big] \\
&=&  - {\bar\zeta}_{(k_1, k_2), (t - 1)} + E\big[
 \prod_{r=1}^{2}
 \{\eta(\widetilde{\beta}_{k_r, (t - 1)}; \theta_{k_r,(t-1)}) -
 \widetilde{\beta}_{k_r, (t - 1)} \}\vert B \big]\\
&& +  {\bar\zeta}_{(k_1, k_2), (t - 1)}
\sum_{r=1}^{2}E\big[\partial_1 \eta(\widetilde{\beta}_{k_r, (t - 1)}; \theta_{k_r, (t-1)})\vert B\big].
\end{eqnarray*}
Replacing the expectations and the covariance $\bar\zeta_{(k_1, k_2), (t -1)}$ with their corresponding empirical versions leads to the unbiased estimator of $(\Sigma_0)_{(k_1, k_2), (t)}$,
\begin{eqnarray*}
\lefteqn{ (\widehat\Sigma_{0, (t)}(p))_{(k_1, k_2)}
 = - {\bar\zeta}_{\textrm{emp}, (k_1, k_2), (t-1)}  }\\
&& + \frac{1}{p} \sum_{j = 1}^p
\prod_{r=1}^{2}
 \{\eta(\widetilde{\beta}_{k_r, (t - 1 ), j} ; \theta_{k_r, (t - 1)}) -
 \widetilde{\beta}_{k_r, (t - 1), j} \}
 \\
 &&
 +\frac{{\bar\zeta}_{\textrm{emp}, (k_1, k_2), (t-1)}}{p}
 \sum_{j = 1}^p \sum_{r=1}^{2}I\{|\widetilde{\beta}_{k_r, (t - 1), j}| \ge \theta_{k_r, (t - 1)})\}.
\end{eqnarray*}
The consistency of the estimator $(\widehat\Sigma_{0, (t)}(p))_{(k_1, k_2)}$ follows since
$$
\lim_{p\to\infty}(\widehat\Sigma_{0, (t)}(p))_{(k_1, k_2)}
 = \lim_{p\to\infty}(\Sigma_{0, (t)}(p))_{(k_1, k_2)}
 = (\Sigma_{(t)})_{(k_1, k_2)}
 $$
 holds with probability one for for all $k_1, k_2 = 1, \ldots, K$. The first equality follows by the unbiasedness of $(\widehat\Sigma_{0, (t)}(p))_{(k_1, k_2)}$ for $(\Sigma_{0, (t)}(p))_{(k_1, k_2)}$, and the second equality holds by Lemma~\ref{thm:amse_mod_avrg}. The proof is completed by realizing the above equality shows almost sure convergence which indicates convergence in probability.
\end{proof}

\subsubsection{Proof of Theorem~\ref{prop:asymp_variance}}
\begin{proof}
 Under the assumption that $n > p$, the model-averaged estimator is unbiased. Hence
 \begin{eqnarray}\label{eq:asymp_var_ma}
 \AMSE(\widehat\beta_{\MA}, \beta) = \lim_{n, p\to \infty} \frac{1}{p}\sum_{j = 1}^p \Var(\widehat \beta_{\rm MA, \it j})
   \stackrel{a.s.}{=} w^\top \Sigma_{(\infty)} w,
 \end{eqnarray}
 where $\Sigma_{(\infty)}$ is a $K\times K$ matrix with $(k_1, k_2)$th component
 \begin{eqnarray}\label{eq:asymp_var_component}
\lefteqn{ (\Sigma)_{(k_1, k_2)} = E\big[\{I(B_0 + \bar\zeta_{k_1} Z_{k_1}) - B_0 \} \{I(B_0 + \bar\zeta_{k_2} Z_{k_2}) - B_0 \}\big]} \nonumber\\
   &=&  \Cov\big(\bar\zeta_{k_1}Z_{k_1}, \bar\zeta_{k_2}Z_{k_2}\big)
  = \Cov(Z_{k_1}, Z_{k_2}) \bar\zeta_{k_1}\bar\zeta_{k_2}.
 \end{eqnarray}

 Combining (\ref{eq:zeta_n>p}) and (\ref{eq:effective_score_n>p}), we obtain that
 \begin{eqnarray}\label{eq:asymp_var_zeta_converge}
\lefteqn{ \bar\zeta_{k} = \delta \{E[\widetilde G (\varepsilon + \bar\zeta_{k}Z_k ; b_k)^2]\}^{1/2} } \nonumber \\
&& = \{E[\widetilde G (\varepsilon + \bar\zeta_{k}Z_k ; b_k)^2]\}^{1/2} \{E[\partial_1 \widetilde G (\varepsilon + \bar\zeta_{k}Z_k ; b_k)]\}^{-1}.
 \end{eqnarray}
 The expressions of the asymptotic variance of the model-averaged estimator in Theorem~\ref{prop:asymp_variance} hold by combining (\ref{eq:asymp_var_ma}), (\ref{eq:asymp_var_component}), and (\ref{eq:asymp_var_zeta_converge}).

\end{proof}

\section*{Acknowledgements}
The authors thank the reviewers for the useful comments which helped improve the paper.
Gerda Claeskens and Jing Zhou acknowledge the support of the Research Foundation Flanders and KU Leuven grant GOA/12/14.
The computational resources and services used in this work were provided by the VSC (Flemish Supercomputer Center), funded by the
Hercules Foundation and the Flemish Government - department EWI.
Jelena Bradic acknowledges the support of the National Science Foundation's Division of Mathematical Sciences grant \#1712481.

\biblist

\end{document}